\documentclass[10pt, leqno]{amsart}

\usepackage{amsmath,amssymb,stmaryrd}
\usepackage[pdftex]{graphicx}
\usepackage[all]{xy}
\usepackage{datetime}

\newcommand{\R}{\mathbb{R}}
\newcommand{\beq}{\begin{equation}}
\newcommand{\ee}{\end{equation}}
\newcommand{\bea}{\begin{eqnarray}}
\newcommand{\eea}{\end{eqnarray}}
\newcommand{\bac}{\begin{array}{c}}
\newcommand{\ea}{\end{array}}
\newcommand{\I}{\mathrm{i}}
\newcommand{\btd}{\nabla}
\newcommand{\btu}{\Delta}

\newcommand{\gm}{\gamma}
\newcommand{\Gm}{\Gamma}
\newcommand{\og}{\omega}
\newcommand{\sg}{\sigma}
\newcommand{\fa}{\forall}

\newcommand{\ie}{{\sl i.e.\/ }}
\newcommand{\cf}{{\sl cf.\/ }}

\newcommand{\dpm}{\displaystyle}
\newcommand{\df}{:=}
\newcommand{\xb}{\boldsymbol{\rm{x}}}
\newcommand{\yb}{\boldsymbol{\rm{y}}}
\newcommand{\Ab}{\boldsymbol{\rm{A}}}
\def\({\left(}
\def\){\right)}
\def\<{\left\langle}
\def\>{\right\rangle}

\newcommand{\pare}[1]{\left(#1\right)}

\def\d{{\partial}}
\newcommand{\eps}{\varepsilon}

\begin{document}
\newtheorem{theorem}{Theorem}[section]
\newtheorem{lemma}{Lemma}[section]
\newtheorem{observation}{Observation}[section]
\newtheorem{definition}{Definition}[section]
\newtheorem{corollary}{Corollary}[section]
\newtheorem{assumption}{Assumption}[section]
\newtheorem{property}{Property}[section]
\newtheorem{remark}{Remark}[section]
\newtheorem{example}{Example}[section]

\numberwithin{equation}{section}


\title{Numerical simulations of  X-rays Free Electron Lasers (XFEL)}
\thanks{Z. Huang was partially supported by the
NSFC Projects No.~11322113, ~91330203, the National Basic Research Program of
China under the grant 2011CB309705. P. Markowich
acknowledges support from the Royal Society through his Wolfson
Research Merit Award and the Award No. KUK-I1-007-43, funded by the
King Abdullah University of Science and Technology (KAUST)}

\author[P. Antonelli]{Paolo Antonelli}
\address[P. Antonelli]
{Centro di Ricerca Matematica Ennio De Giorgi, Scuola Normale Superiore, Piazza dei Cavalieri, 3, 56100 Pisa, Italy}
\email{paolo.antonelli@sns.it}

\author[A. Athanassoulis]{Agissilaos Athanassoulis}
\address[A. Athanassoulis]
{Department of
Applied Mathematics, University of Crete, Heraklion 71409, Greece}
\email{athanassoulis@tem.uoc.gr}

\author[Z.Y. Huang]{Zhongyi Huang$^*$}
\address[Z.Y. Huang]
{Dept. of Mathematical Sciences, Tsinghua University, Beijing
100084, China}
\email{zhuang@math.tsinghua.edu.cn}

\author[P. Markowich]{Peter A. Markowich}
\address[P. Markowich]
{King Abdullah University of Science and Technology (KAUST), MCSE
Division, Thuwal 23955-6900, Saudi Arabia}
\email{P.Markowich@kaust.edu.sa}

\thanks{$^*$ Corresponding author}
\begin{abstract}
We study a nonlinear Schr\"odinger equation which arises as an effective single particle
model in X-ray Free Electron Lasers (XFEL). This equation appears as a first-principles
model for the beam-matter interactions that would take place in an XFEL molecular imaging
experiment in \cite{frat1}. Since XFEL is more powerful by several orders of magnitude
than more conventional lasers, the systematic investigation of many of the standard
assumptions and approximations has attracted increased attention.

In this model the electrons move under a rapidly oscillating electromagnetic
field,  and the convergence of the problem to an effective time-averaged one is examined.
We use an operator splitting pseudo-spectral method to investigate numerically the
behaviour of the model versus its time-averaged version in complex situations, namely the
energy subcritical/mass supercritical case, and in the presence of a periodic lattice.

We find the time averaged model to be an effective approximation, even close to blowup,
for fast enough oscillations of the external field. This work extends previous analytical
results for simpler cases \cite{xfel1}.
\end{abstract}
\maketitle
\tableofcontents

\newpage

\section{Introduction}

In this paper we study a first principles model for beam-matter interaction in X-ray free
electron lasers (XFEL) \cite{Ch,frat1}. Recent developments using XFEL include the
observation of the motion of atoms \cite{Br}, measuring the dynamics of atomic vibrations
\cite{Fri}, biomolecular imaging \cite{Neu} etc. The fundamental model for XFEL is the
following nonlinear Schr\"odinger equation 
\beq \label{eq000} \left\{ \begin{aligned}
&i\hbar\d_t\psi=\frac{1}{2m}(i\hbar\nabla-e\Ab)^2\psi+Ze\frac{1}{|\xb|}\psi+e^2(|\cdot|^{-
1}\ast|\psi|^2)\psi -\lambda|\psi|^\sigma\psi,\\ &\psi(t,\xb)|_{t=0}=\psi_0(\xb),
\end{aligned}\right. \ee 
on $\mathbb{R}^3$. Here $m_e$ is the electron mass, $e$ the
electron charge, $Z$ the atomic number and $\hbar$ is the scaled (by $2\pi$) Planck
constant.  The constant $\lambda>0$ measures the strenght of the local attractive
nonlinearity with power $\sigma>0$, in particular $\sigma=\frac{2}{3}$ and
$\lambda=3e^2\left(\frac{3}{8\pi}\right)^{1/3}$ for the local Hartree-Fock approximation.

A solution $\psi=\psi(t, \xb)$  of this Schr\"odinger equation can be considered as the
wave function of an electron beam, interacting self-consistently through the repulsive
Coulomb (Hartree) force, the attractive local Fock-type approximation with strength
$\lambda$ and exponent $\sigma$ and interacting repulsively with an atomic nucleus,
located at the origin. The vector field ${\Ab}={\Ab}(t)$ represents an external
electromagnetic field, which we shall assume to depend on time $t$ only (not on position
$\xb$).

In the experimental setup of interest \cite{Ch,frat1} the magnetic field is rapidly
oscillating. This effect of oscillating external magnetic field appears often in the
modelling of undulators; here the aim is to justify systematically the use of an
effective, time-averaged magnetic field, in the presence of beam-matter interaction.

 In \cite{xfel1} it was shown that, under appropriate assumptions, the  solution of
 \eqref{eq000} can be approximated, in the asymptotic regime $|\omega| \gg 1$, by the
 solution of an effective, time-averaged Schr\"odinger initial value problem. The present
 work is a continuation of that line of investigation in more general contexts. More
 specifically, with the help of numerical simulations, we are able to tackle more detailed
 questions with respect to the role of the fast magnetic field, the mode of convergence,
 and attack problems with blowup. Moreover, we investigate the effect of a periodic
 lattice on the model.

\subsection{The model}

We scale the equation (1.1) by choosing a characteristic time scale $t_c$ and a characteristic spatial scale $x_c$. Defining the semiclassical parameter $\eps=\frac{\hbar t_c}{2m_ex_c^2}$ we obtain after rescaling

\begin{equation*}
i\eps\d_t\psi=(i\eps\nabla-\Ab)^2\psi+\frac{c}{|x|}\psi+C_1(|\cdot|^{-1}\ast|\psi|^2)\psi-a|\psi|^\sigma\psi
\end{equation*}
where $c=eZ\frac{t_c^2}{2m_ex_c}$, $C_1=e^2\frac{t_c^6}{8m_e^3x_c^4}$,
$a=\lambda\left(\frac{t_c^2}{2m_ex_c^2}\right)^{\sigma+1}$.
Here also $\psi$ and $\Ab$ were scaled appropriately.

A standard transform is used to simplify the original equation, and bring it to a form more amenable to analysis as well as computation. By setting
\beq
\label{eqlenrd}
u(t,\xb) = \psi\(t, \xb + b(t)\) \exp\(\frac{i}{\eps}\dpm\int_0^t{|\Ab(s)|^2ds}\),
\ee
where $b(t)=2\int_0^t\Ab(s)ds$,
one readily checks that $u$ satisfies the initial value problem
\begin{equation}
\label{eq:fast-1}
\left\{\begin{aligned}
&i\eps\partial_tu=-\eps^2\Delta u+\frac{c}{|\xb-b(t)|} u+C_1(|\cdot|^{-1}\ast|u|^2)u
-a|u|^\sigma u\\
&u^\omega(t,\xb)\big|_{t=0}=u_0(\xb),
\end{aligned}\right.
\end{equation}
where $u_0:=\psi_0$.
Now we choose the scales $x_c, t_c$ such that $\eps=1$ (thus relating $x_c$ and $t_c$ by $2m_ex_c^2=\hbar t_c$), since we are mostly interested in the fully quantum mechanical case, and we fix
\begin{equation*}
b(t)=e(t)\sin(2\pi\omega t).
\end{equation*}
where $e:\R\to\R^3$ is a smooth vector field, slowly varying in time.
We consider the large frequency case
\begin{equation*}
|\omega|\gg\frac{{\rm sec}}{t_c}=\frac{\hbar}{2m_e}\left(\frac{{\rm meter}}{x_c}\right)^2\approx 5\times10^{-3}\left(\frac{{\rm meter}}{x_c}\right)^2.
\end{equation*}
and study the asymptotic regime $|\omega|\to\infty$. In particular we shall investigate the convergence of solutions of
\begin{equation}\label{eq:fast}
\left\{\begin{aligned}
&i\partial_tu^\omega=-\Delta u^\omega+V^\omega u^\omega+C_1(|\cdot|^{-1}\ast|u^\omega|^2)u^\omega-a|u^\omega|^\sigma u^\omega\\
&u^\omega(t,\xb)\big|_{t=0}=u_0(\xb),
\end{aligned}\right.
\end{equation}
where the fast oscillatory potential is given by
\begin{equation}\label{eq:fast-p}
V^\omega(t, x):=\frac{c}{|\xb-e(t)\sin(2\pi\omega t)|},
\end{equation}
to solutions of the approximate problem (see \cite{xfel1} and Theorem 1 there)
\begin{equation}\label{eq:avr}
\left\{\begin{aligned}
&i\partial_tu=-\Delta u+\langle V\rangle u+C_1(|\cdot|^{-1}\ast|u|^2)u-a|u|^\sigma u\\
&u(t,\xb)\big|_{t=0}=u_0(\xb),
\end{aligned}\right.
\end{equation}
with $\langle V\rangle$ is the time \emph{averaged} potential given by
\begin{equation}\label{eq:avr-p}
\langle V\rangle(t, \xb):=\int_0^1\frac{c}{|\xb-e(t)\sin(2\pi s)|}ds.
\end{equation}
We remark that \eqref{eq:avr-p} is well defined as the weak limit of $V^\omega$ defined in \eqref{eq:fast-p}
in the space $L^\infty(\R;L^{3-}(\R^3)+L^{3+}(\R^3))$, as $|\omega|\to\infty$ (see \cite{xfel1} for more details).
Passing back and forth between $u$ and $\psi$ (models (\ref{eq:fast}) and (\ref{eq000})) is completely straightforward, by virtue of the transformation (\ref{eqlenrd}).
\begin{remark}
In molecular imaging (spatial characteristic scale of the order of magnitude of thousand nanometers, temporal characteristic scale of some femtoseconds), we compute
$\omega\gg10^{15}$. Therefore it makes sense to investigate the asymptotic regime. We shall demonstrate later on that this is reached already for much smaller values of $\omega$.
\end{remark}

\subsection{Main objectives}

The main result of \cite{xfel1} was the following

\begin{theorem}
Let $0<\sigma<4/3$, $u_0\in L^2(\R^3)$, $u^\omega, u\in C(\R; L^2(\R^3))$ be the unique global solutions of
\eqref{eq:fast}, \eqref{eq:avr}, respectively. Then
for each finite time $0<T<\infty$ and for each admissible Strichartz index pair $(q, r)$, we have
\begin{equation*}
\|u^\omega-u\|_{L^q([0, T]; L^r(\R^3))}\to0\qquad\textrm{as}\;|\omega|\to\infty.
\end{equation*}
\end{theorem}

We use this averaging result as a benchmark (example \ref{ex1}); the energy for the fast problem is observed to converge weakly to the energy for the
averaged problem (and not strongly $L^\infty$ in time). We proceed to investigate mass supercritical problems in example \ref{ex2}, i.e. the exponent $\sigma$ is taken to be larger than $\frac{4}3$. The point is to investigate whether the time-averaged model still gives a good description near blowup. We proceed to an example that highlights the interaction with the (rapidly oscillating) external electromagnetic potential in example \ref{ex:shot}, and with a two-time-scale dependent electromagnetic potential in example \ref{ex:td}. Finally we investigate the interaction with a periodic lattice in example \ref{ex:pl}.

\begin{figure}[htb!]
\centering
\includegraphics[width=0.35\textwidth]{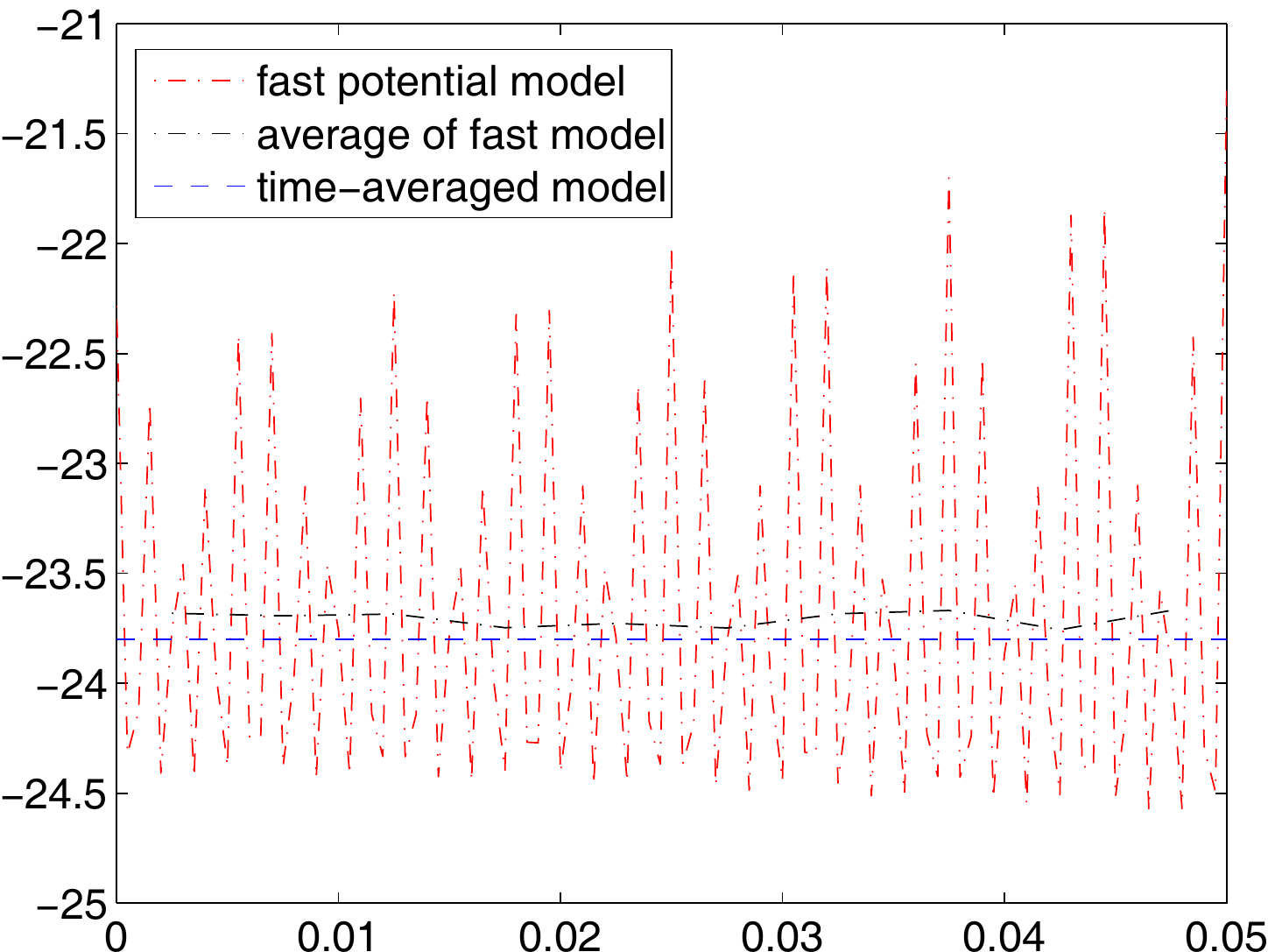}
\includegraphics[width=0.35\textwidth]{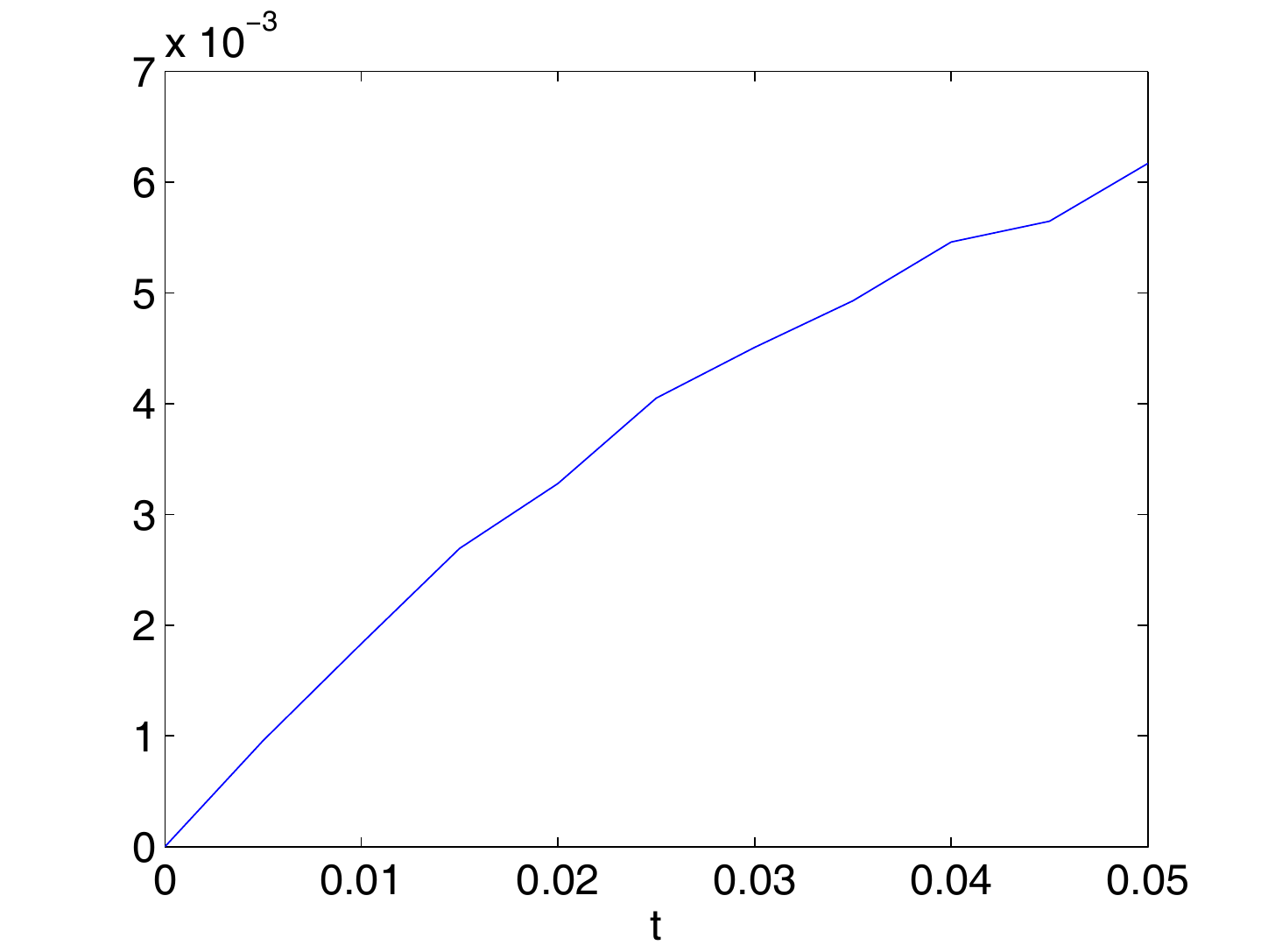}
\caption{
Convergence of the energy for problem (\ref{eq:fast}) to that of (\ref{eq:avr}) is a good way to keep track of the time-averaging,
as well as a diagnostic tool for the computation. Here we show the differences of the energies (left) and the wave functions (right)
between the fast potential model and time-averaged model. We take the time local average of $E^\og(t)$ over $[t,t+0.005]$.
}

\label{Fig222}
\end{figure}

\section{Numerical methods for nonlinear Schr\"odinger equations \eqref{eq:fast} and \eqref{eq:avr}}

As we have shown in \cite{BJM, HJMS, HJMS-2,HJMSZ} the time-splitting pseudo-spectral/Bloch-decomposition based
pseudo-spectral methods works well for (non)linear Schr\"odinger equation (with periodic potential), thus
we shall focus on these two kind of methods in this paper. We have shown that these two methods have spectral convergency in spatial discretization
\cite{BJM, HJMS, HJMS-2,HJMSZ}.
\subsection{Time-splitting spectral algorithm}\label{sec:2.1}
We solve \eqref{eq:fast} or \eqref{eq:avr}
by a classical time-splitting spectral scheme as follows \cite{BJM, HJMSZ}:

\textbf{Step 1.} We solve the equation
\begin{equation}
\label{ts1}
\begin{aligned}
&\, \I \partial _t u_1 = -\, \btu u_1, \\
\end{aligned}
\end{equation}
on a fixed time interval $\Delta t$, relying on the
\emph{pseudo-spectral} method.
Thus we can get the value of $u_1$ at new time
$t^{n+1}$ by
\[u^{n+1}_1=\(e^{-\I\btu t |\xi|^2}\hat{u}_1^{n}(\xi)\)^\vee,\]
here `$\wedge$' and `$\vee$' denote the Fourier and inverse Fourier
transform respectively.

Let $\rho(t^{n+1},\xb)=|u_1(t^{n+1},\xb)|^2$, then we compute the Hartree potential $V_{Hartree}$ by
\[V_{Hartree}(t^{n+1},\xb)=\(\frac{1}{|\xi|^2}\hat{\rho}(t^{n+1}, \xi)\)^\vee.\]

\textbf{Step 2.} We solve the ordinary differential equation
\begin{equation}
\label{ts2}
\begin{aligned}
& \,\I \partial _tu_2  = \left(C_1V_{Hartree}+V(t,\xb)-a|u_2|^\sg\right)u_2, \quad u_2(t^{n},\xb)=u_1(t^{n+1},\xb),
\end{aligned}
\end{equation}
on the same time-interval, here $V(t,\xb)$ is given by \eqref{eq:fast-p} or \eqref{eq:avr-p}.
It is easy to check that in this step the density $\rho(t,\xb)\equiv|u_2|^2(t,\xb)$ and consequently the Hartree potential do not change in time, hence
\[ u_2(t+\btu t,\xb)=u_2(t,\xb)\,\exp\({\dpm-\I \int_{t}^{t+\btu t}
\left(C_1V_{Hartree}(\xb)+V(s,\xb)-a\rho^\frac{\sg}{2}(\xb)\right) ds}\)
\] with $\rho(\xb)=|u_1(t^{n+1},\xb)|^2$, $V_{Hartree}(x)=\(\frac{1}{|\xi|^2}\hat{\rho}(\xi)\)^\vee$.

The key is to calculate the oscillatory integral,
\beq\int_{t}^{t+\btu t} \frac{1}{\left|\xb-b(s)\right|}ds\equiv \int_{t}^{t+\btu t} \frac{1}{\left|\xb-e(s)\sin(2\pi\og s)\right|}ds,\label{2.3}\ee
where $|\og|\gg1$, particularly due to the occurrence of the singularity in the kernel.
Therefore, we use a regularization technique. For example, we could take
\beq\int_{t}^{t+\btu t} \frac{1}{\left|\xb-b(s)\right|}ds\approx \int_{t}^{t+\btu t}\Phi^\eta \(\frac{1}{\left|\xb-b(s)\right|}\)ds,\label{2.4}\ee
where
\[\Phi^\eta \(\frac{1}{\left|\xb-b(s)\right|}\)=\(\frac{2}{\eta}\)^{\frac{3}{2}}
\int_{\R^3}e^{-2\pi\frac{|\xb-\xb'|^2}{\eta}}\frac{1}{\left|\xb'-b(s)\right|}d\xb'.\]

In Section \ref{anal:sect} we will discuss the stability with respect to perturbations of the potential.

The energies for two models are
\beq\label{eq:energy1}
E^\og(t)\df\int_{\R^3}\(|\btd u^\omega|^2+\frac{C_1}{2}|\btd V_{Hartree}^\omega|^2
+cV^\og|u^\omega|^2-a\frac{2}{\sg+2}|u^\omega|^{\sg+2}\)d\xb,
\ee
and
\beq\label{eq:energy2}
E(t)\df\int_{\R^3}\(|\btd u|^2+\frac{C_1}{2}|\btd V_{Hartree}|^2+{c}\langle V\rangle|u|^2-a\frac{2}{\sg+2}|u|^{\sg+2}\)d\xb,
\ee
respectively. The total mass
\beq\label{eq:mass}
M(t)\df\int_{\R^3}|u(t,\xb)|^2d\xb,
\ee
is conserved in both models.

\subsection{Bloch-decomposition based time-splitting spectral algorithm}\label{sec:2.2}
If we consider an external potential $V_{ext}\(\xb\)$ in the equation \eqref{eq000} \cite{Fri,MT},
we have
\[
i\eps\partial_t\psi=(i\eps\nabla-A)^2\psi+V_{ext}(\xb)\psi+c\frac{1}{|\xb|}\psi+C_1(|\cdot|^{-1}\ast|\psi|^2)\psi
-a|\psi|^\sigma\psi,
\]
by the same transformation \eqref{eqlenrd} and by fixing $\eps=1$, we obtain the following equation
\begin{equation}\label{ext-fast}
i\partial_tu^\omega=-\Delta u^\omega+V_{ext}\(\xb-b(t)\)u^\omega+V^\omega u
+C_1(|\cdot|^{-1}\ast|u^\omega|^2)u^\omega-a|u^\omega|^\sigma u^\omega,
\end{equation}
and its time-averaged problem is
\begin{equation}\label{ext-avr}
i\partial_tu=-\Delta u+\langle V_{ext}\rangle u+\langle V\rangle u+C_1(|\cdot|^{-1}\ast|u|^2)u-a|u|^\sigma u,
\end{equation}
with
\[\langle V_{ext}\rangle(t,\xb) = \dpm\int\limits_{-\frac{1}2}^{{\frac{1}2}} V_{ext}\left(\xb-e(t) \sin(2\pi  s)\right) ds.\]

In particular, if we consider a case with an periodic potential $V_\Gm\(\xb\)$ \cite{Fri,MT},
\ie $V_\Gm(\yb)$ is periodic w.r.t to a regular lattice $\Gm$,
\[V_\Gm(\yb+\gm)=V_\Gm(\yb), \quad \fa \gm\in\Gm, \yb\in \R^3,\]
the above two equations \eqref{ext-fast}--\eqref{ext-avr} become
\begin{equation}\label{bloch-fast}
i\partial_tu^\omega=-\Delta u^\omega+V_\Gm\(\xb-b(t)\)u^\omega+V^\omega u
+C_1(|\cdot|^{-1}\ast|u^\omega|^2)u^\omega-a|u^\omega|^\sigma u^\omega,
\end{equation}
and its time-averaged problem is
\begin{equation}\label{bloch-avr}
i\partial_tu=-\Delta u+\langle V_\Gm\rangle\(t,\xb\) u+\langle V\rangle u+C_1(|\cdot|^{-1}\ast|u|^2)u-a|u|^\sigma u,
\end{equation}
with
\[\langle V_{\Gm}\rangle(t,\xb) = \dpm\int\limits_{-\frac{1}2}^{{\frac{1}2}} V_{\Gm}\left(\xb-e(t) \sin(2\pi  s)\right) ds.\]

We solve \eqref{bloch-fast} or \eqref{bloch-avr} by using the \emph{Bloch-decomposition based time-splitting method}, see
\cite{HJMS, HJMS-2}.

\textbf{Step 1.} We solve the equation
\begin{equation}
\label{bd1}
\begin{aligned}
&\, \I \partial _t u_1 = -\, \btu u_1 +V_\Gm\(t,\xb\) u_1, \\
\end{aligned}
\end{equation}
on a fixed time interval $\Delta t$, relying on the
\emph{Bloch-decomposition based pseudo-spectral} method.

Let $\rho(t^{n+1},\xb)=|u_1(t^{n+1},\xb)|^2$, then we obtain the Hartree potential $V_{Hartree}$ by
\[V_{Hartree}\(t^{n+1},\xb\)=\(\frac{1}{|\xi|^2}\hat{\rho}(t^{n+1}, \xi)\)^\vee.\]

\textbf{Step 2.} Then, we solve the ordinary differential equation
\begin{equation}
\label{bd2}
\begin{aligned}
& \,\I \partial _tu_2  = \left(C_1V_{Hartree}+V(t,\xb)-a|u_2|^\sg\right)u_2, \quad u_2\(t^{n},\xb\)=u_1\(t^{n+1},\xb\),
\end{aligned}
\end{equation}
on the same time-interval, here $V(t,\xb)$ is given by \eqref{eq:fast-p} or \eqref{eq:avr-p}.
It is easy to check that  in this step the density $\rho(t,\xb)\equiv|u_2|^2(t,\xb)$ and consequently the Hartree potential do not change in time, hence
\[ u_2(t+\btu t,\xb)=u_2(t,\xb)\,\exp\({\dpm-\I \int_{t}^{t+\btu t}
\(C_1V_{Hartree}(\xb)+V(s,\xb)-a\rho^\frac{\sg}{2}(\xb)\) ds}\)
\] with $\rho(\xb)=\left|u_1(t^{n+1},\xb)\right|^2$, $V_{Hartree}(x)=\(\frac{1}{|\xi|^2}\hat{\rho}(\xi)\)^\vee$.

\begin{remark} Certainly, if there is no lattice potential $V_\Gamma$, the above two algorithms are coincided with each other.
But we must include the lattice potential $V_\Gamma$ when the interaction with laser or other lattice potenital are not eligible.
We will see that the Bloch-decomposition based method is much more efficient in the case with lattice potential $V_\Gamma$
than the traditional pseudo-spectral method (\cf Fig. \ref{Fig012}-\ref{Fig013}), \ie we could use much larger mesh size
than the traditional pseudo-spectral method to get the same accuracy.
\end{remark}
\section{Numerical results}


\subsection{Numerical experiments}

%

\subsubsection{Convergence to the time-averaged model}
\begin{example}\label{ex1} Here we set the semiclassical parameter $\eps=1$, $e(t)=(0,0,1)^T$, $\sigma=\frac{2}{3}$, $a=50$, $C_1=20$, $c=1$.
We compare the time-averaged model with the fast-potential model for different values of $\og$.
We choose appropriate constants $C_1$ and $c$ to make the two potential components in the energy comparable.
\end{example}
\begin{table}[h t]
\begin{center}
\caption{Example \ref{ex1}: errors between two models vs frequency $\og$.}\label{tab:ex1}\vspace{1mm}

\begin{tabular}{c|cccc}\hline
 $\og$  & 5  & 10  & 20 & 40     \\ \hline
$\|u-u^\og\|_{L^2}$ & 1.1E-1 & 5.4E-2 & 2.9E-2 & 1.8E-2  \\ \hline
convergence rate &  &  1.0 & 0.9 & 0.7  \\ \hline
$\|E-E^\og\|_{L^1}$ & 7.1E-1  & 4.5E-1 & 2.9E-1 & 1.9E-1 \\ \hline
convergence rate &  & 0.6 & 0.6 & 0.6  \\ \hline
\end{tabular}
\end{center}
\end{table}

\begin{figure}[!ht]
\centering
\includegraphics[width=40mm]{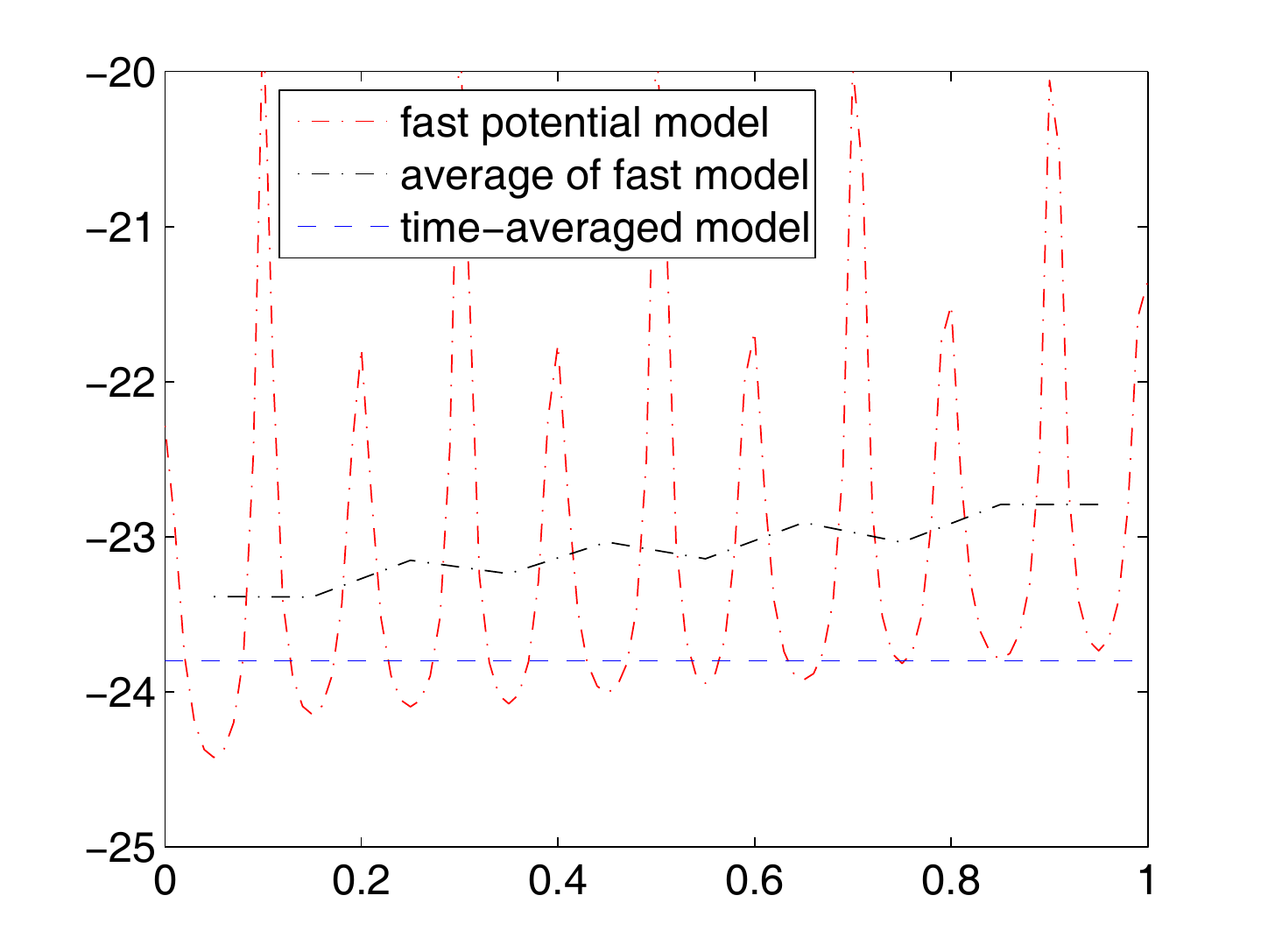}
\includegraphics[width=40mm]{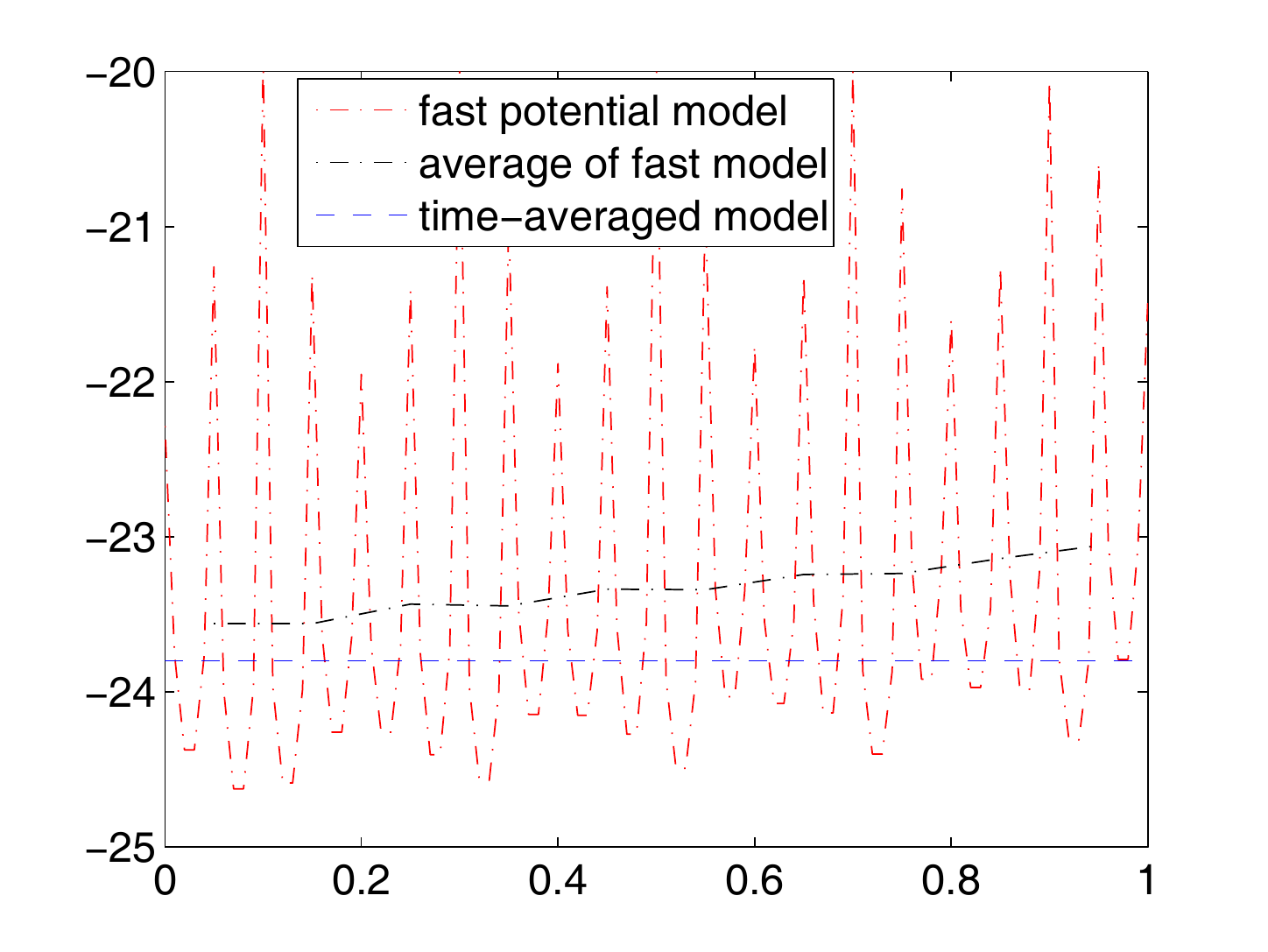}
\includegraphics[width=40mm]{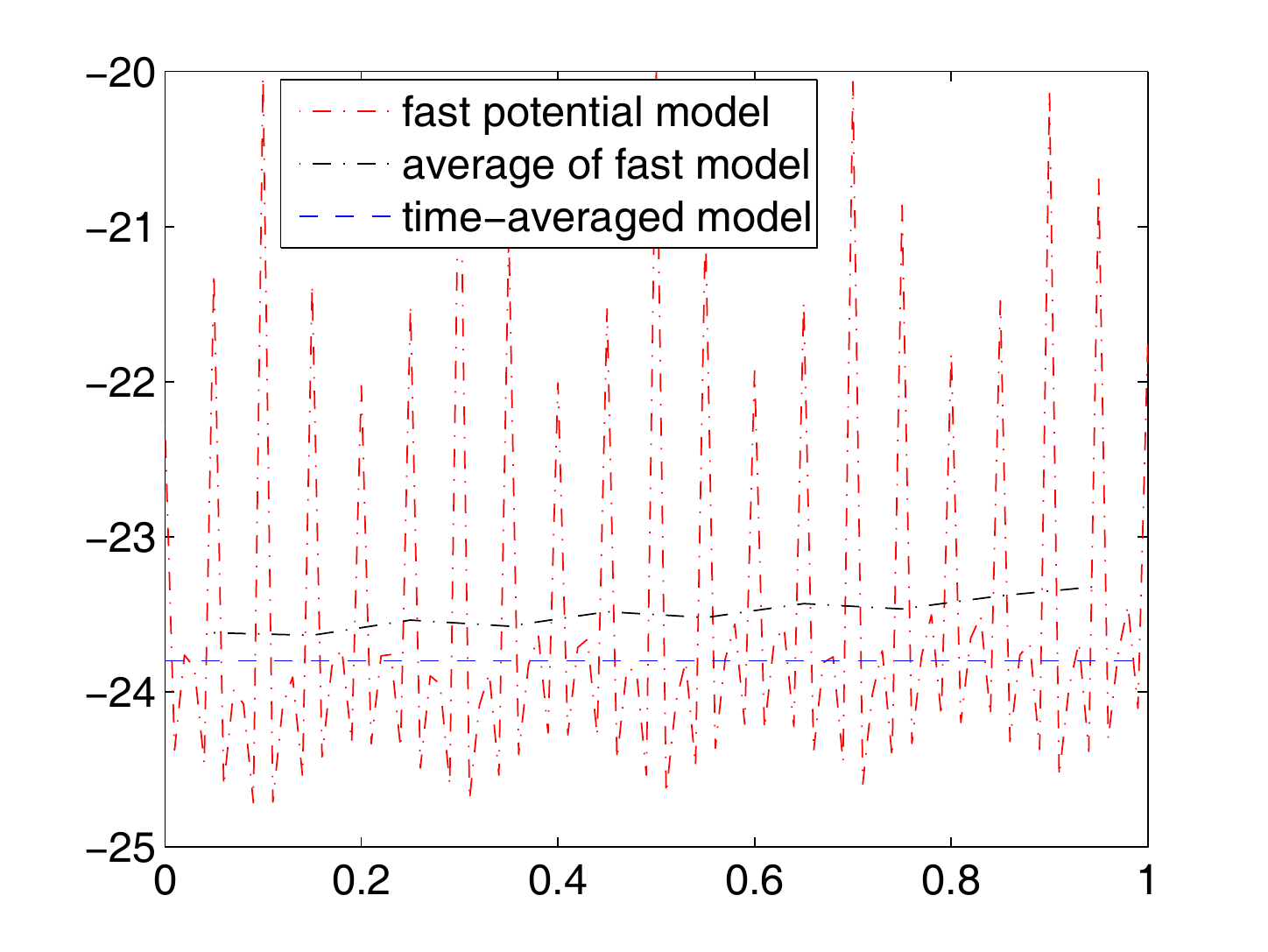}
\includegraphics[width=40mm]{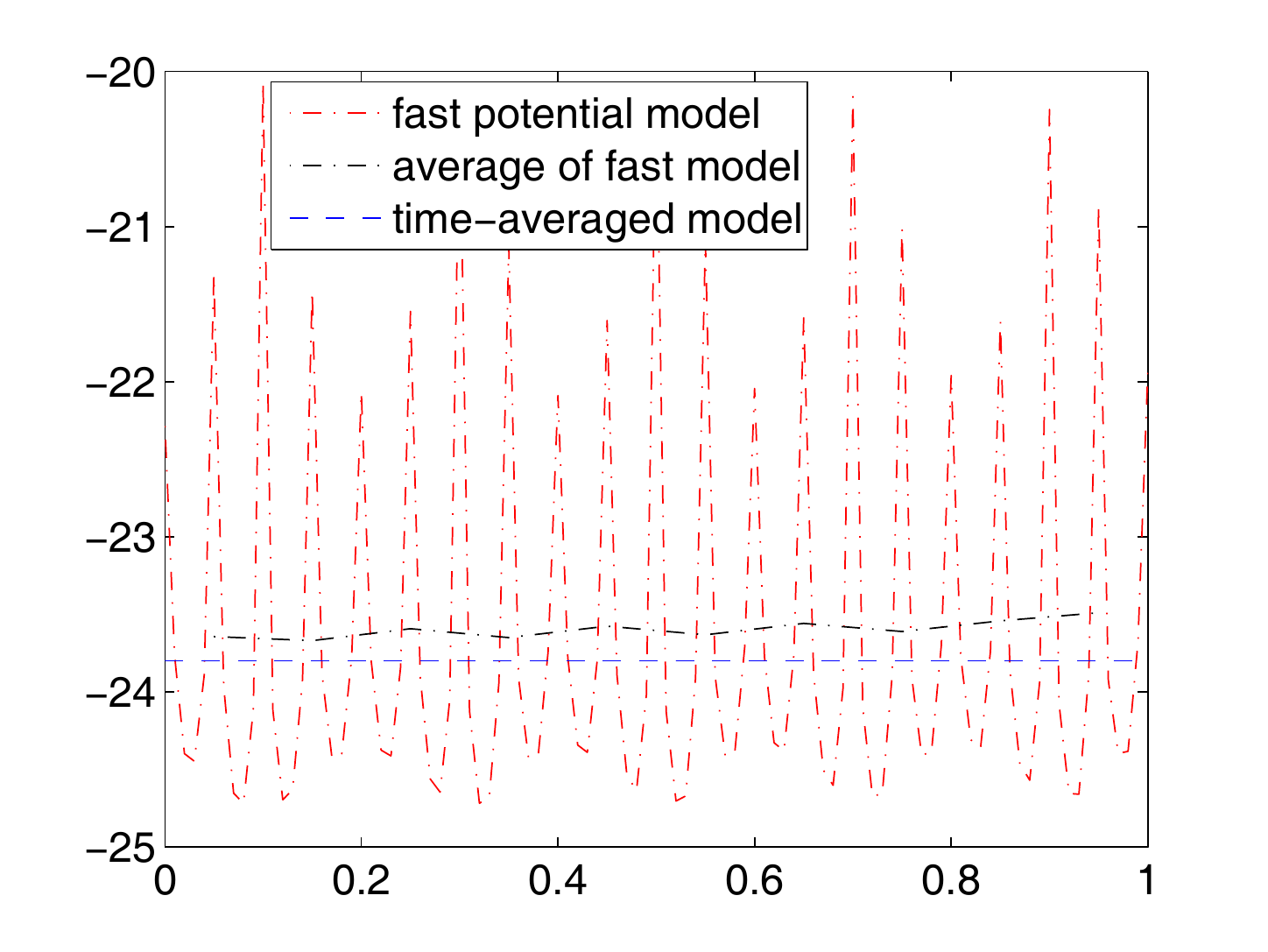}

\hspace{0.3cm}$\og=5$ \hspace{2.95cm} $\og=10$ \hspace{2.9cm} $\og=20$\hspace{2.9cm} $\og=40$
\caption{
Example \ref{ex1}: Comparison of the energies of two models with different $\omega$. Here we take the average of $E^\og(t)$ over $[t,t+0.1]$.}
\label{Fig001}
\end{figure}

\begin{figure}[!ht]
\centering
\includegraphics[width=40mm]{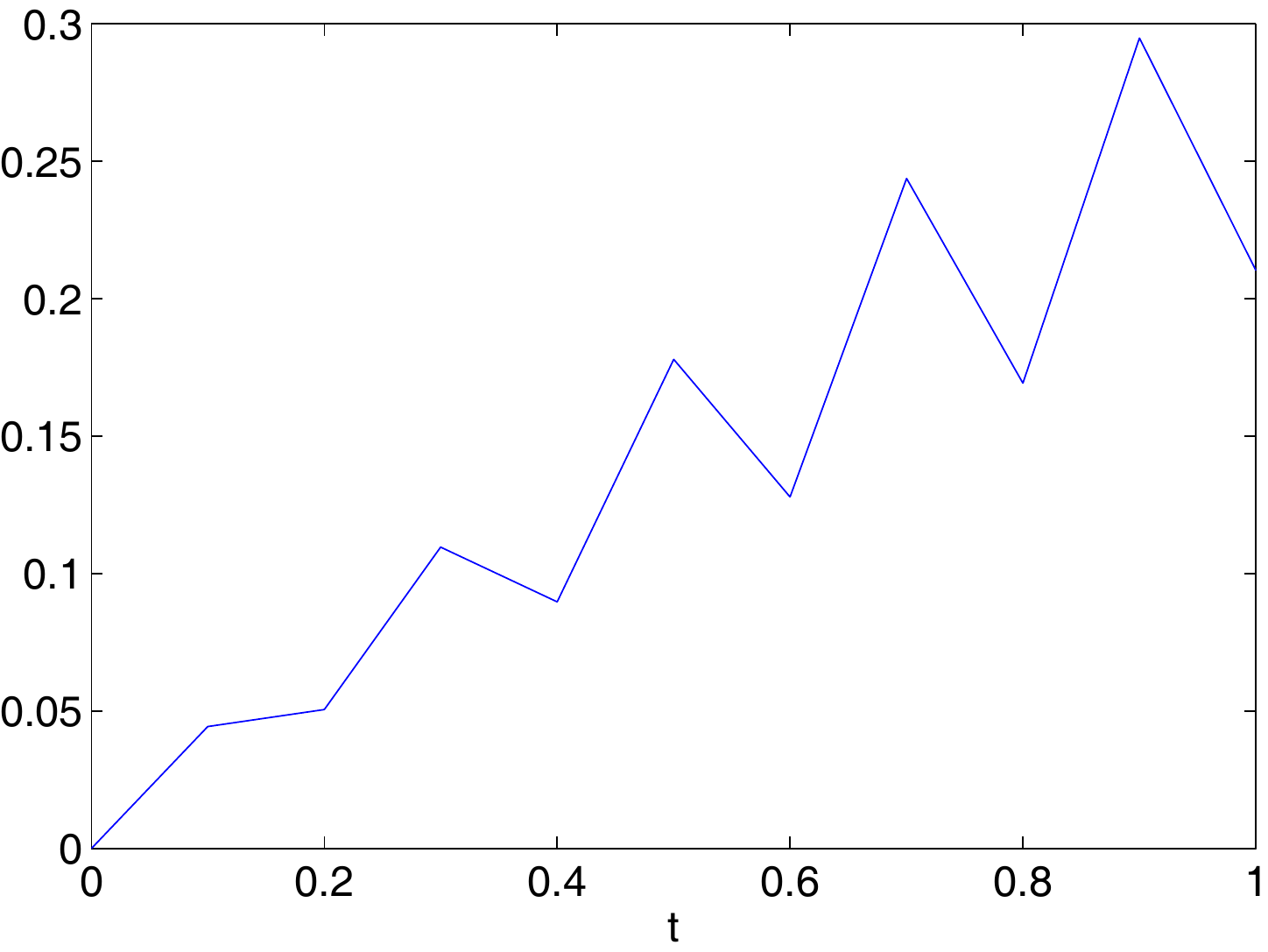}
\includegraphics[width=40mm]{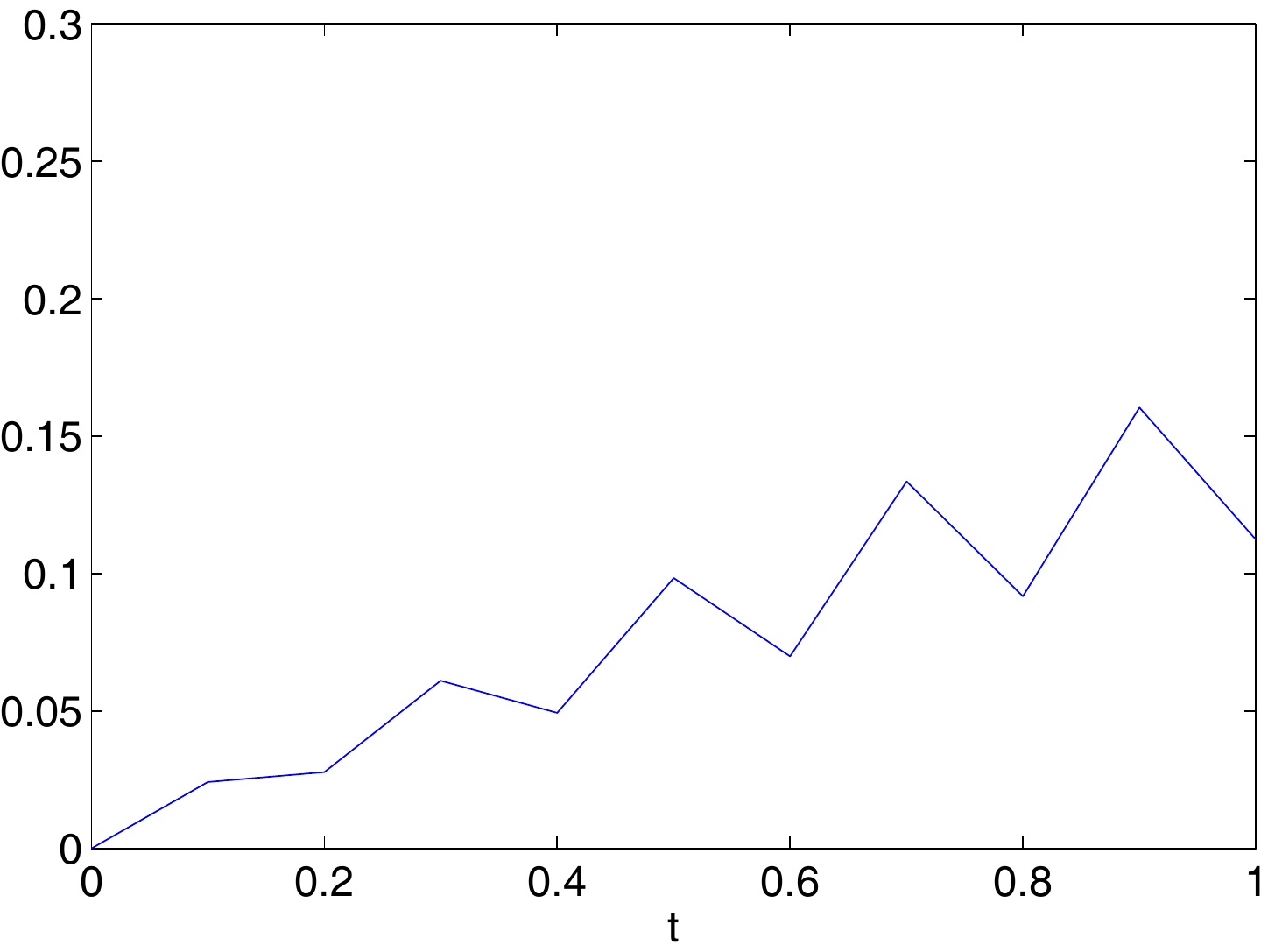}
\includegraphics[width=40mm]{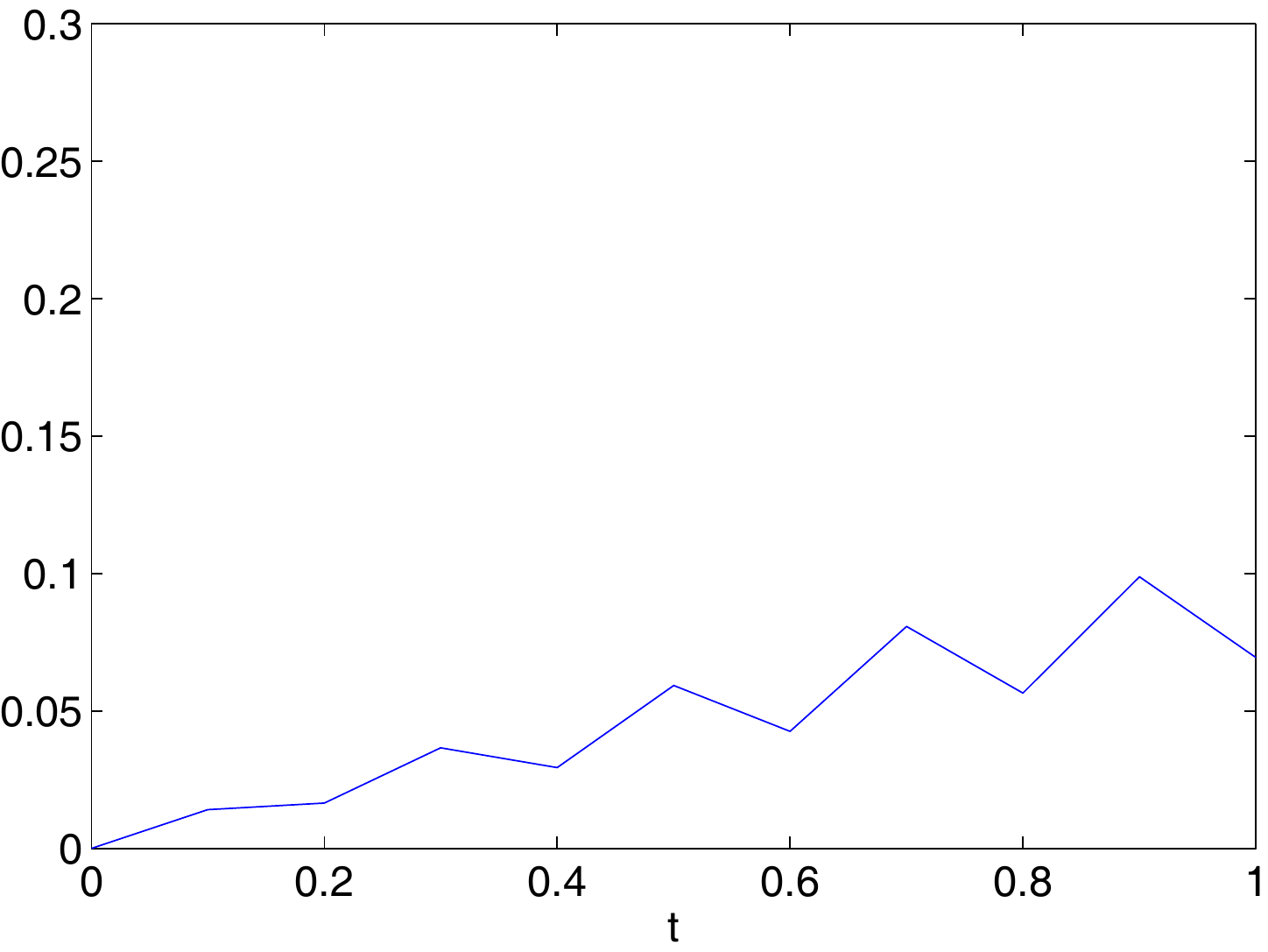}
\includegraphics[width=40mm]{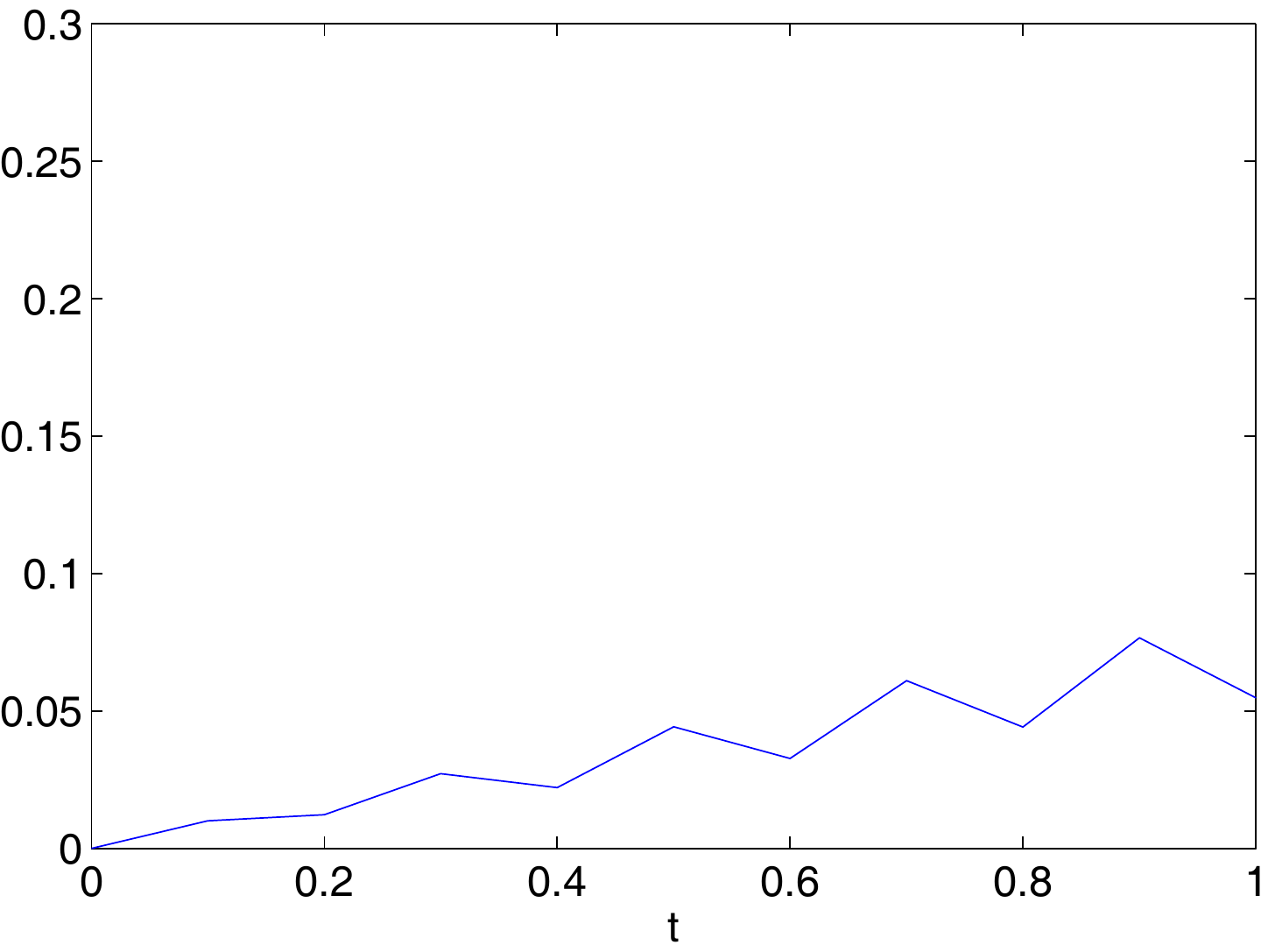}

\hspace{0.3cm}$\og=5$ \hspace{2.95cm} $\og=10$ \hspace{2.9cm} $\og=20$\hspace{2.9cm} $\og=40$
\caption{
Example \ref{ex1}:
The relative distance between the wavefunctions of two models with different $\omega$: $\frac{\|u^\og-u\|_{L^2}}{\|u\|_{L^2}}$.}
\label{Fig001-1}
\end{figure}

\begin{figure}[!ht]
\centering
\includegraphics[width=50mm]{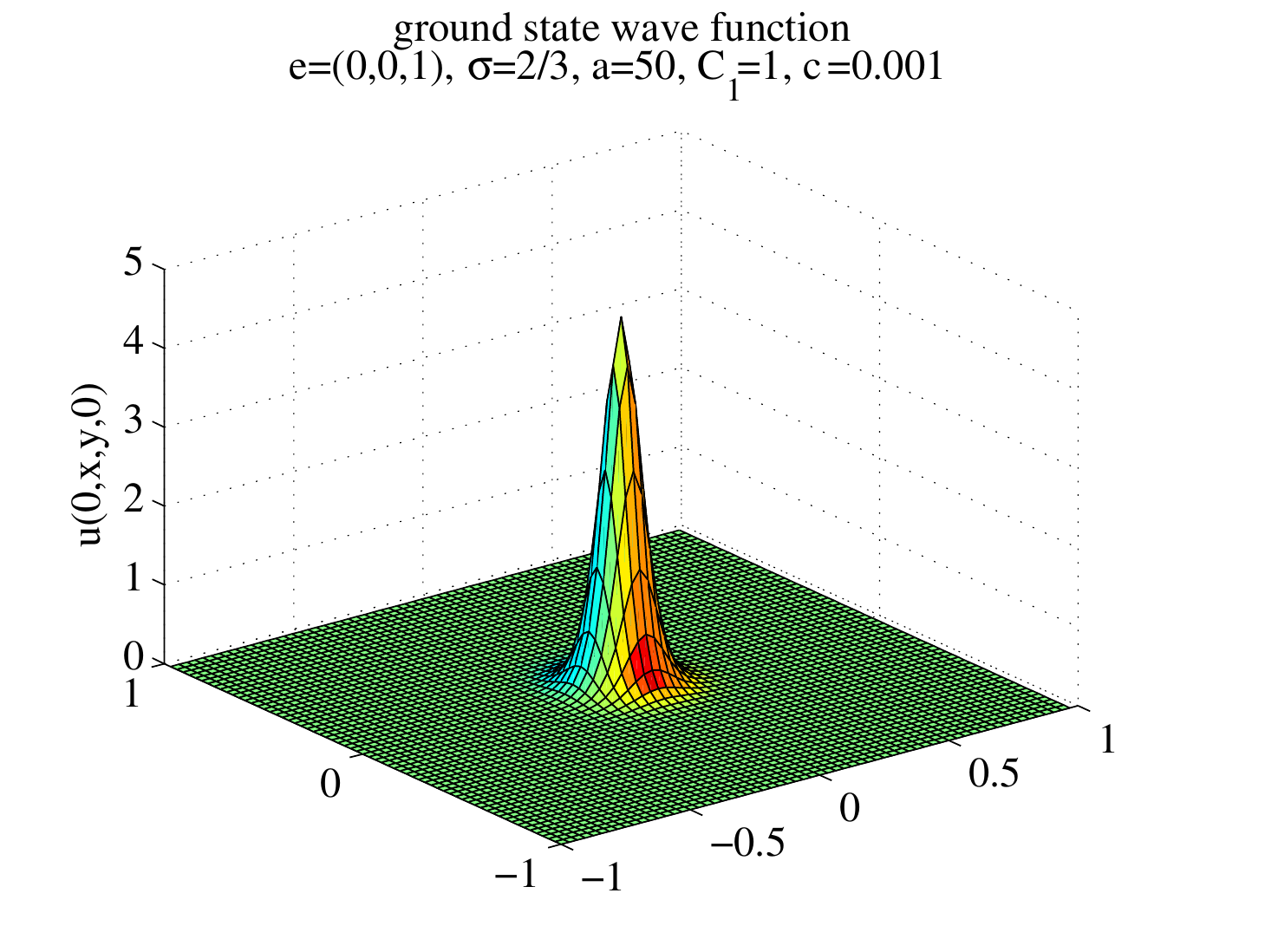}
\includegraphics[width=50mm]{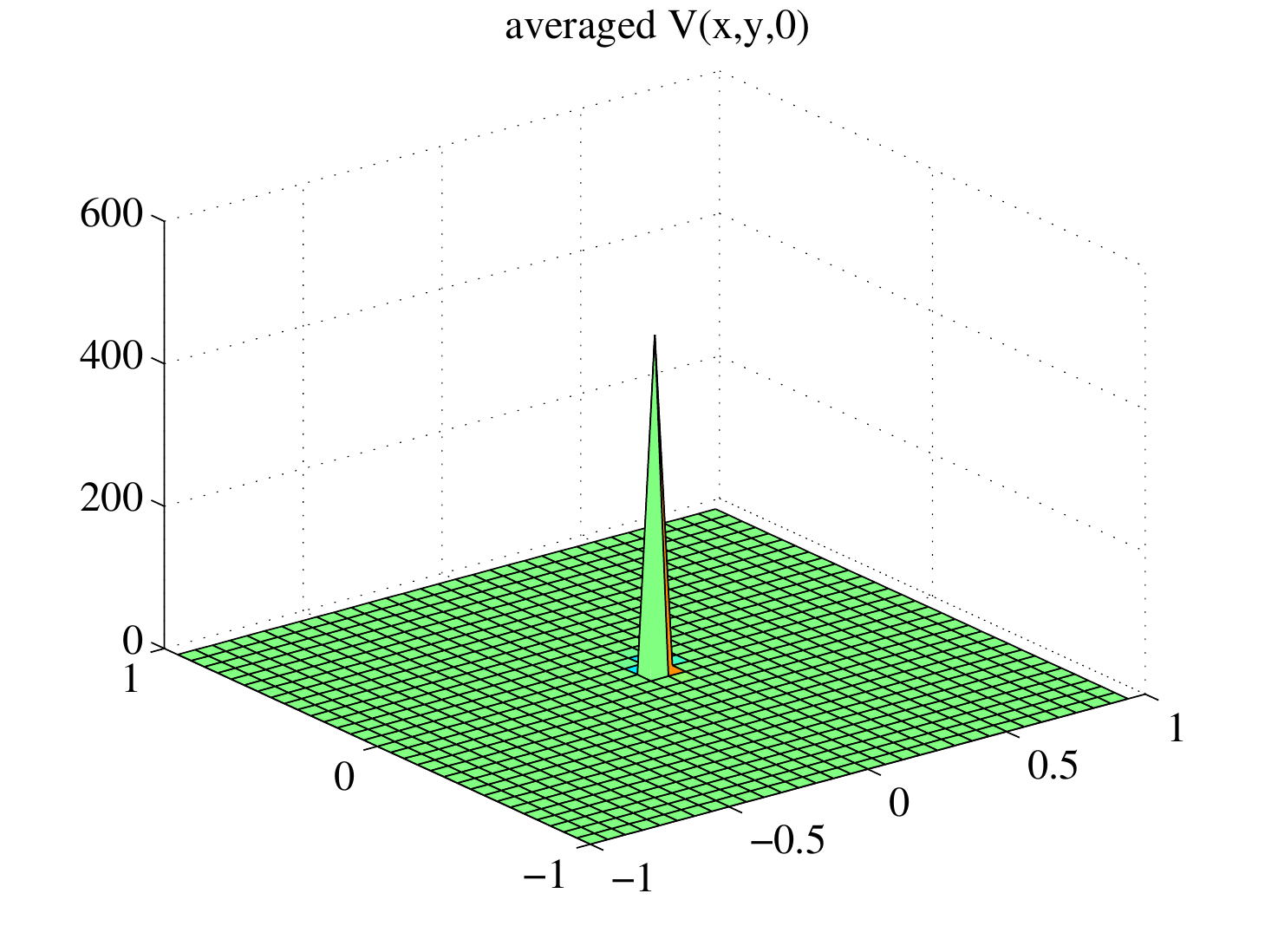}\vspace{-2mm}

Initial datum $u_0(x,y,0)$ \hspace{2cm} $\overline{V}(x,y,0)$ \hspace{1cm}\vspace{2mm}


\caption{Example \ref{ex1}: Graphs of the time-averaged model, the parameters are:
$\eps=1$, $e(t)=(0,0,1)^T$, $\sigma=\frac{2}{3}$, $a=50$, $C_1=1$, $c=0.001$.
Our numerical results show the conservation of the energy $E(t)$ (given by \eqref{eq:energy2})
and mass $M(t)$ (given by \eqref{eq:mass}) in both time-averaged model and fast-potential model.
}
\label{Fig002}
\end{figure}

\begin{figure}[!ht]
\centering
\includegraphics[width=40mm]{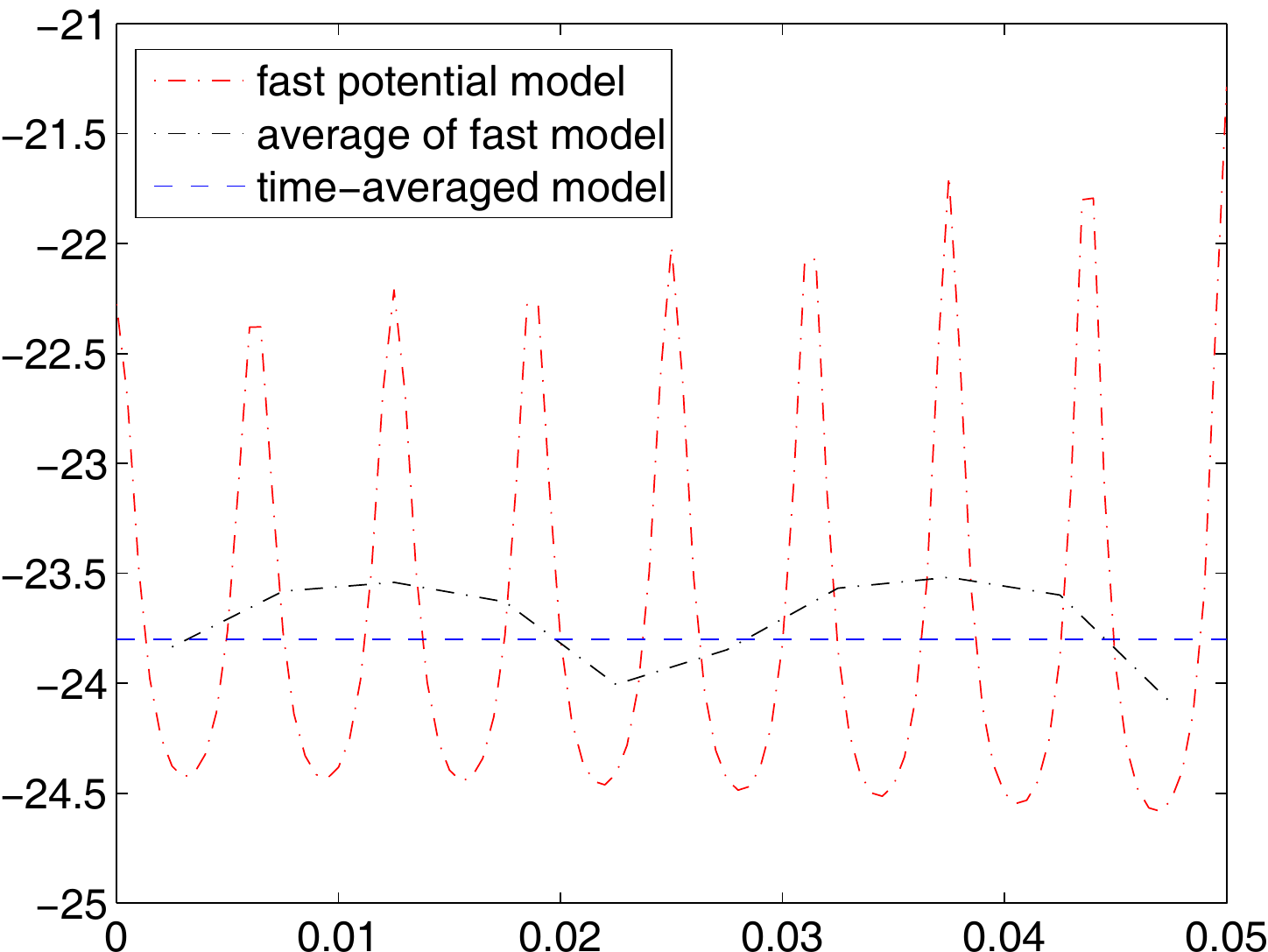}
\includegraphics[width=40mm]{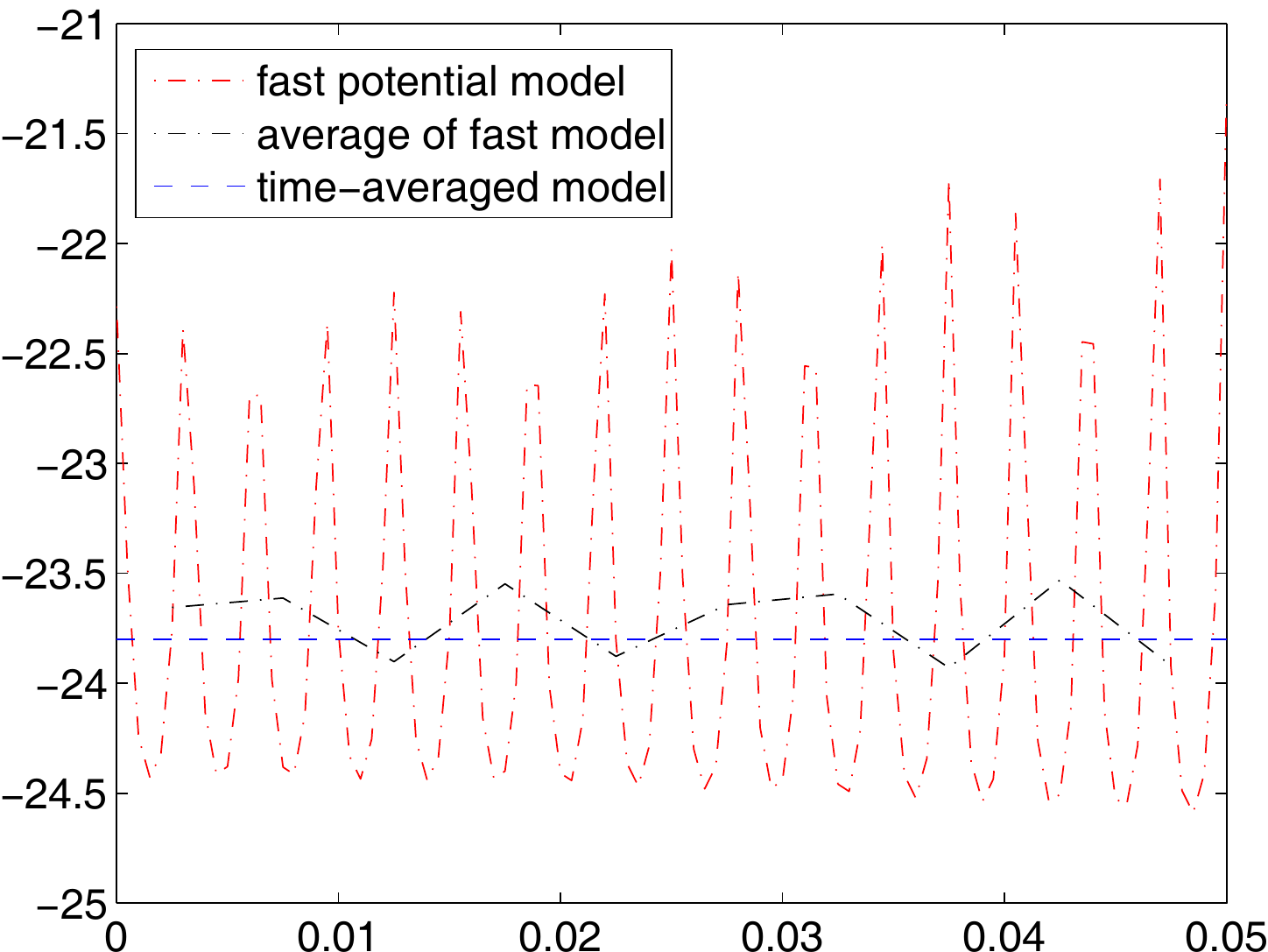}
\includegraphics[width=40mm]{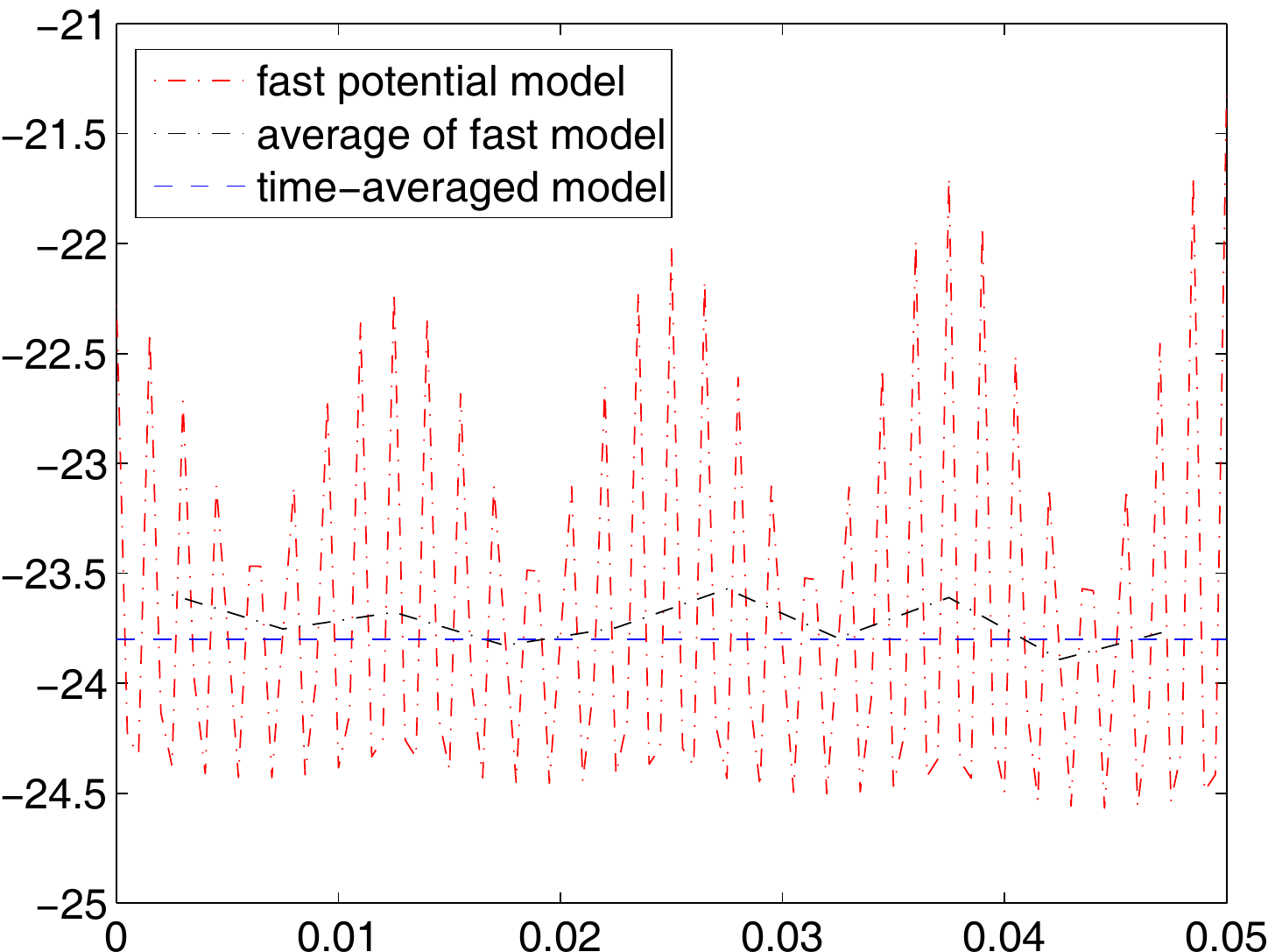}
\includegraphics[width=40mm]{omega640.pdf}

\hspace{0.1cm} $\og=80$ \hspace{2.6cm} $\og=160$ \hspace{2.6cm} $\og=320$ \hspace{2.6cm} $\og=640$
\caption{
Example \ref{ex1}:
The distance between the energies of two models with larger $\omega$ vs. time.
Here we take the average of $E^\og(t)$ over $[t,t+0.005]$.}
\label{Fig001-add}
\end{figure}

\begin{remark} From Figures \ref{Fig001}--\ref{Fig002}, we can see that as $|\og|\to\infty$, the average of the energy $E^\og(t)$ converges to
the energy $E(t)$ of time-averaged model. The corresponding wave-function $u^\og(t,\xb)$ also converges to $u(t,\xb)$ at least in $L^2$-norm.
Figure \ref{Fig002} also shows the conservation of $E(t)$ and $M(t)$.
\end{remark}

\subsubsection{Blow-up tests}
\begin{figure}[!ht]
\centering
\includegraphics[width=50mm]{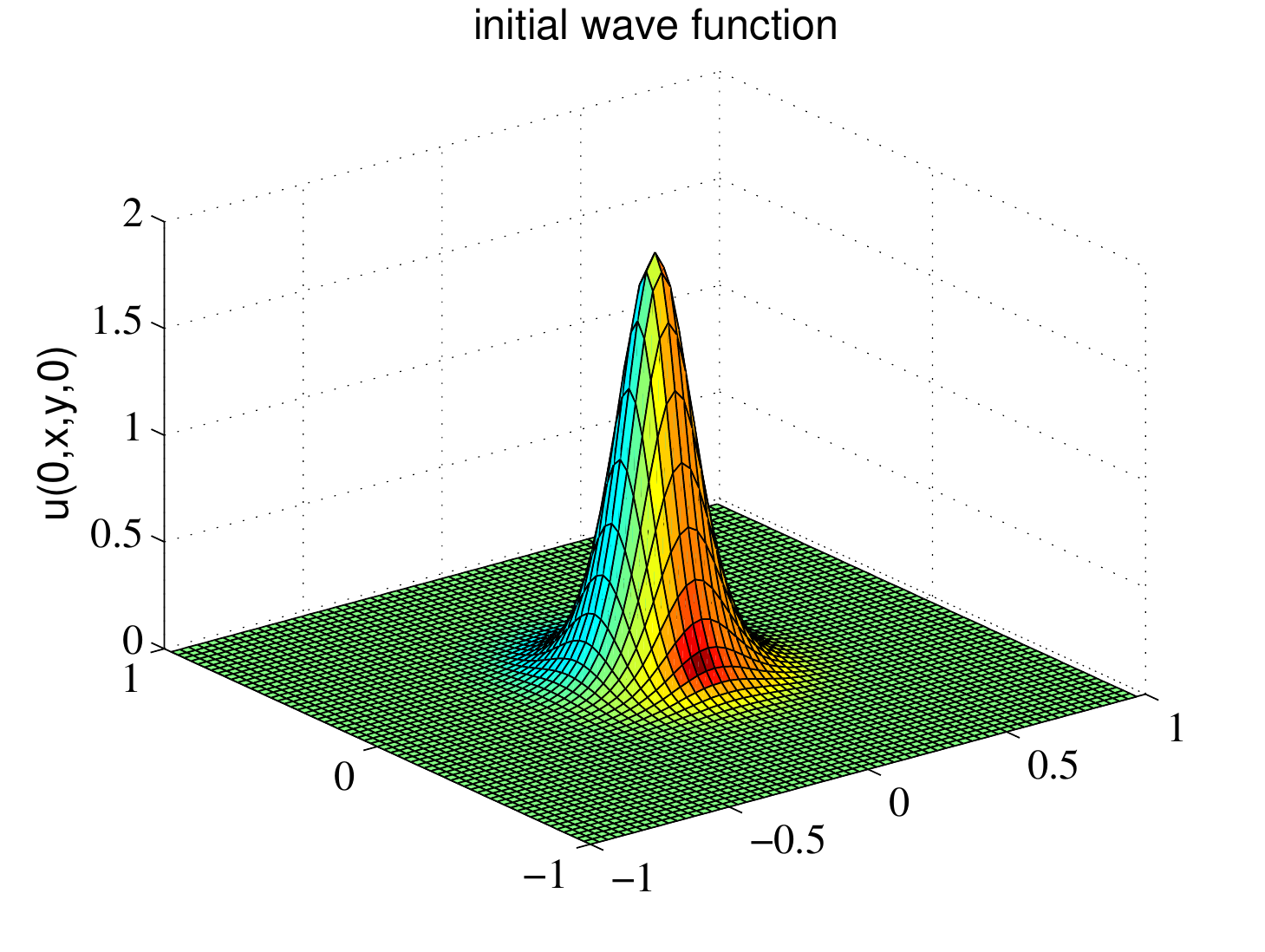}
\includegraphics[width=50mm]{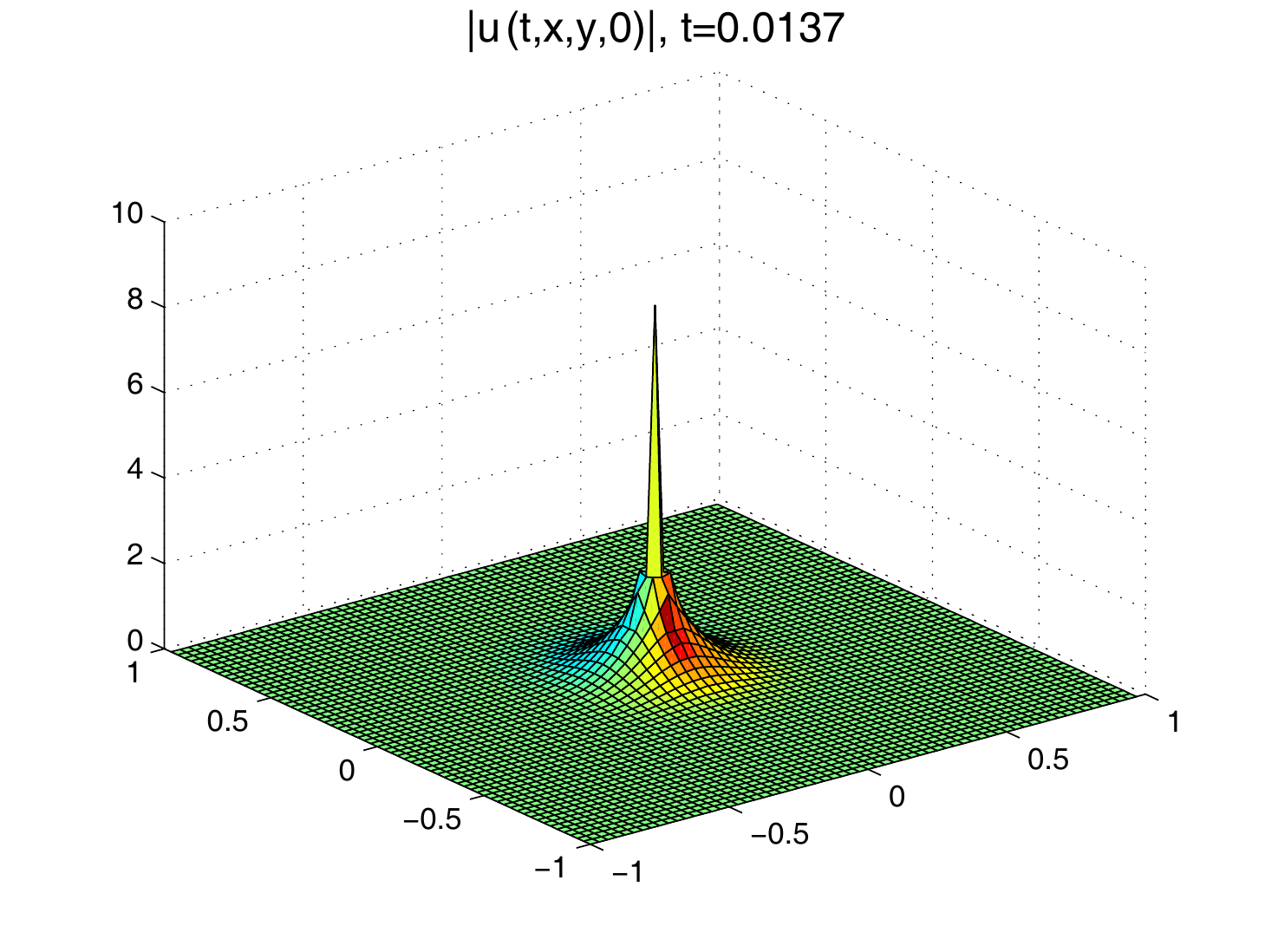}\vspace{-2mm}

Initial datum $u_0(x,y,0)$ \hspace{2cm} $|u(t,x,y,0)|$ \hspace{1cm}\vspace{2mm}

\includegraphics[width=50mm]{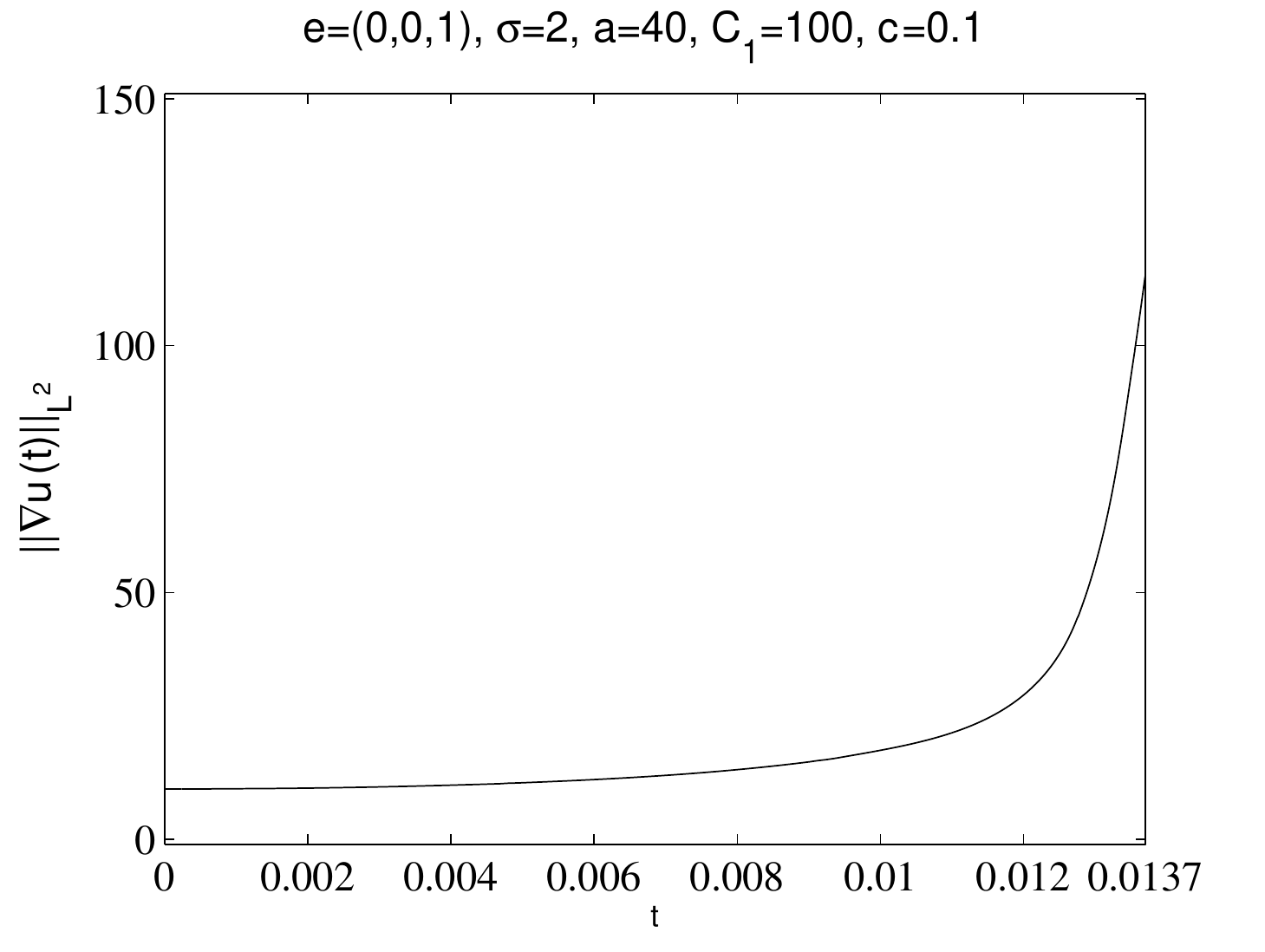}
\includegraphics[width=50mm]{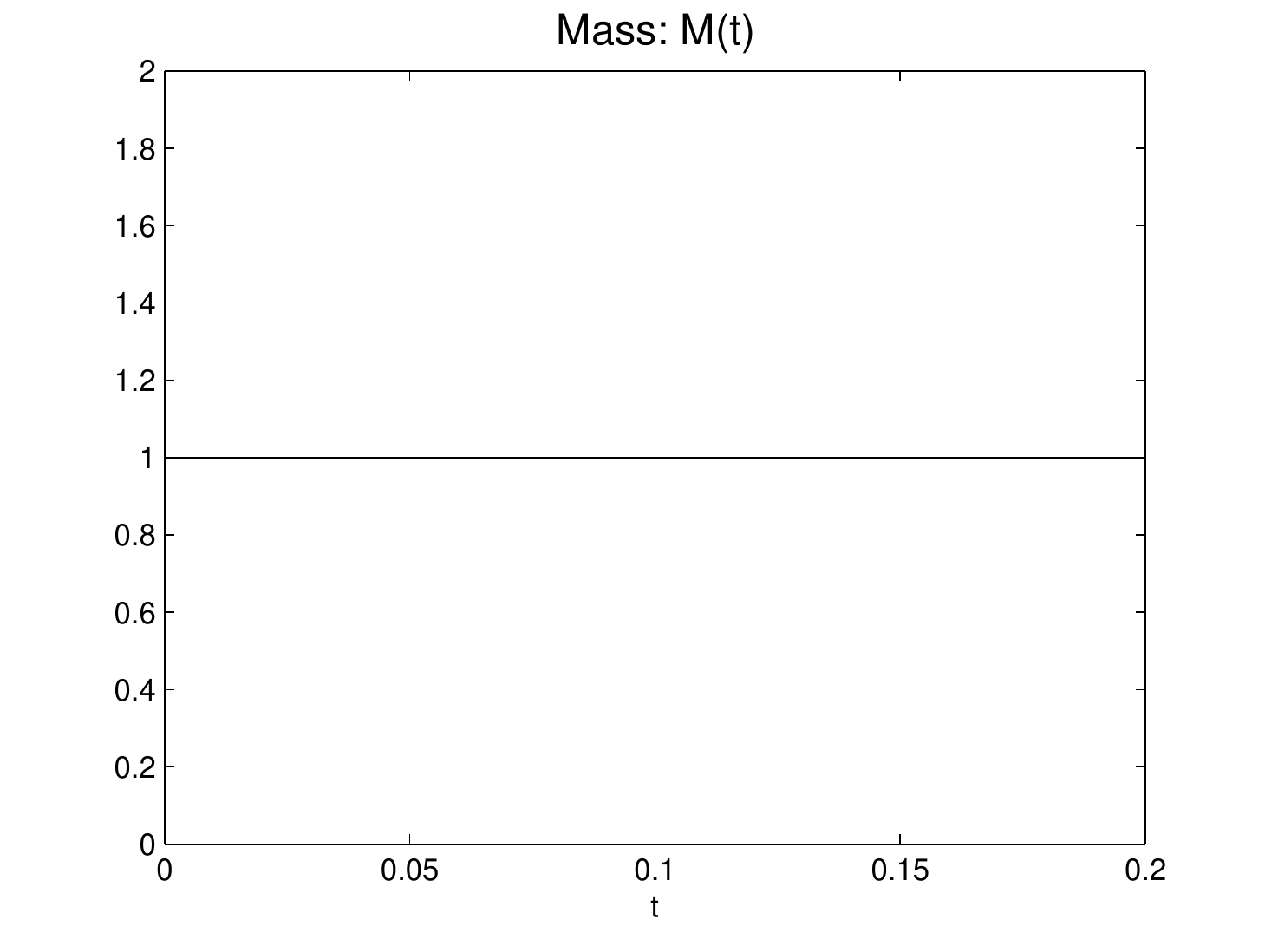}

$\|\btd u\|_{L^2}$ (given by \eqref{eq:H1}) \hspace{1.5cm} $M(t)$ (given by \eqref{eq:mass})
\caption{Example \ref{ex2}: Graphs of the time-averaged model, the parameters are:
$\eps=1$, $e(t)=(0,0,1)^T$, $\sigma=2$, $a=50$, $C_1=100$, $c=0.1$.
Blow-up case.}
\label{Fig003}
\end{figure}

\begin{example} \label{ex2} Here we consider an example with $\eps=1$, $e(t)=(0,0,1)^T$, $a=50$, $C_1=100$, $c=0.1$.
We start the time-averaged model from a ground state for $\sg=2/3$ (\cf plot (a) in Fig. \ref{Fig003}).
Then we consider \eqref{eq:fast} with a mass supercritical power-type nonlinearity, i.e. $\sg>4/3$, we analyse the blow-up of the norm
\beq\label{eq:H1}\|\btd u\|_{L^2}=\left(\int_{\R^3}|\btd u(t,\xb)|^2d\xb\right)^{1/2}.\ee
\end{example}

\begin{figure}[h!t]
\centering
\includegraphics[width=50mm]{initial_40_100_64.pdf}
\includegraphics[width=50mm]{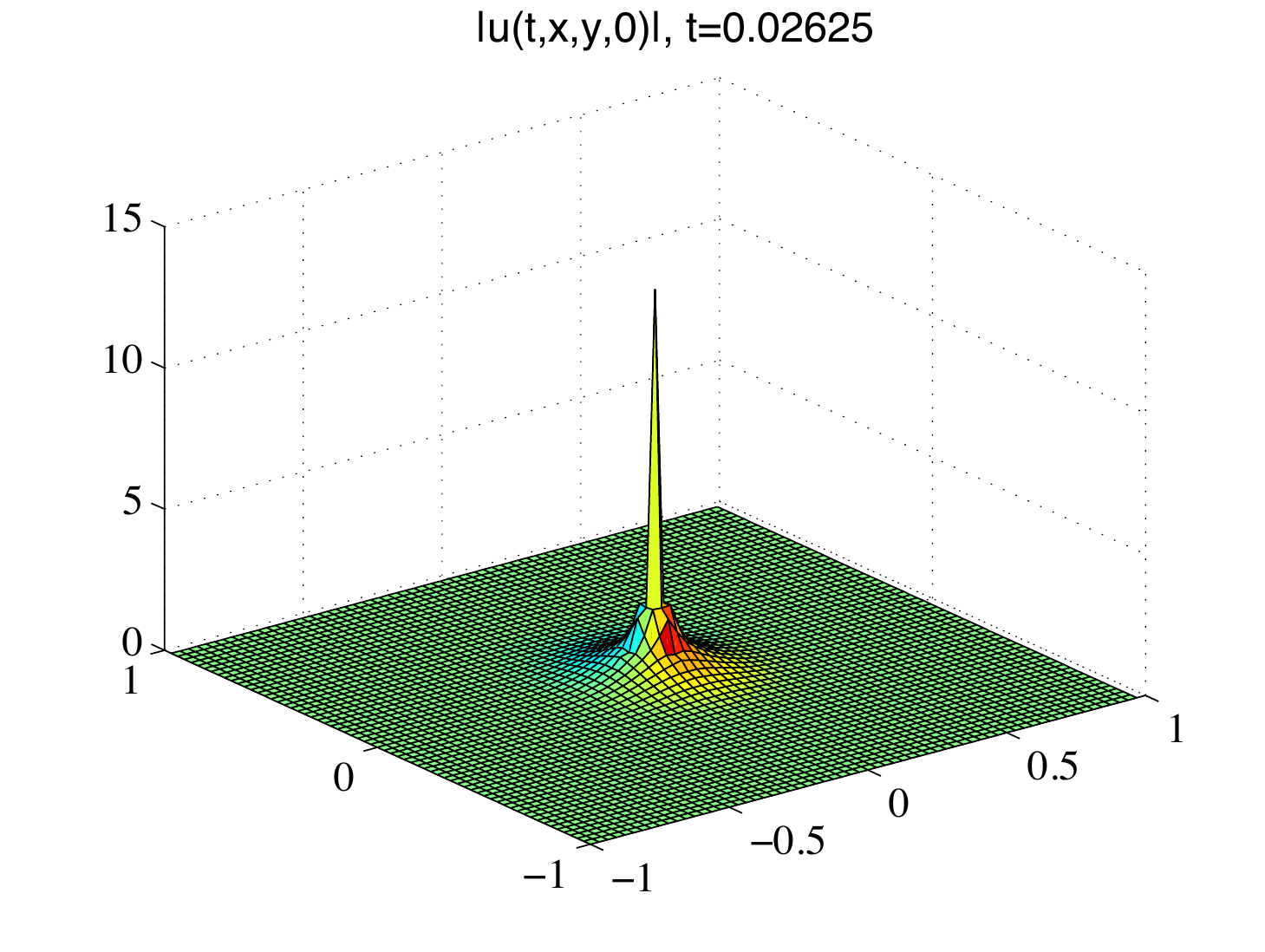}\vspace{-2mm}

Initial datum $u_0(x,y,0)$ \hspace{2cm} $|u(t,x,y,0)|$ \hspace{1cm}\vspace{2mm}

\includegraphics[width=50mm]{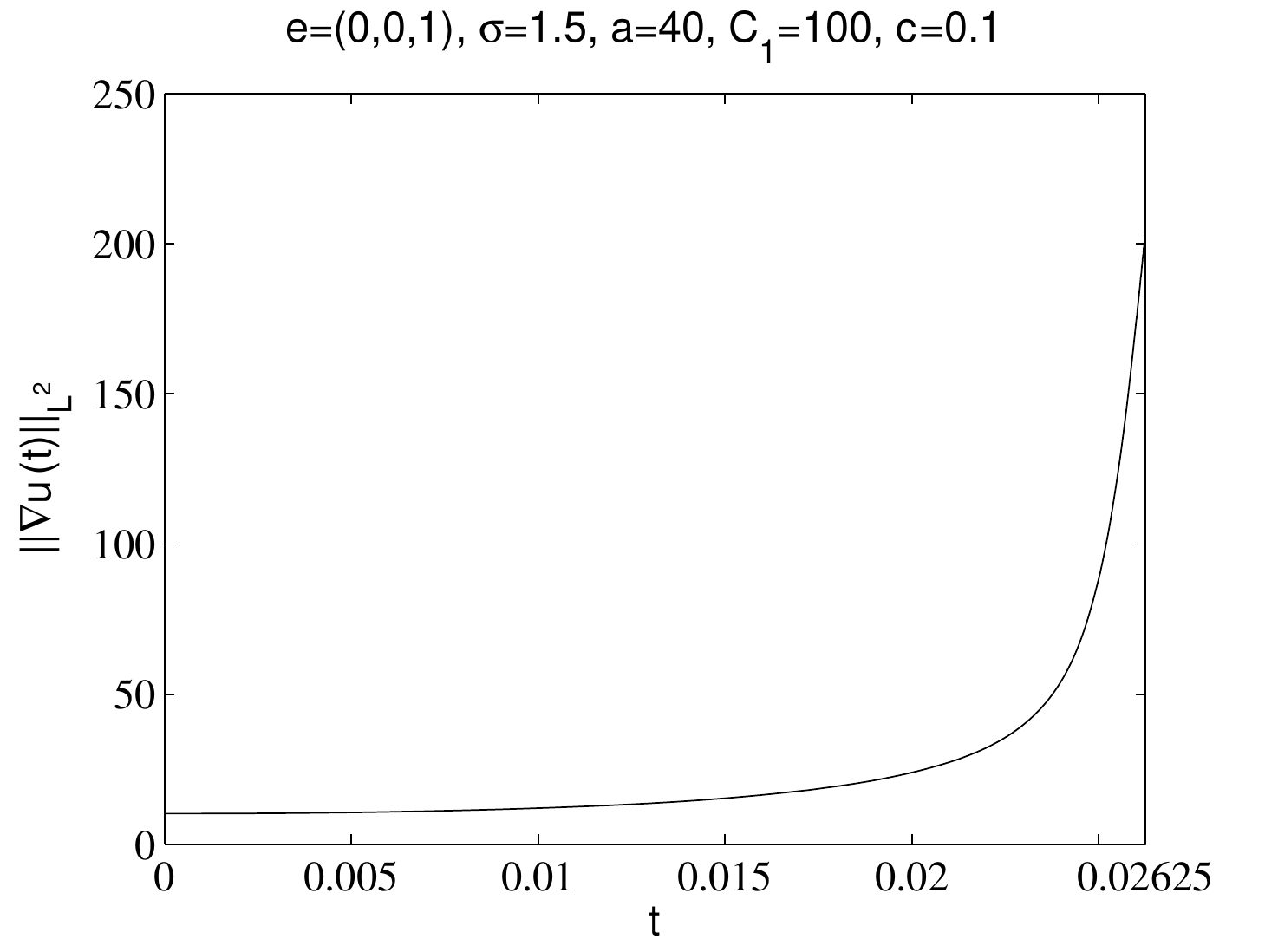}
\includegraphics[width=50mm]{gs_mass.pdf}

$\|\btd u\|_{L^2}$ (given by \eqref{eq:H1}) \hspace{1.5cm} $M(t)$ (given by \eqref{eq:mass})
\caption{Example \ref{ex2}: Graphs of the time-averaged model, the parameters are:
$\eps=1$, $e(t)=(0,0,1)^T$, $\sigma=1.5$, $a=50$, $C_1=100$, $c=0.1$.
Blow-up case.}
\label{Fig004}
\end{figure}

In Fig. \ref{Fig003}--\ref{Fig004}, we can see that the blow-up time decreases as $\sg$ increases.

\begin{figure}[!ht]
\centering
\includegraphics[width=50mm]{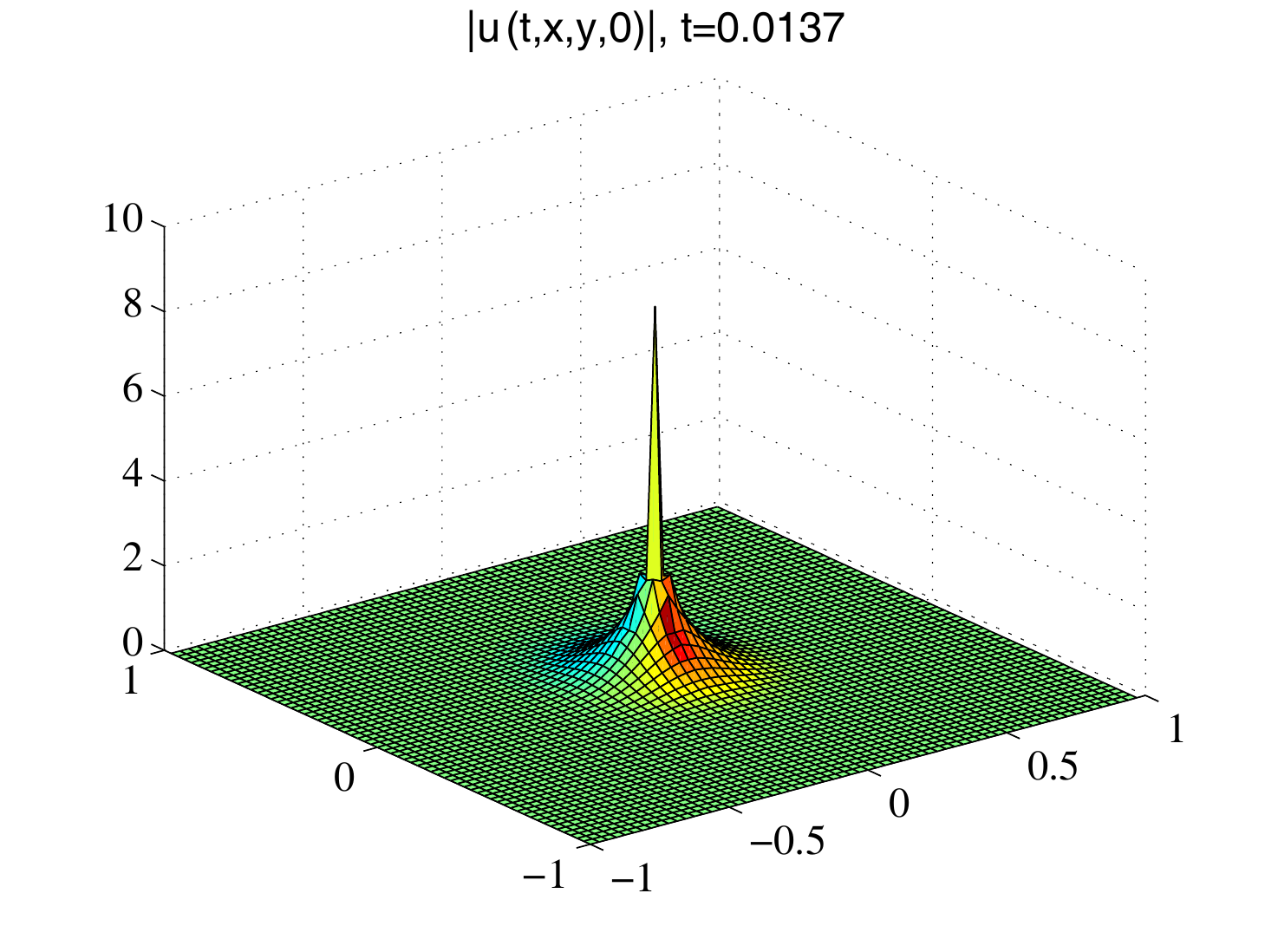}
\includegraphics[width=50mm]{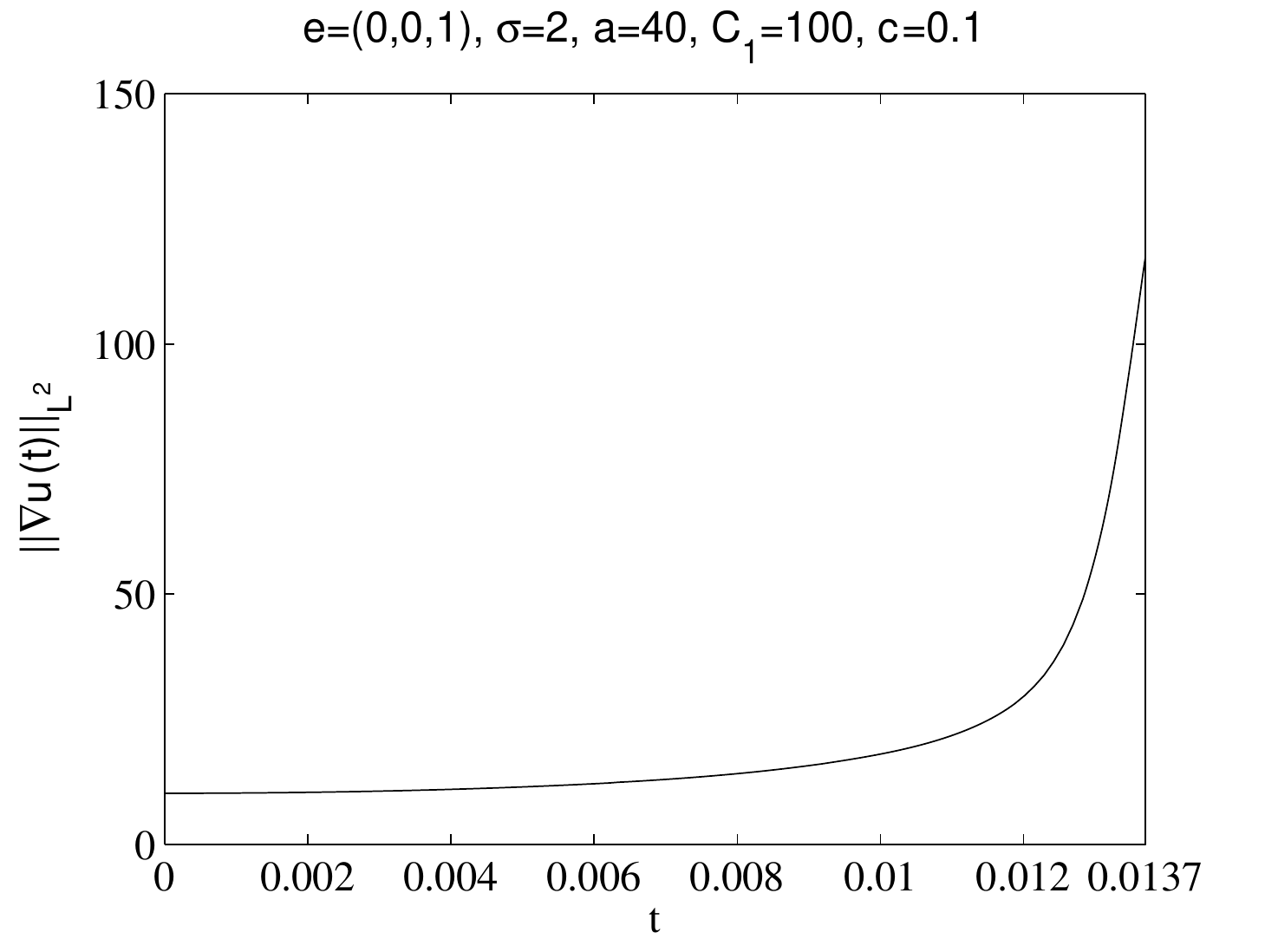}

\hspace{1cm}$u(t,x,y,0)$ \hspace{2.5cm} $\|\btd u\|_{L^2}$ (given by \eqref{eq:H1}) 
\caption{Example \ref{ex2}: Graphs of the fast-potential model, the parameters are:
$\eps=1$, $e(t)=(0,0,1)^T$, $\sigma=2$, $a=50$, $C_1=100$, $c=0.1$, $\og=10^4$.
Blow-up case.}
\label{Fig005}
\end{figure}

\begin{figure}[h!t]
\centering
\includegraphics[width=50mm]{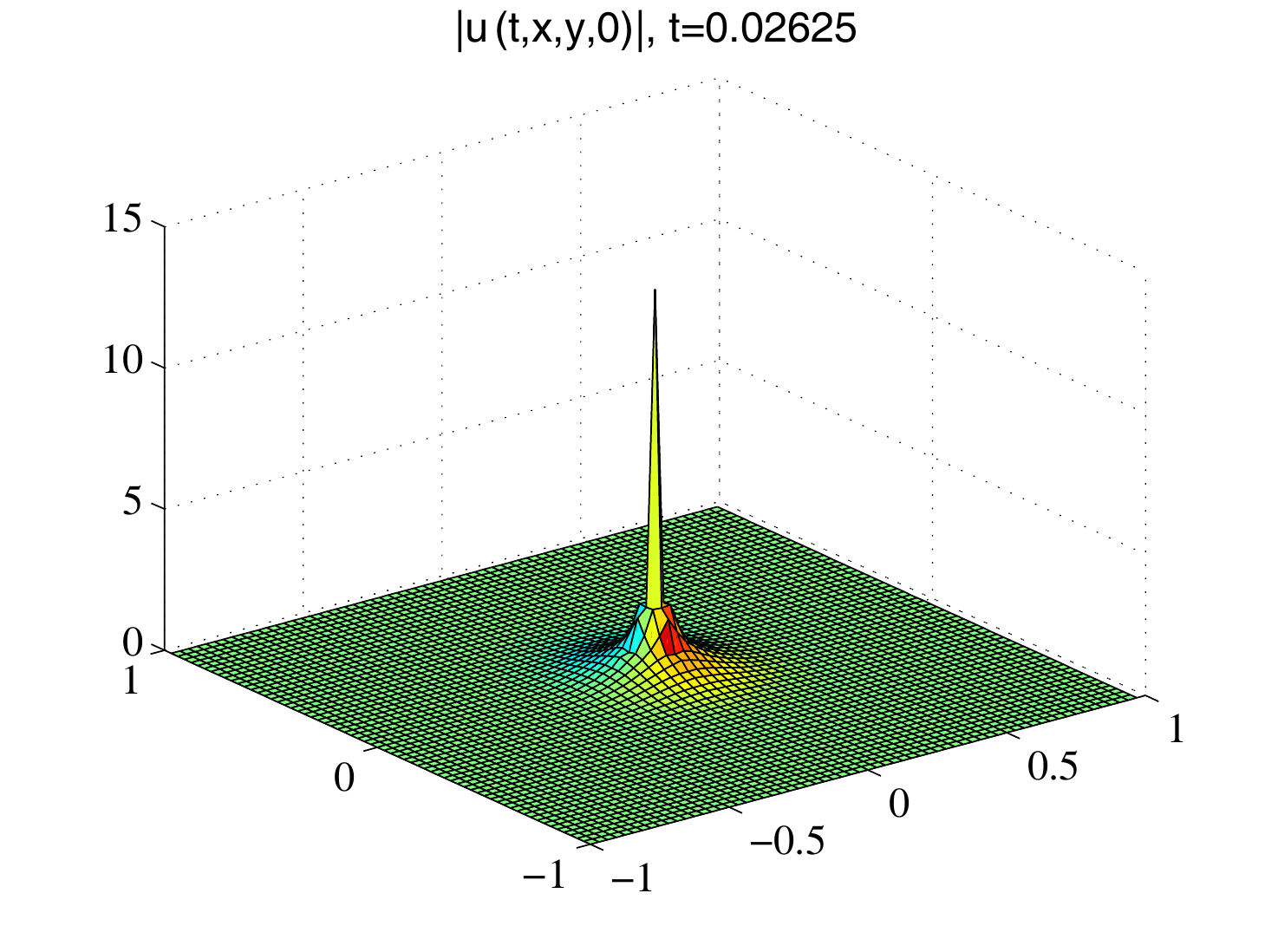}
\includegraphics[width=50mm]{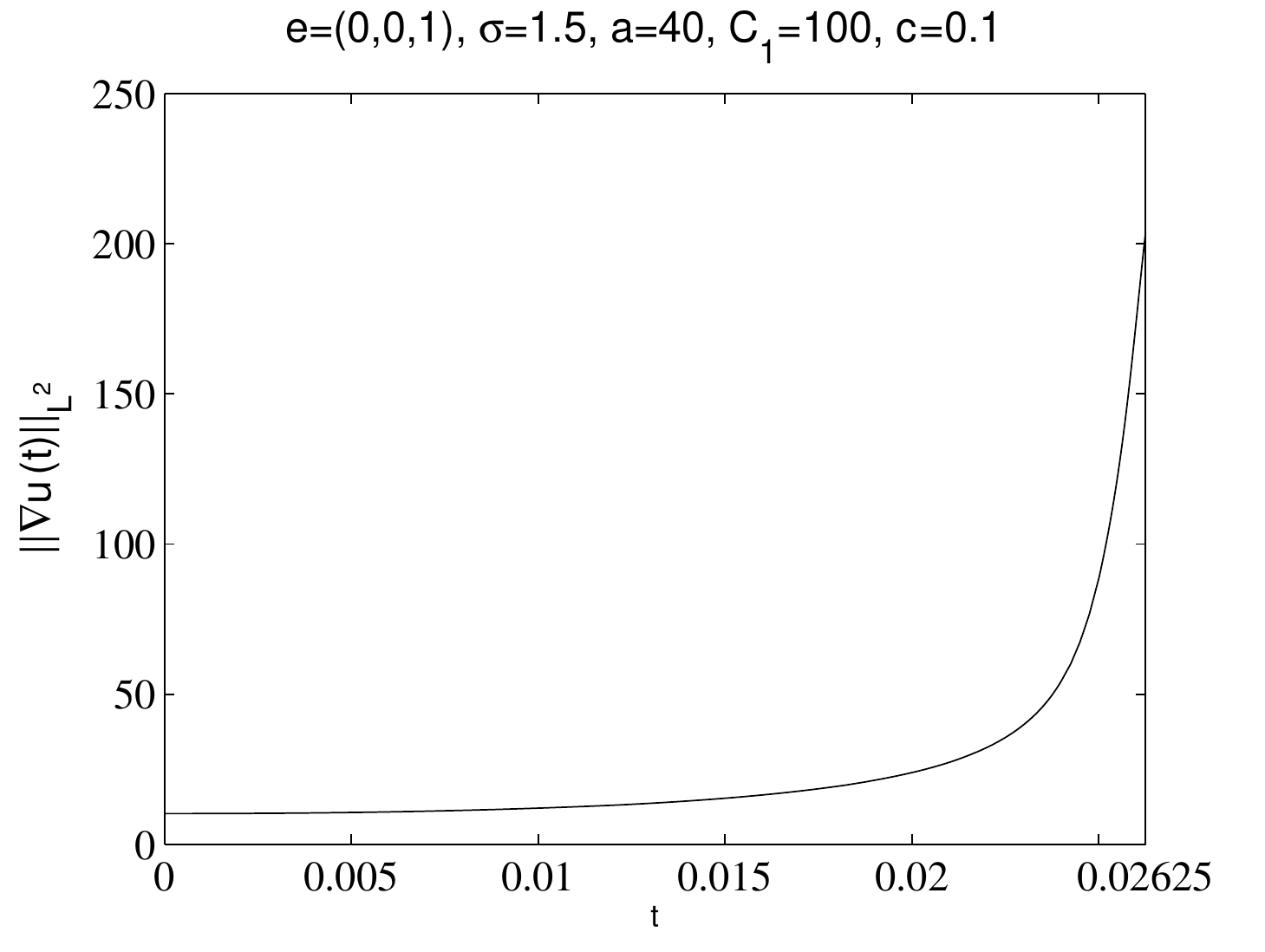}

\hspace{1cm}$u(t,x,y,0)$ \hspace{2.5cm} $\|\btd u\|_{L^2}$ (given by \eqref{eq:H1}) 
\caption{Example \ref{ex2}: Graphs of the fast-potential model, the parameters are:
$\eps=1$, $e(t)=(0,0,1)^T$, $\sigma=1.5$, $a=50$, $C_1=100$, $c=0.1$, $\og=10^4$.}
\label{Fig006}
\end{figure}

\begin{figure}[h!t]
\centering
\includegraphics[width=50mm]{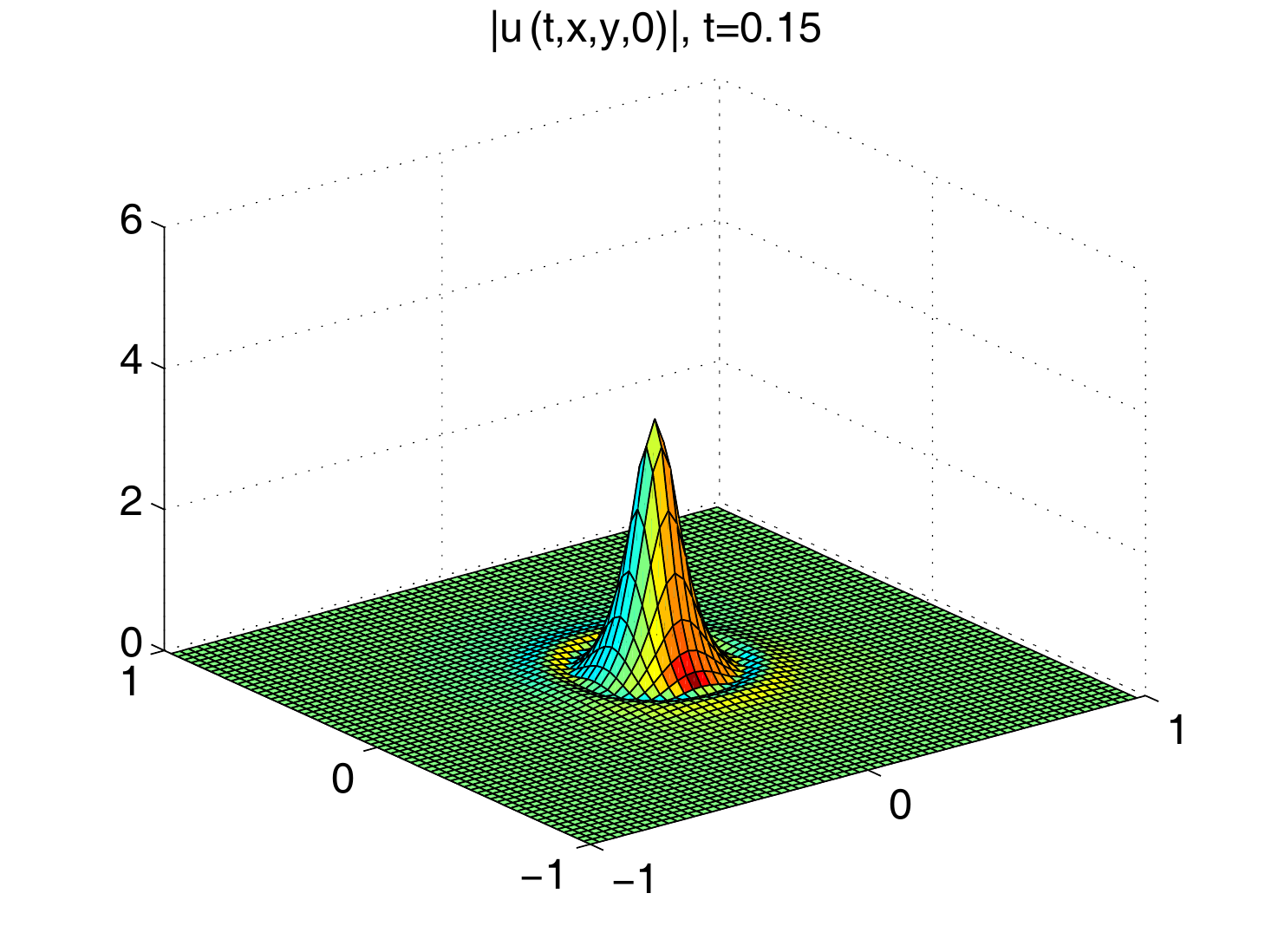}
\includegraphics[width=50mm]{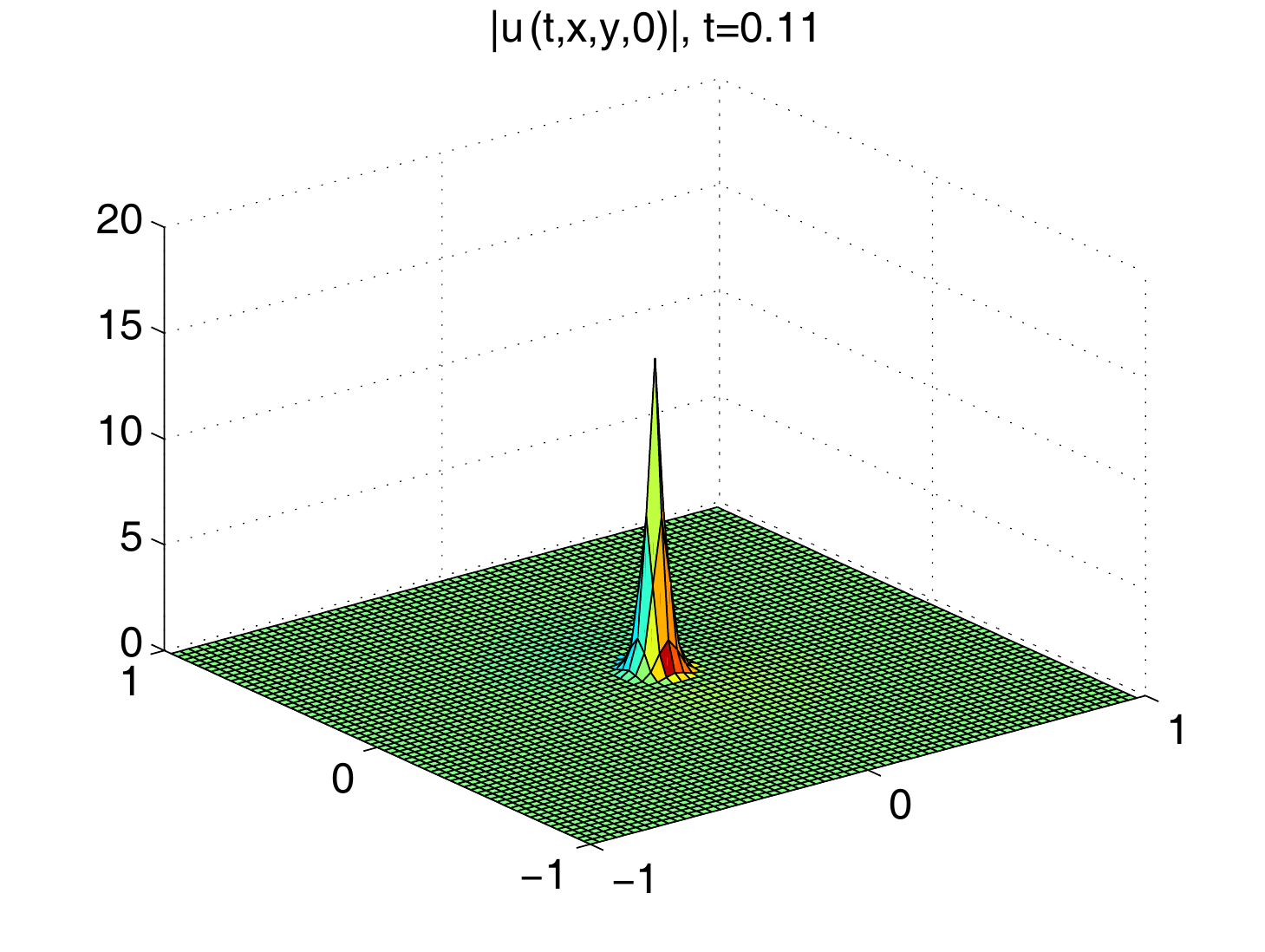}
\includegraphics[width=50mm]{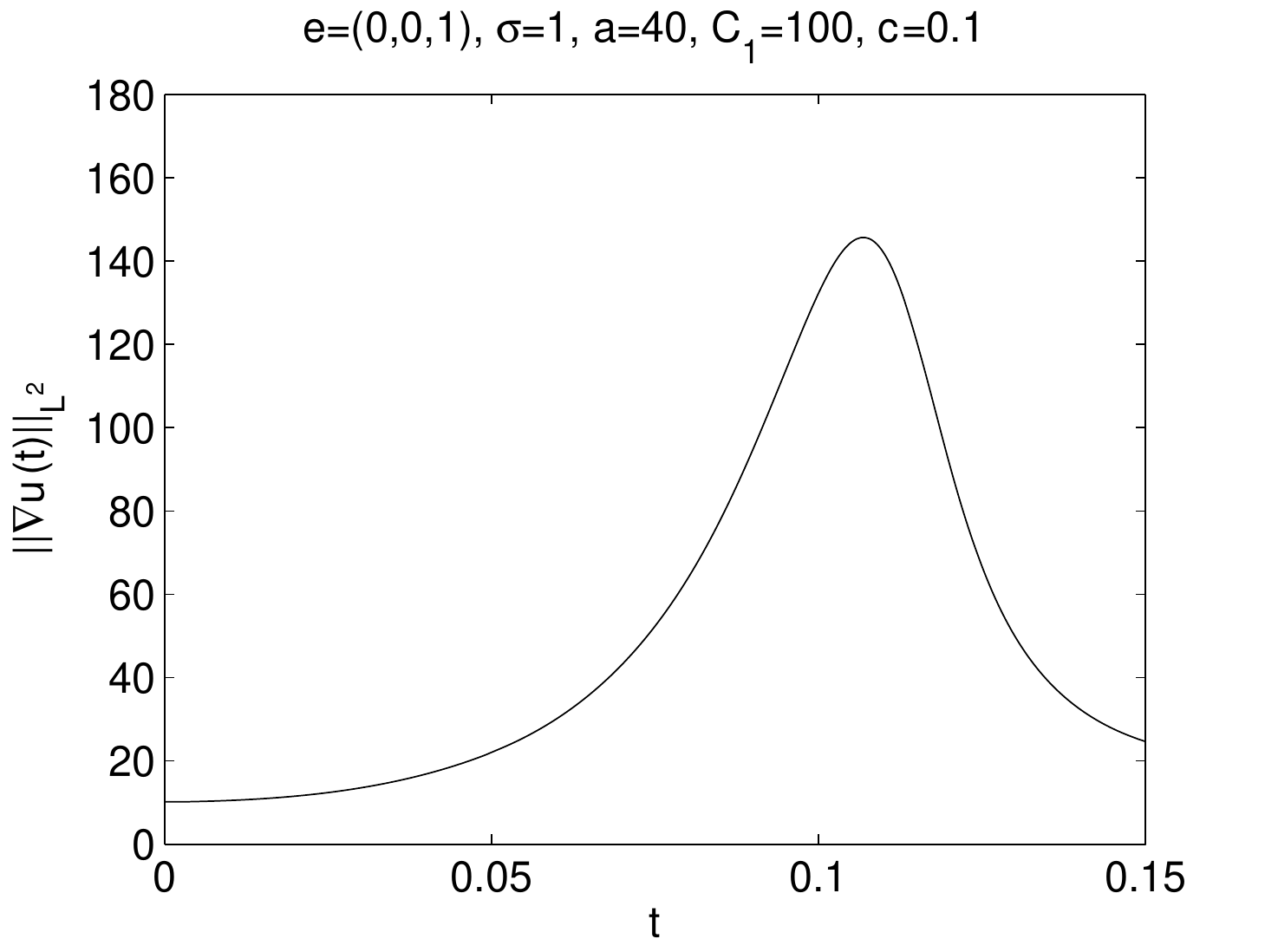}

\hspace{0.6cm}$|u(t,x,y,0)|_{t=0.15}|$ \hspace{2cm} $|u(t,x,y,0)|_{t=0.11}|$ \hspace{1.8cm} $\|\btd u\|_{L^2}$ (given by \eqref{eq:H1})\vspace{2mm}

\includegraphics[width=50mm]{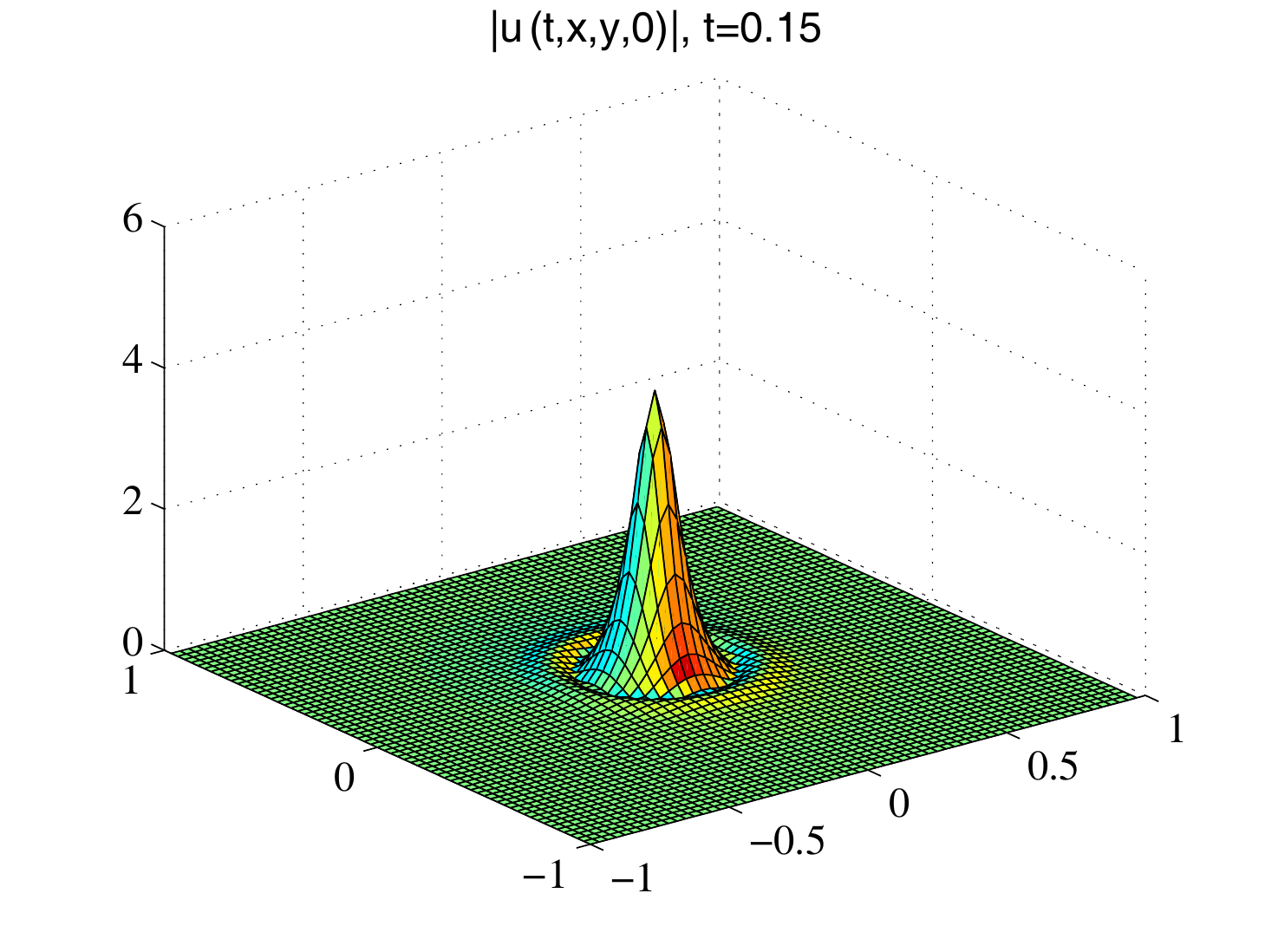}
\includegraphics[width=50mm]{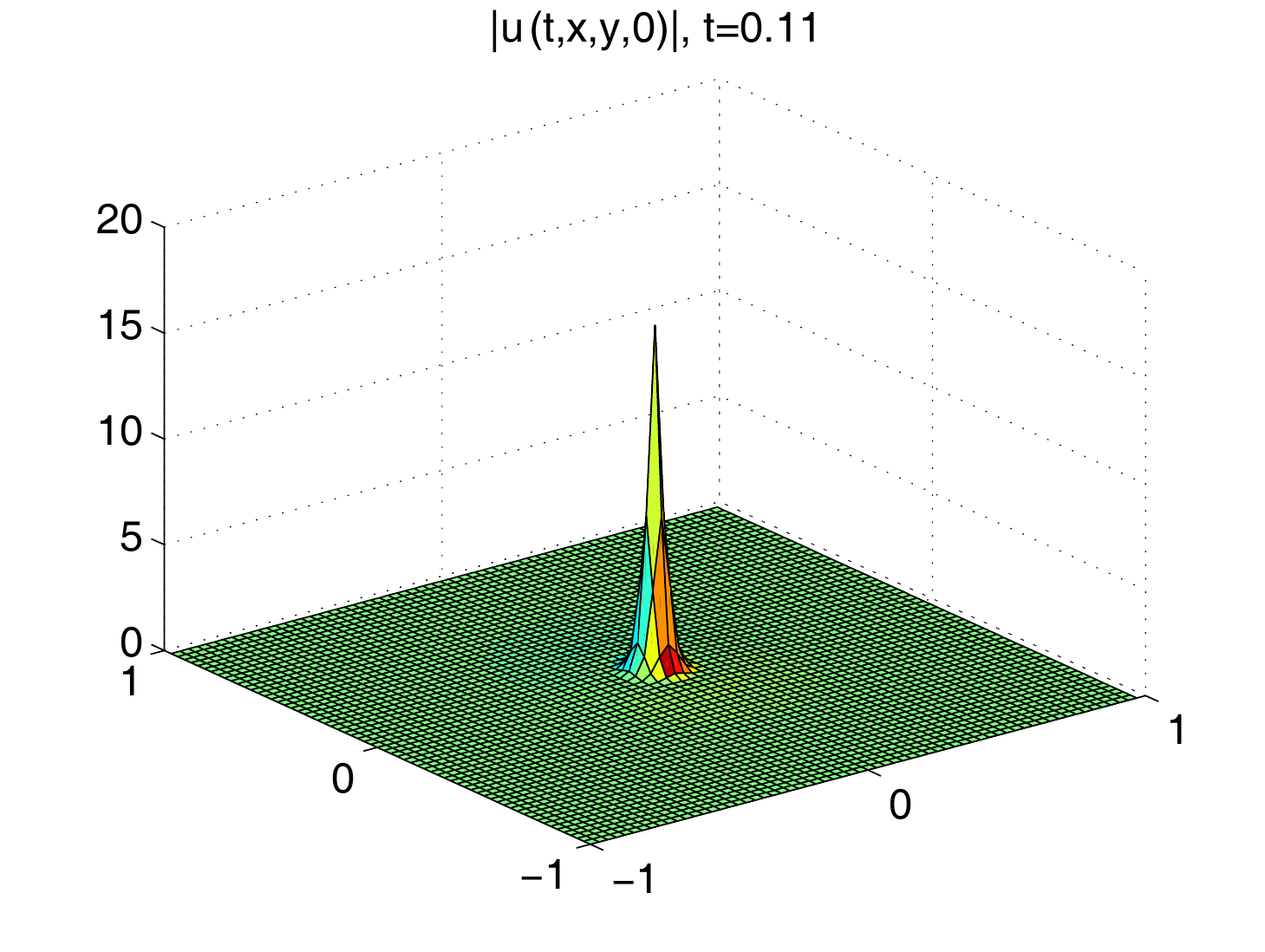}
\includegraphics[width=50mm]{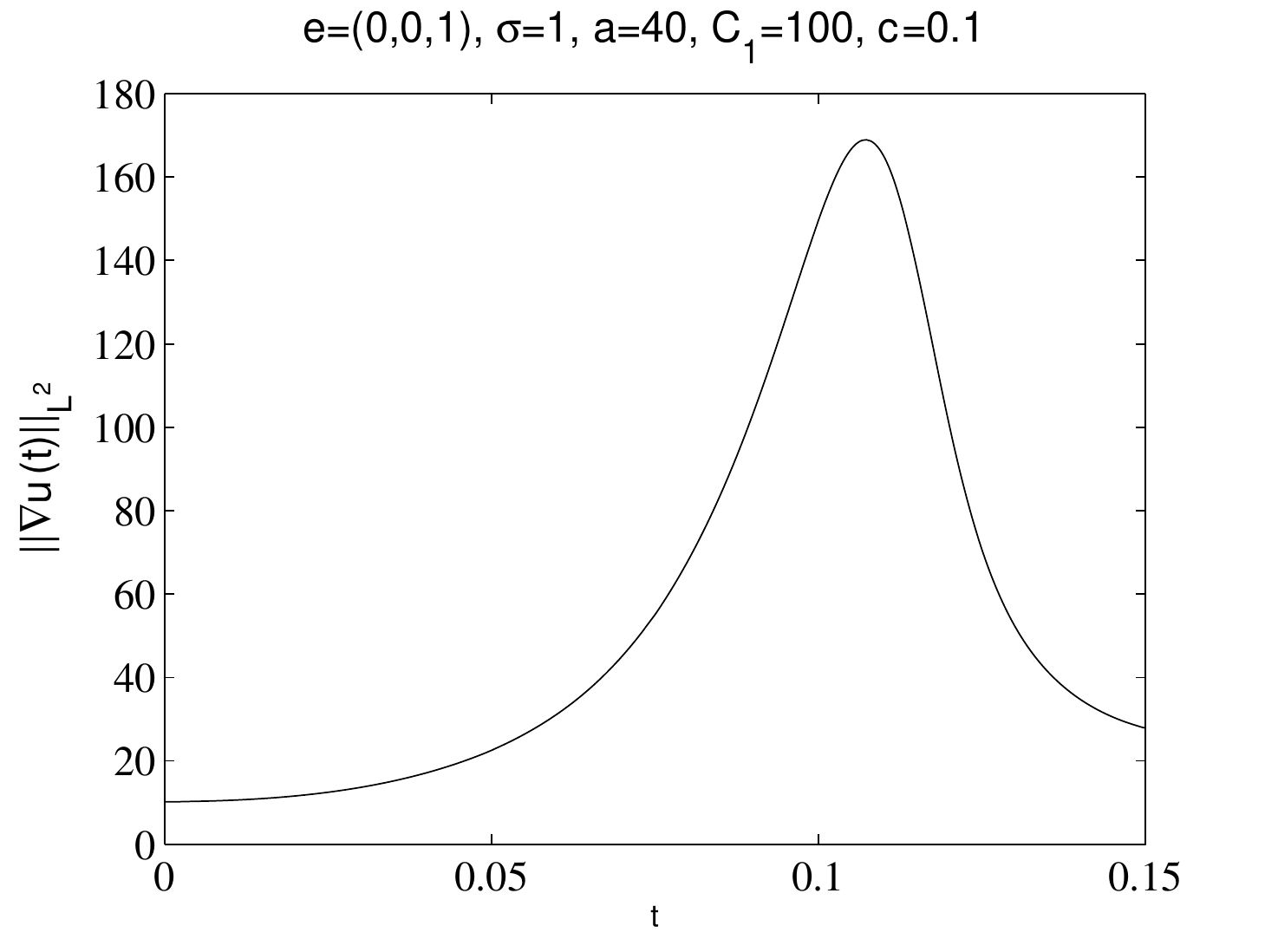}

\hspace{0.6cm}$|u(t,x,y,0)|_{t=0.15}|$ \hspace{2cm} $|u(t,x,y,0)|_{t=0.11}|$ \hspace{1.8cm} $\|\btd u\|_{L^2}$ (given by \eqref{eq:H1})

\caption{Example \ref{ex2}: Graphs of the fast-potential model (first line) and time-averaged model (second line), the parameters are:
$\eps=1$, $e(t)=(0,0,1)^T$, $\sigma=1$, $a=50$, $C_1=100$, $c=0.1$, $\omega=10^4$.}
\label{Fig007}
\end{figure}

\begin{figure}[h!t]
\centering
\includegraphics[width=50mm]{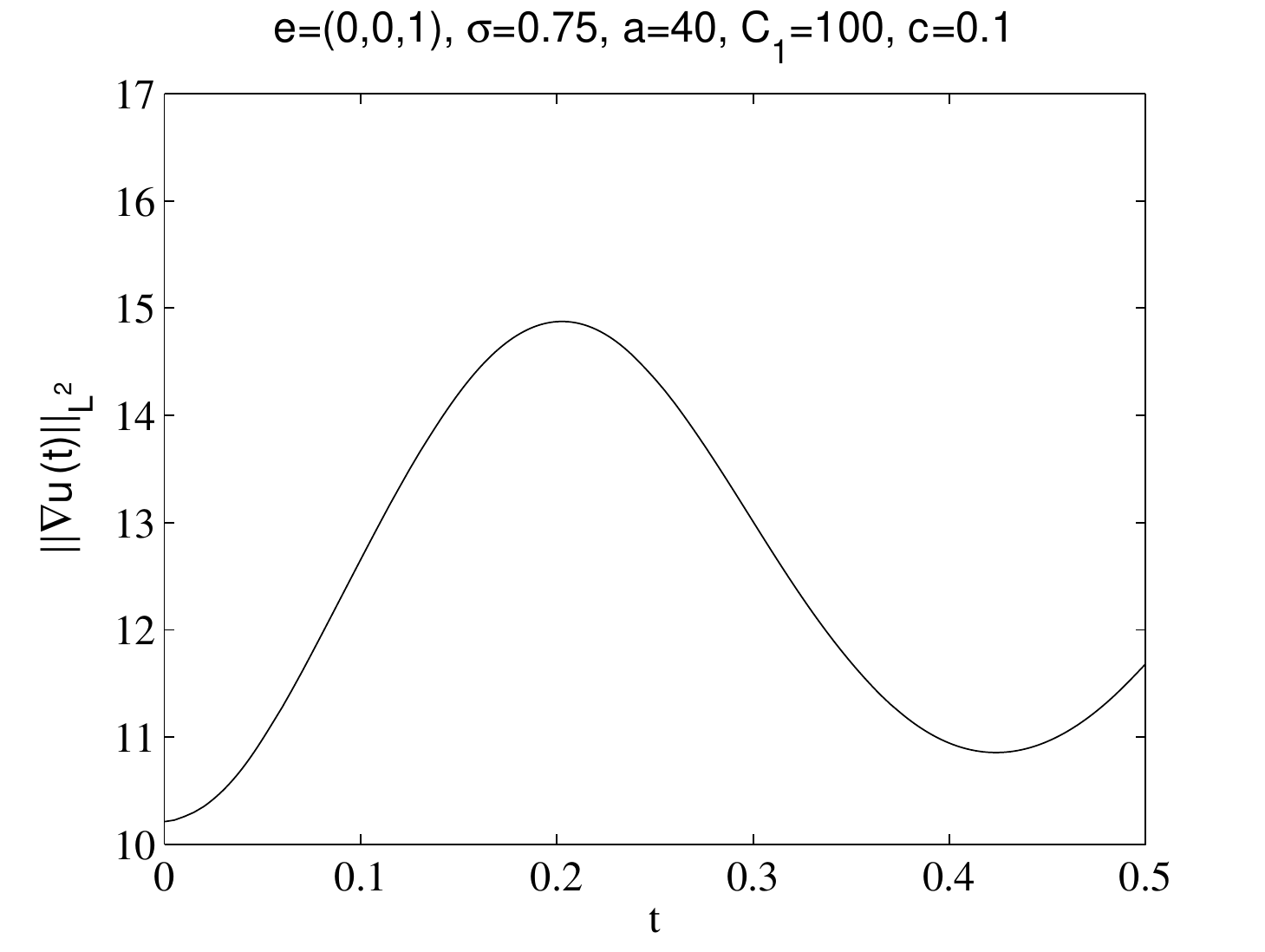}
\includegraphics[width=50mm]{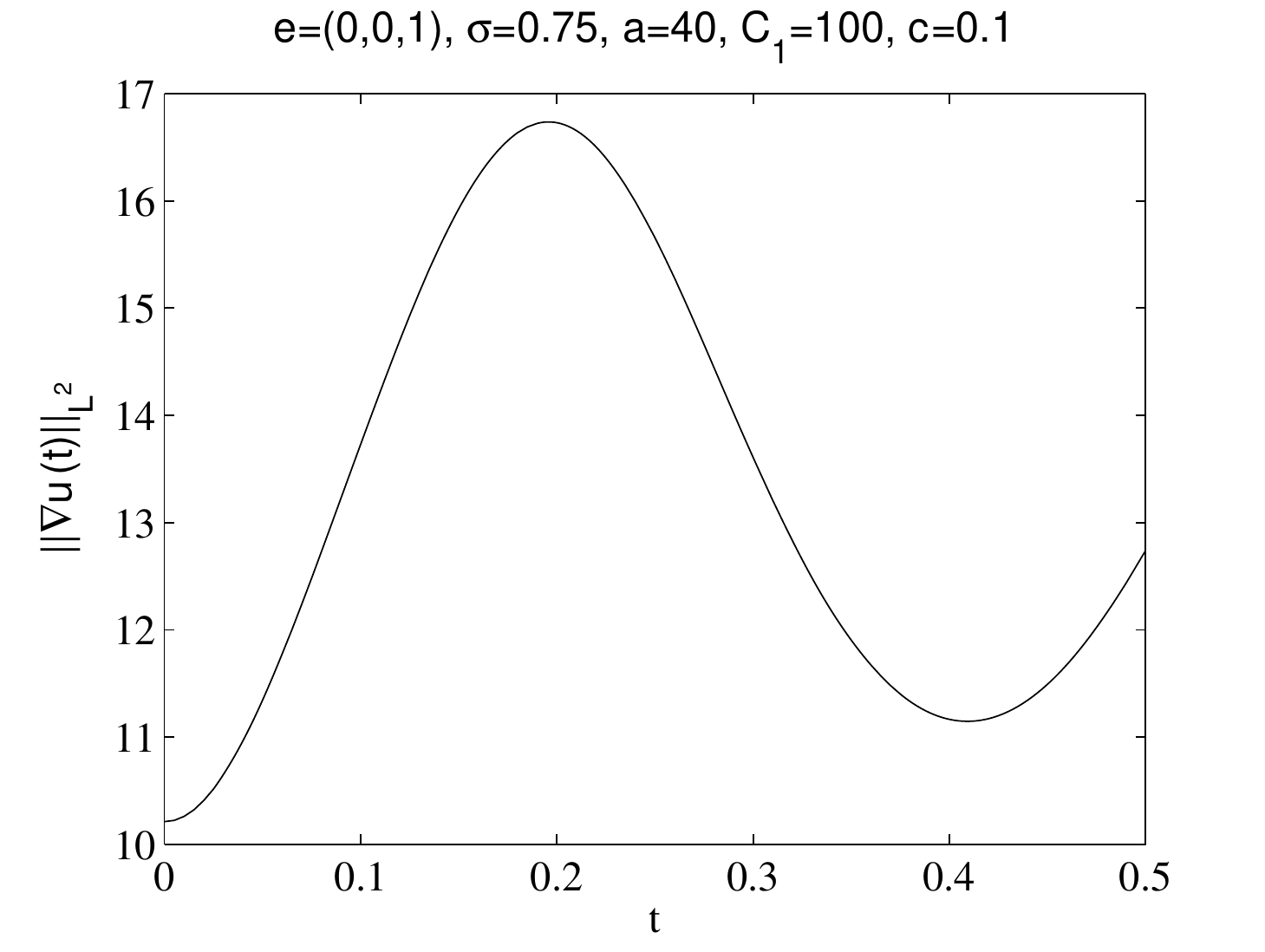}

$\|\btd u^\og\|_{L^2}$ (given by \eqref{eq:H1}) \hspace{1.1cm} $\|\btd u\|_{L^2}$ (given by \eqref{eq:H1})

\caption{Example \ref{ex2}: Graphs of the fast-potential model (left) and time-averaged model (right), the parameters are:
$\eps=1$, $e(t)=(0,0,1)^T$, $\sigma=0.75$, $a=50$, $C_1=100$, $c=0.1$, $\omega=10^4$.}
\label{Fig008}
\end{figure}

Comparing the $L^2$ norms of the gradient of the wave-function for the time averaged model and the fast model in the blow up test we can see that they are very similar in the super-critical case
(\cf Figures \ref{Fig005}--\ref{Fig006} and Figures \ref{Fig003}--\ref{Fig004}).
They exhibit more differences in the subcritical case (see Fig. \ref{Fig007}--\ref{Fig008}).

\begin{figure}[h!t]
\centering
\includegraphics[width=50mm]{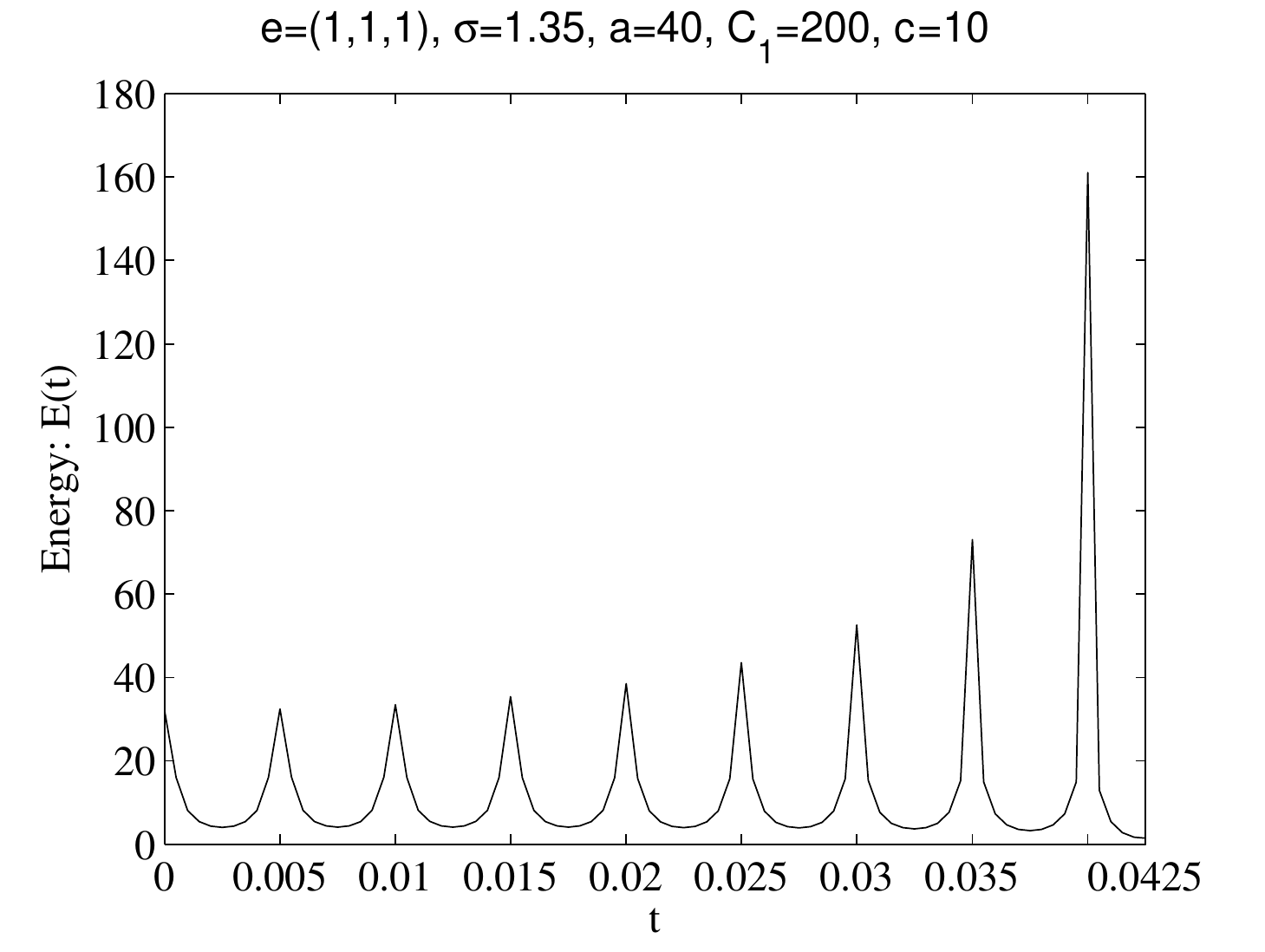}
\includegraphics[width=50mm]{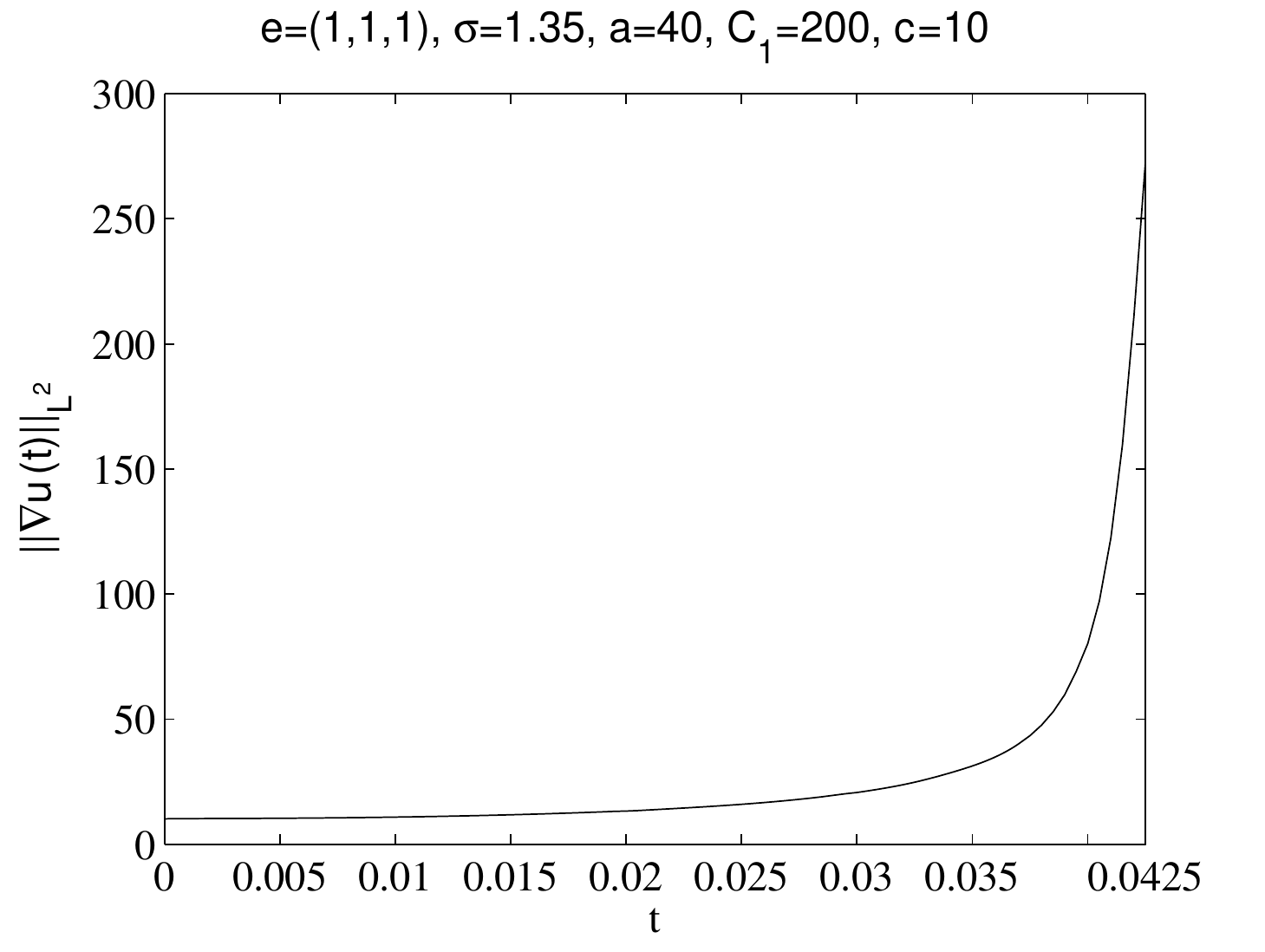}
\includegraphics[width=50mm]{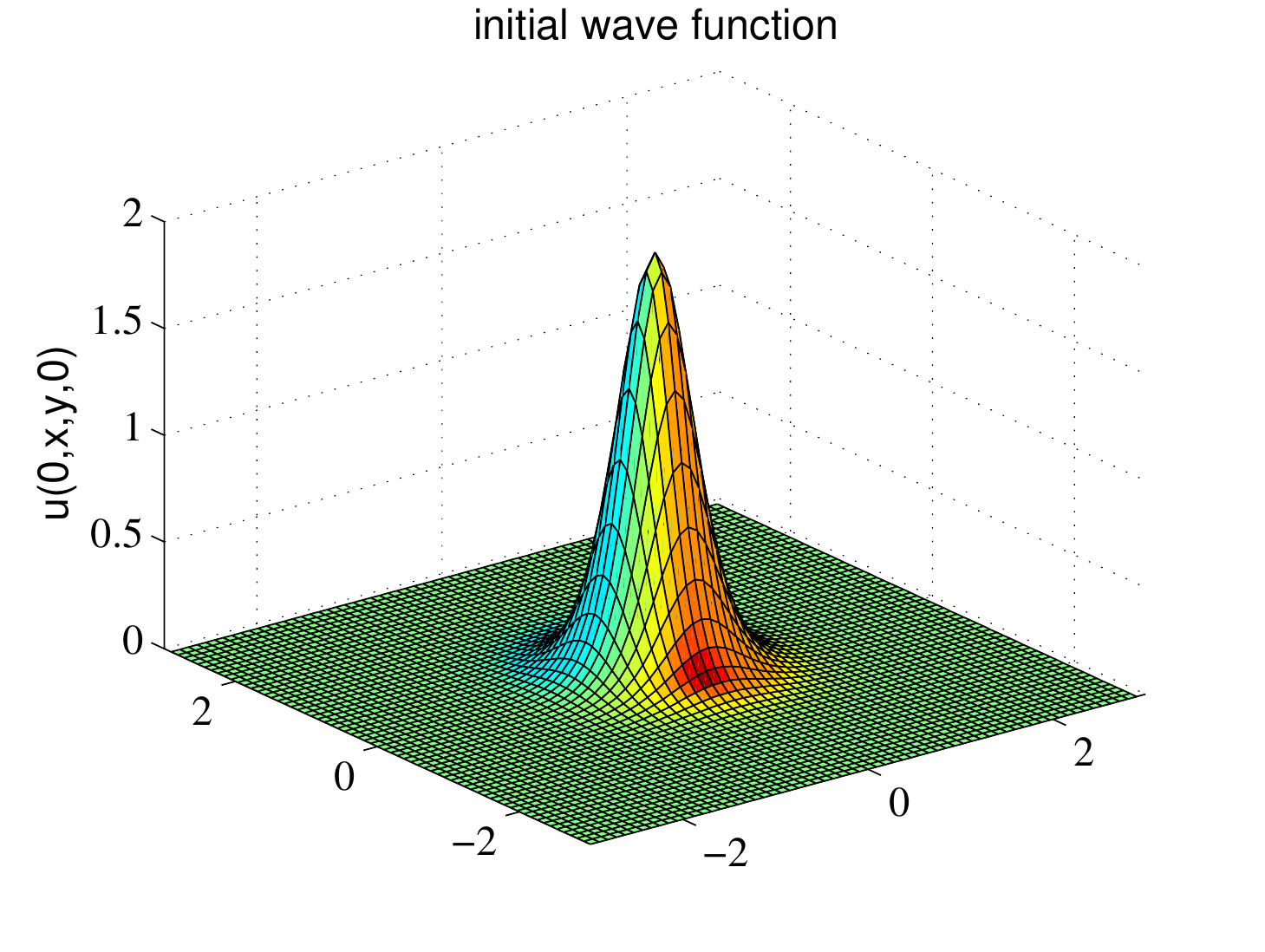}

\hspace{0.6cm}$E(t)$ (given by \eqref{eq:energy1}) \hspace{1.6cm} $\|\btd u\|_{L^2}$ (given by \eqref{eq:H1})\hspace{2.2cm} $|u(t,x,y,0)|_{t=0}|$ \vspace{4mm}

\includegraphics[width=50mm]{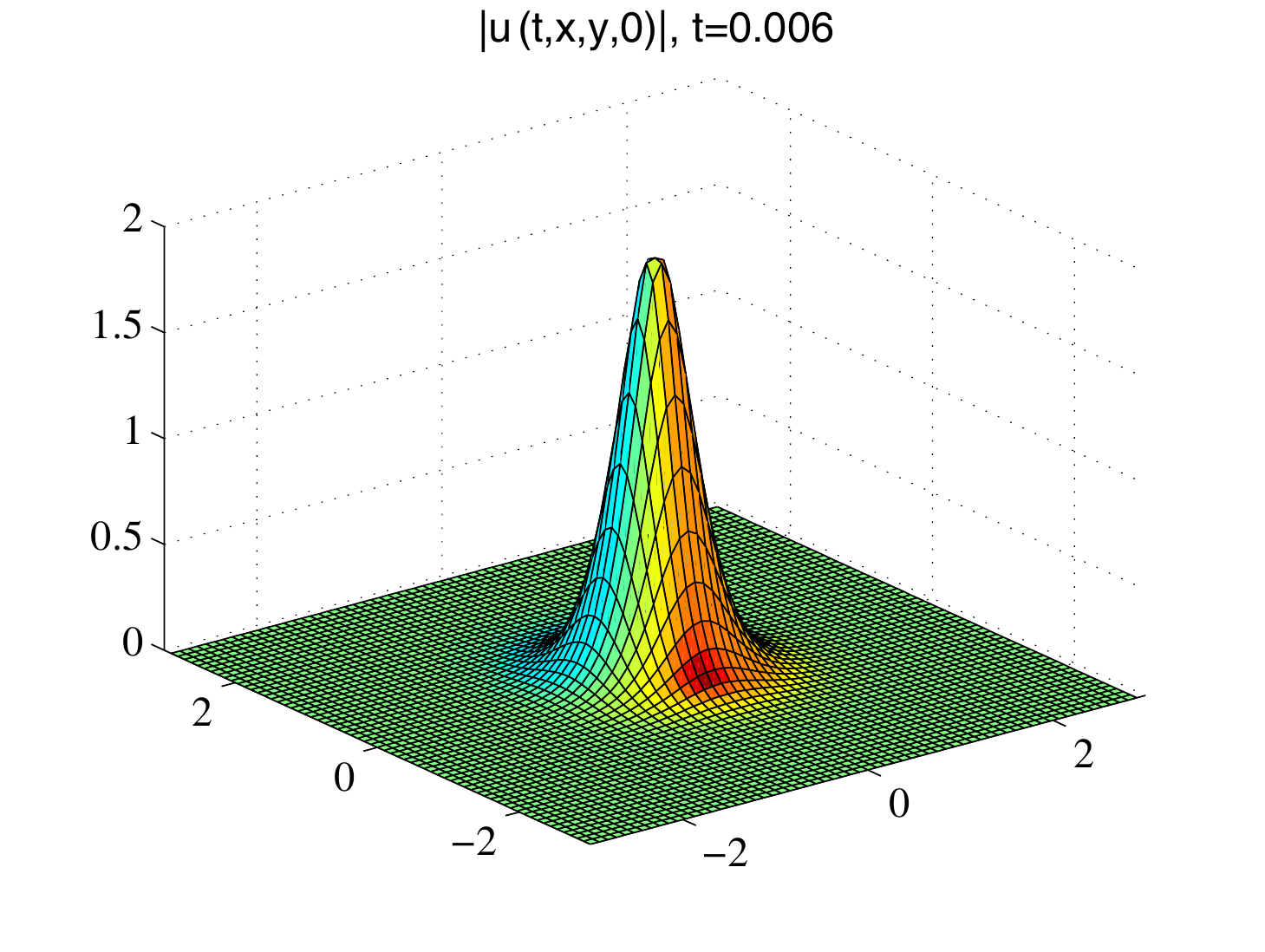}
\includegraphics[width=50mm]{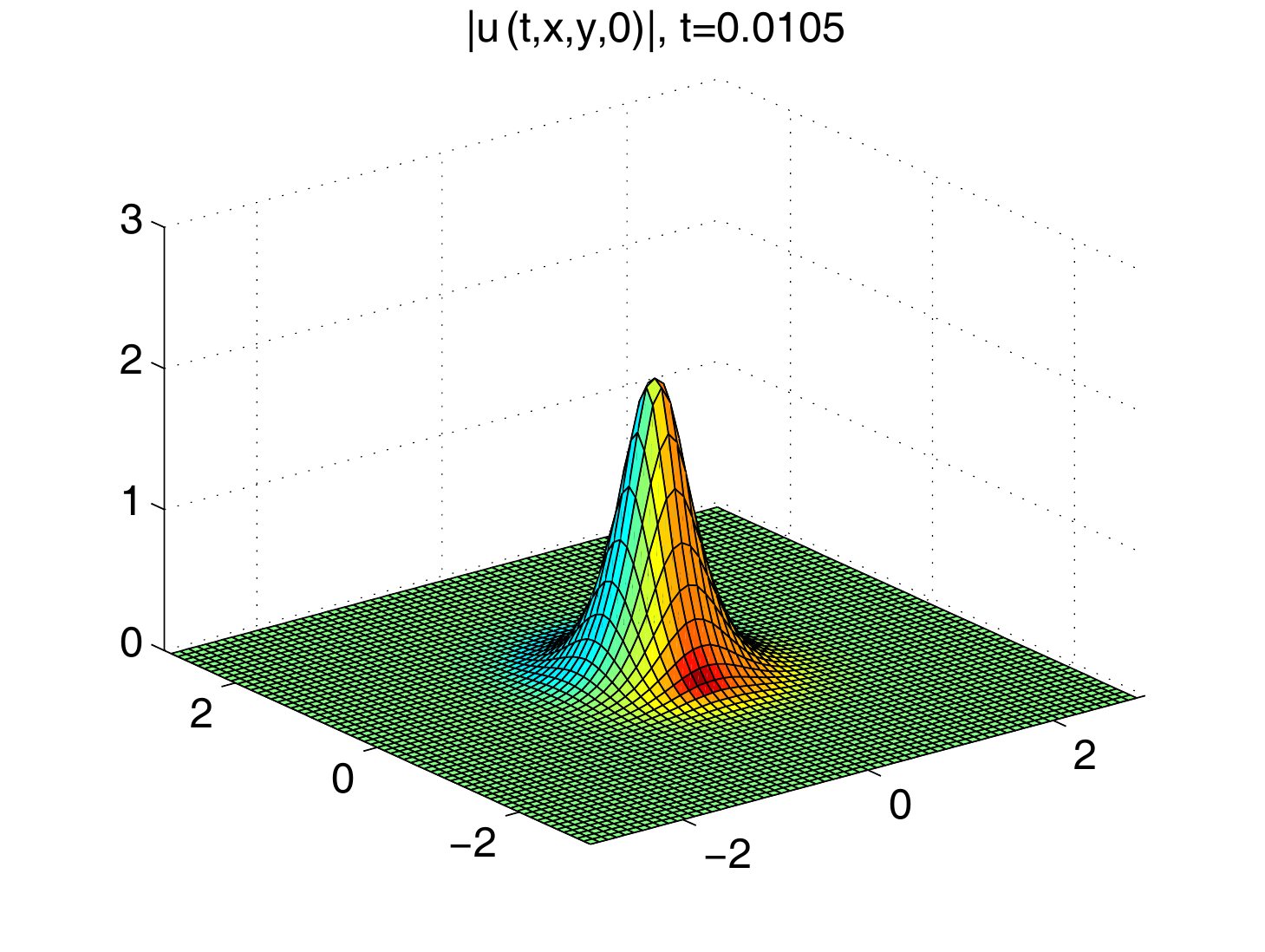}
\includegraphics[width=50mm]{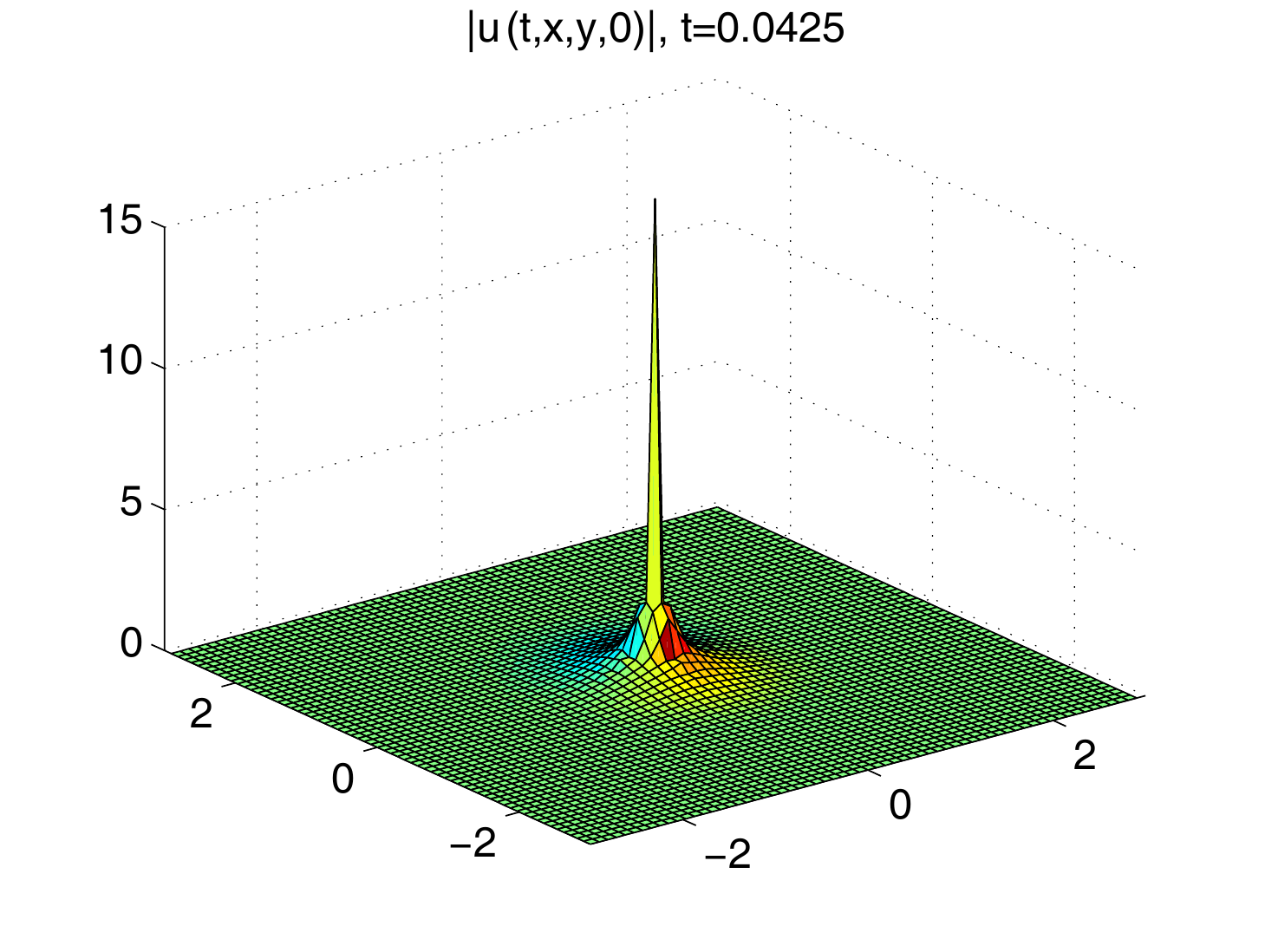}

\hspace{0.6cm}$|u(t,x,y,0)|_{t=0.006}|$ \hspace{1.8cm} $|u(t,x,y,0)|_{t=0.0105}|$ \hspace{1.8cm} $|u(t,x,y,0)|_{t=0.0425}|$ \hspace{2mm}

\caption{Example \ref{ex2}: Graphs of the fast-potential model, the parameters are:
$\eps=1$, $e(t)=(1,1,1)^T$, $\sigma=1.35$, $a=40$, $C_1=200$, $c=10$, $\omega=10^2$.}
\label{Fig009}
\end{figure}

\begin{figure}[h!t]
\centering
\includegraphics[width=40mm]{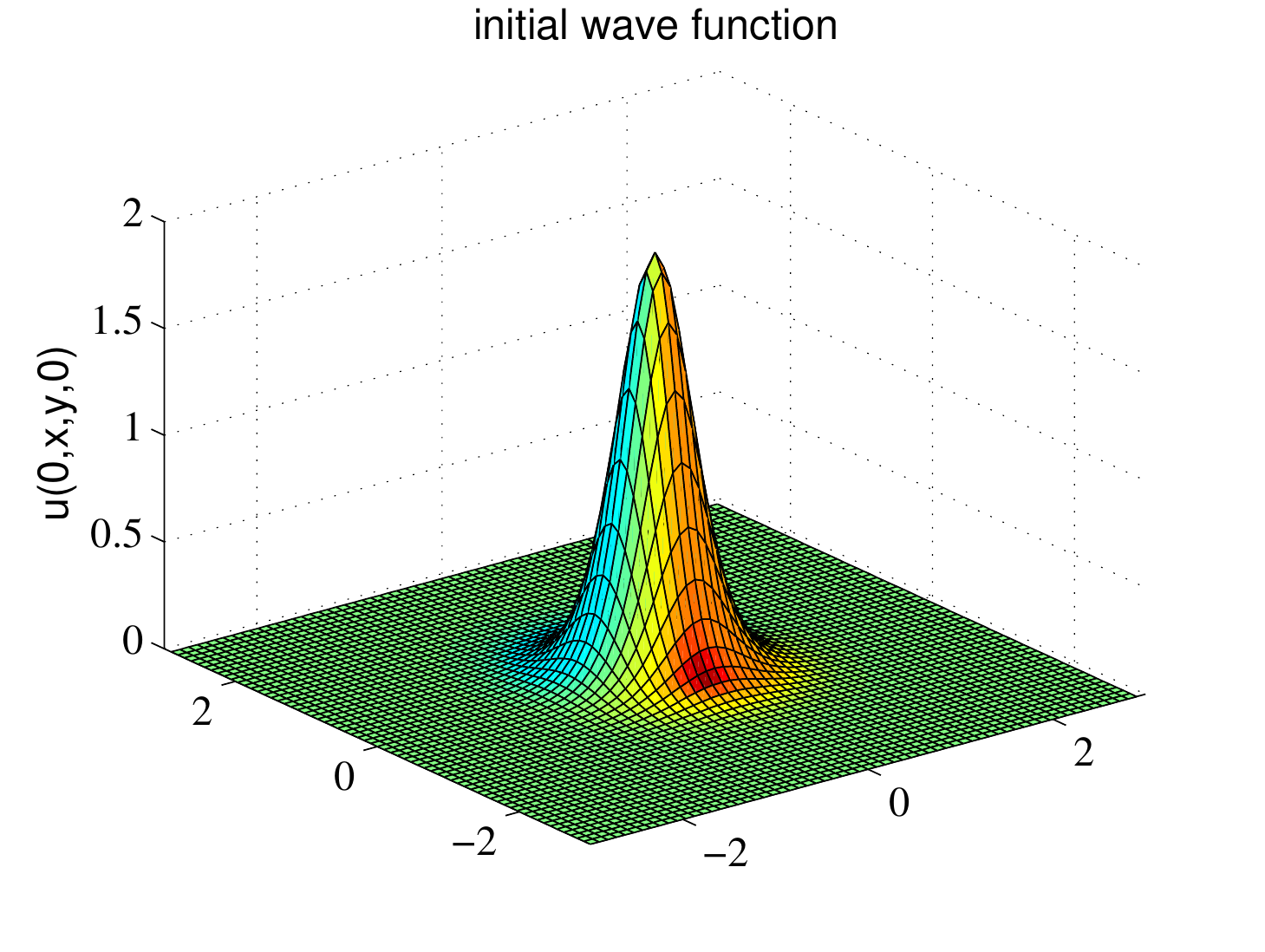}
\includegraphics[width=40mm]{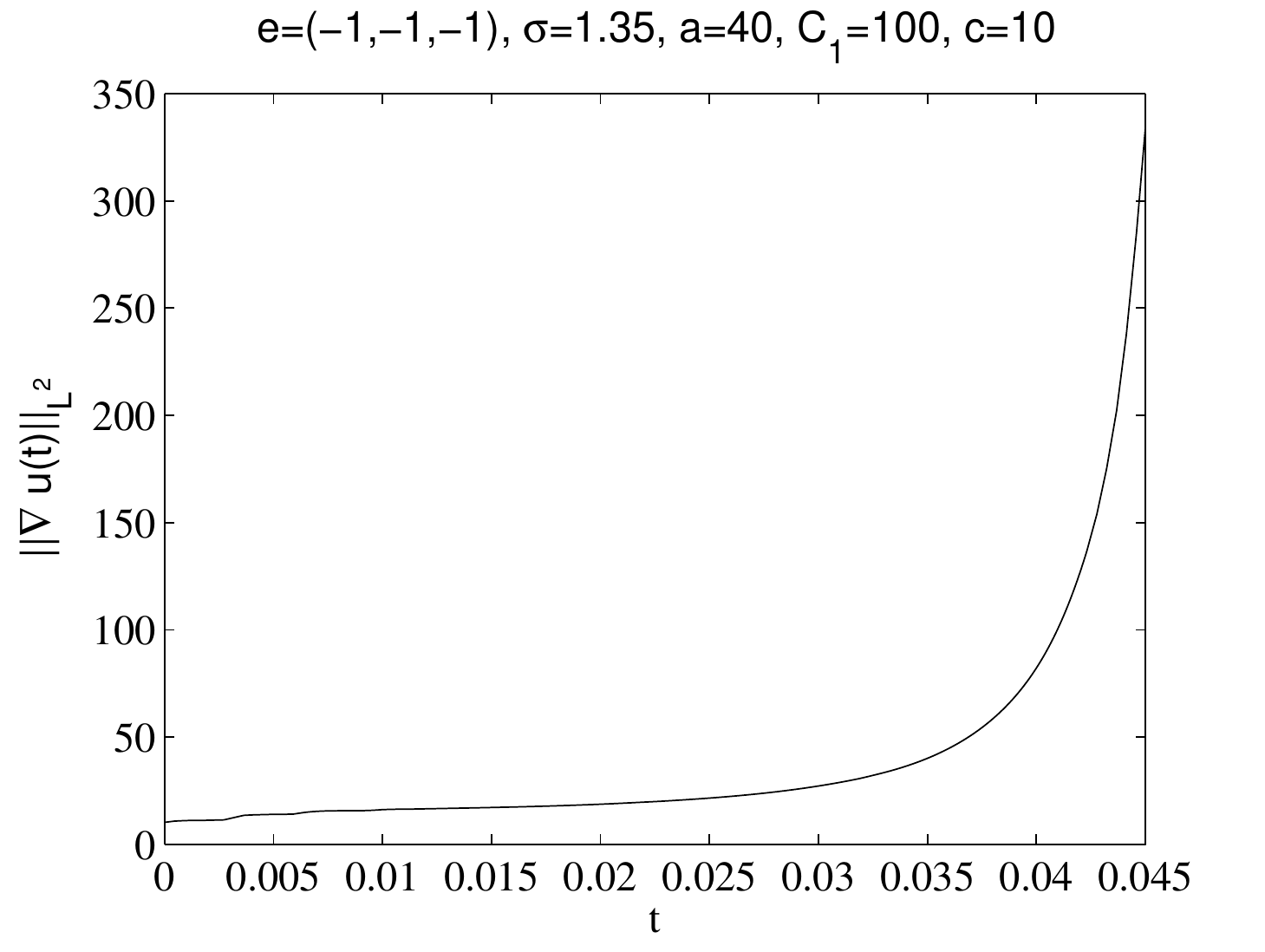}
\includegraphics[width=40mm]{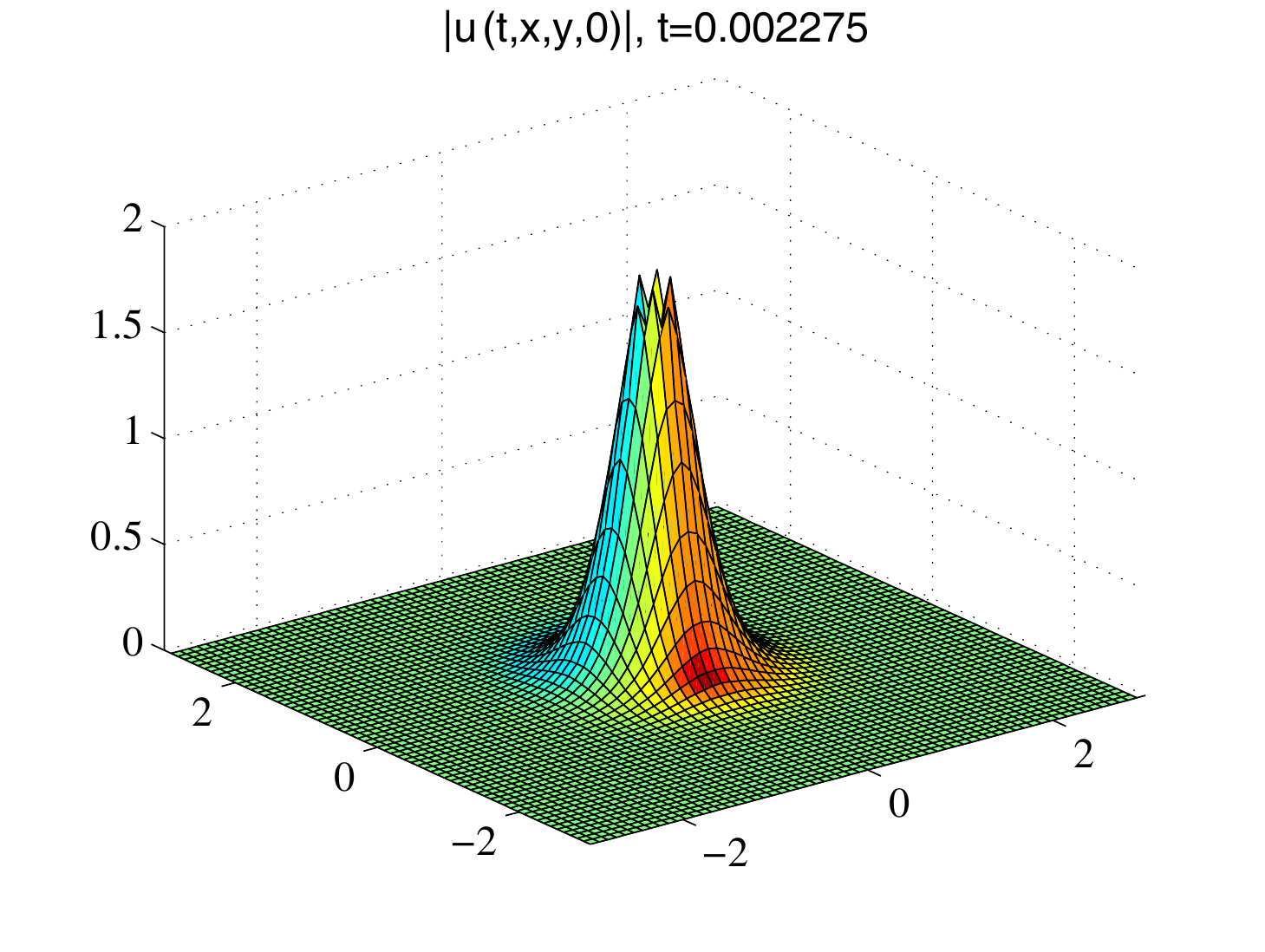}
\includegraphics[width=40mm]{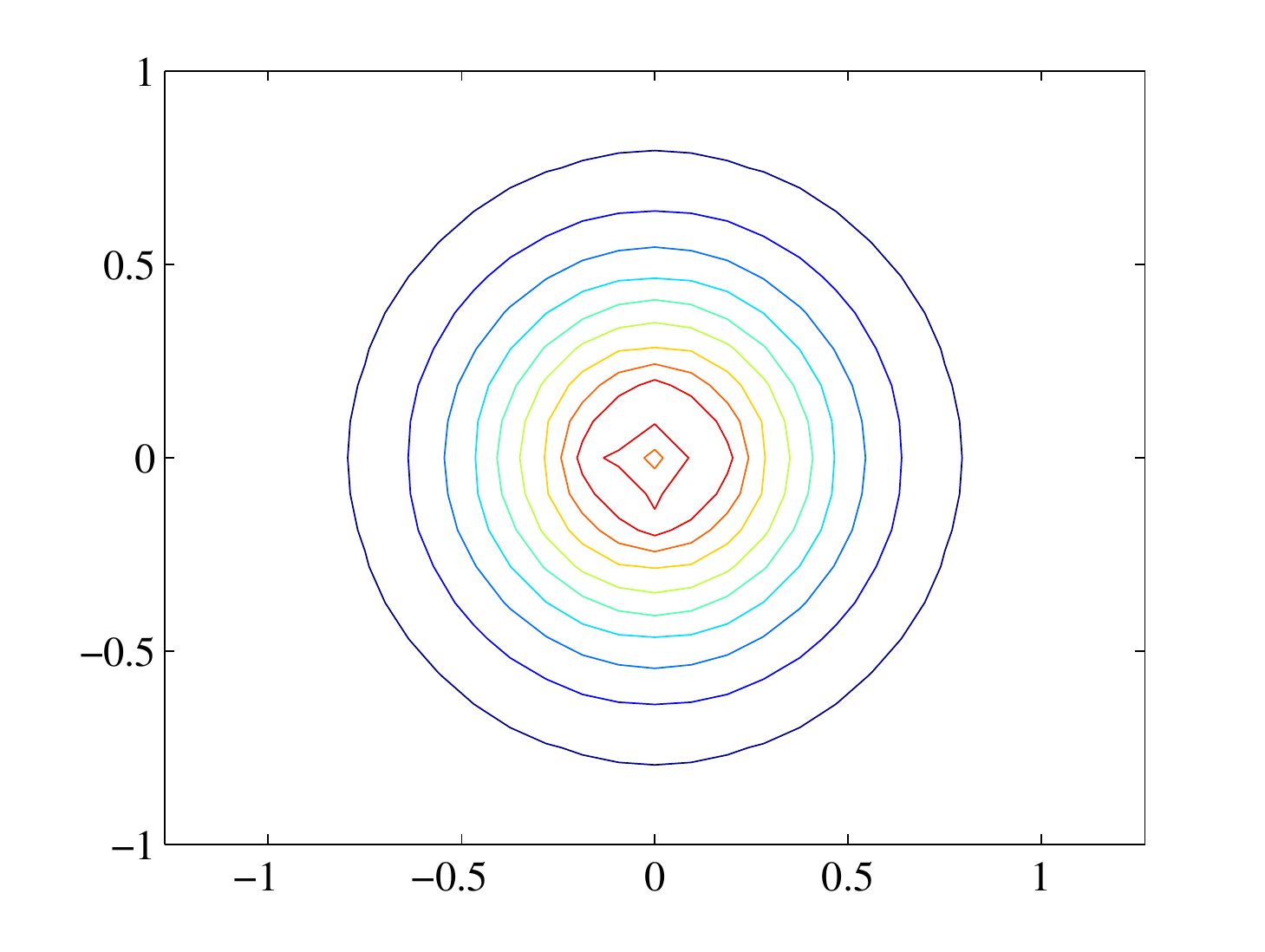}

$|u(0,x,y,0)|$ \hspace{1.2cm} $\|\btd u\|_{L^2}$ (given by \eqref{eq:H1})\hspace{0.5cm}  $|u(t,x,y,0)|_{t=0.002275}|$\hspace{1.2cm} contour plot\vspace{5mm}

\includegraphics[width=40mm]{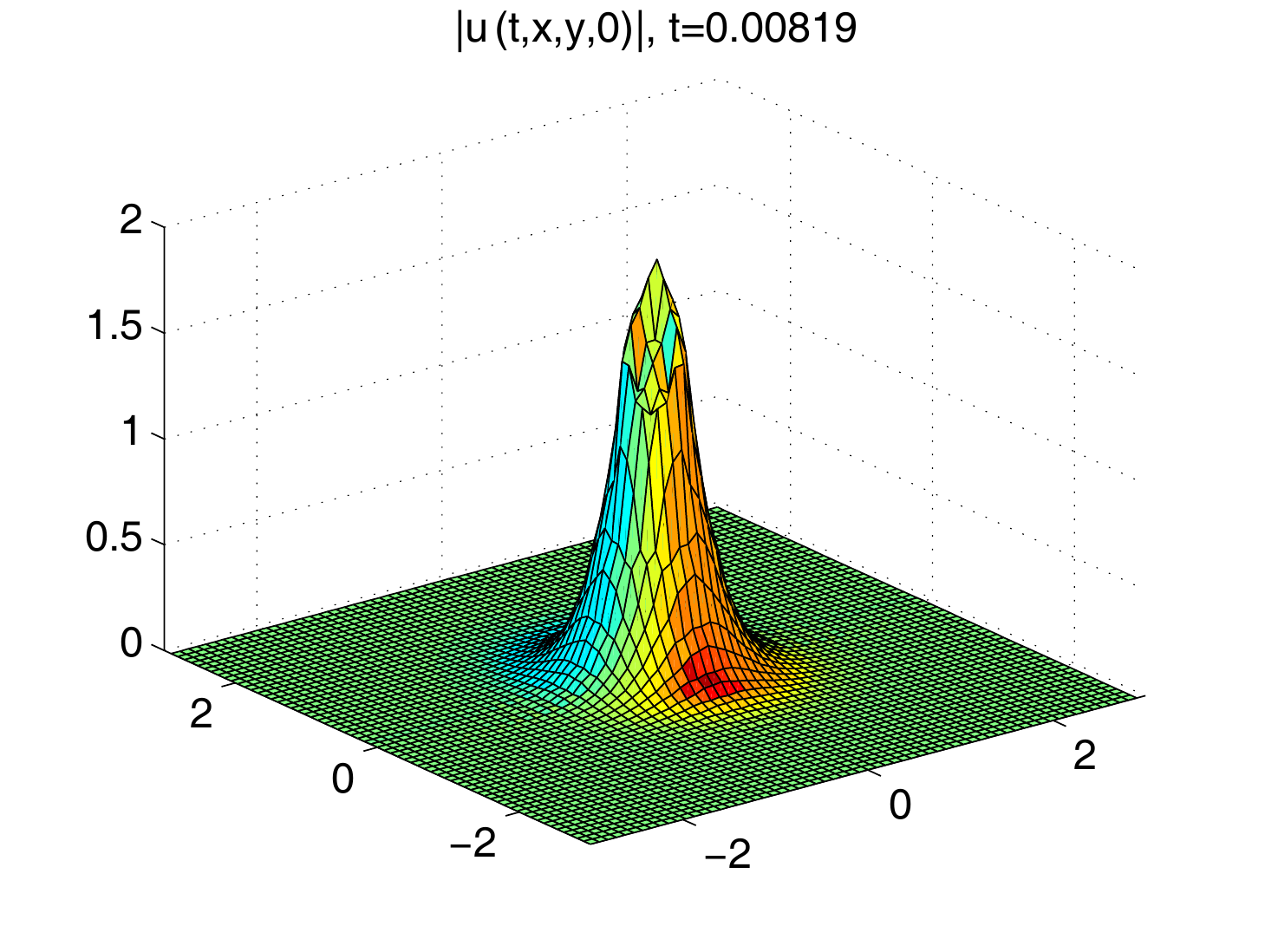}
\includegraphics[width=40mm]{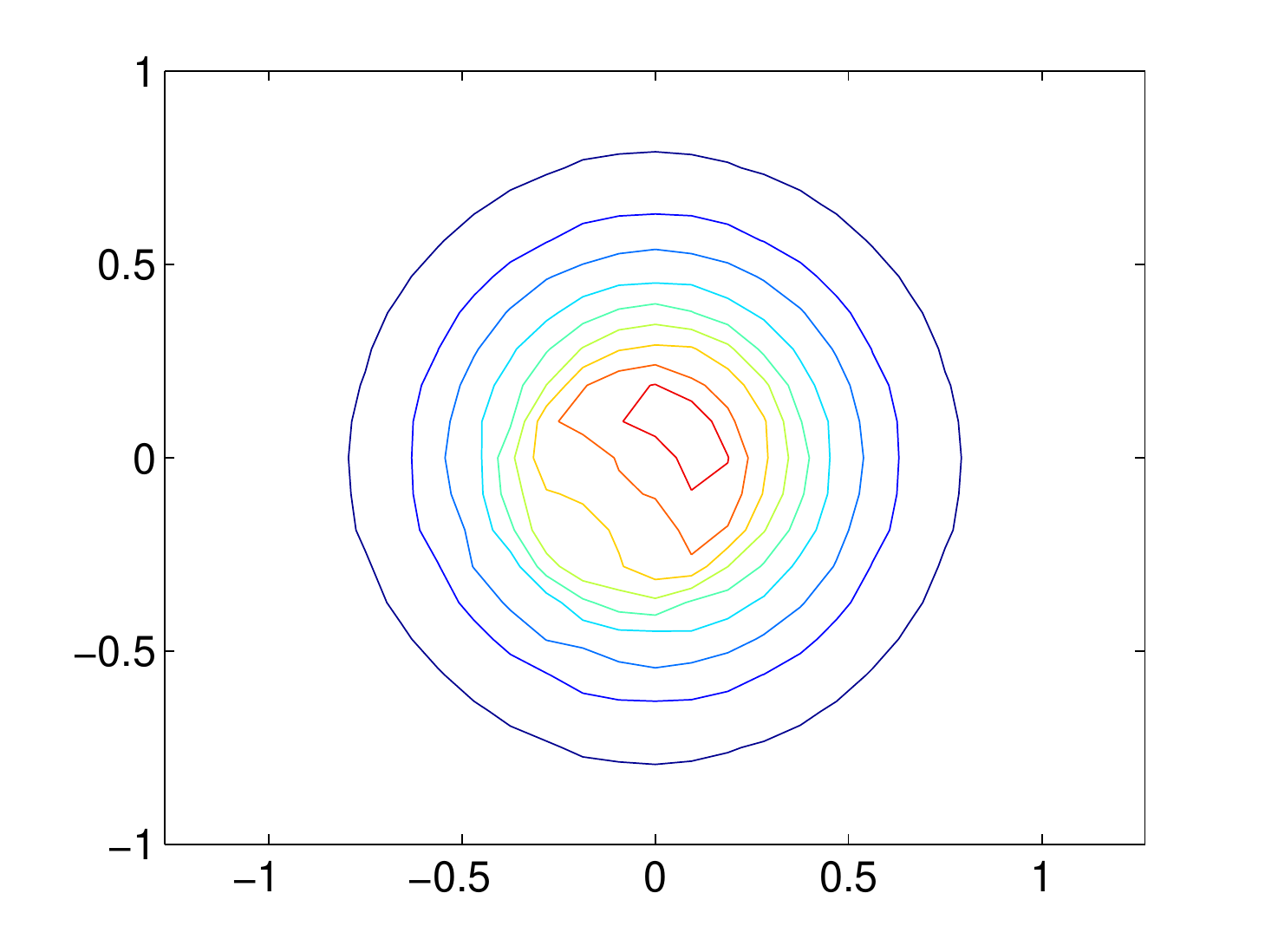}
\includegraphics[width=40mm]{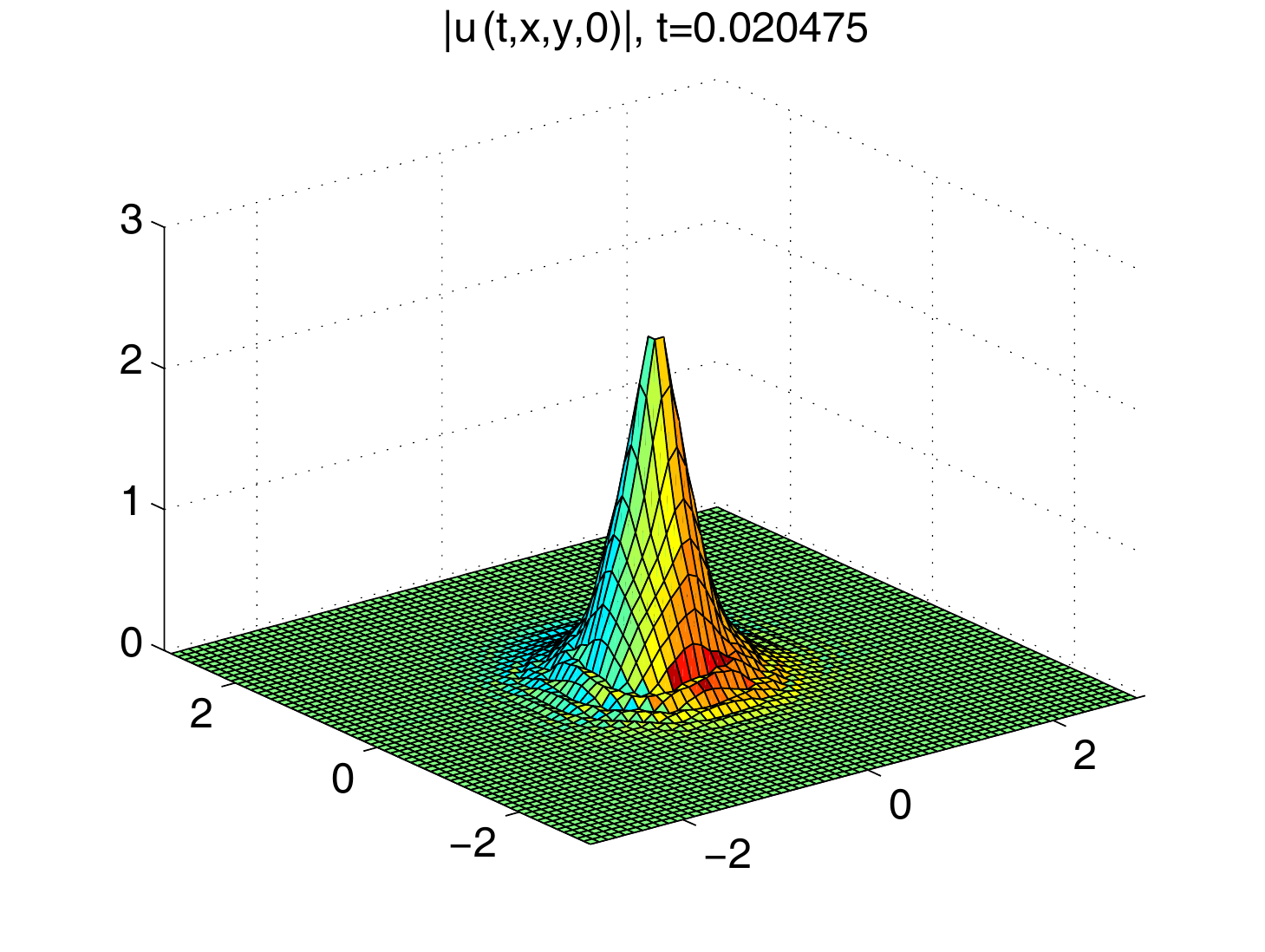}
\includegraphics[width=40mm]{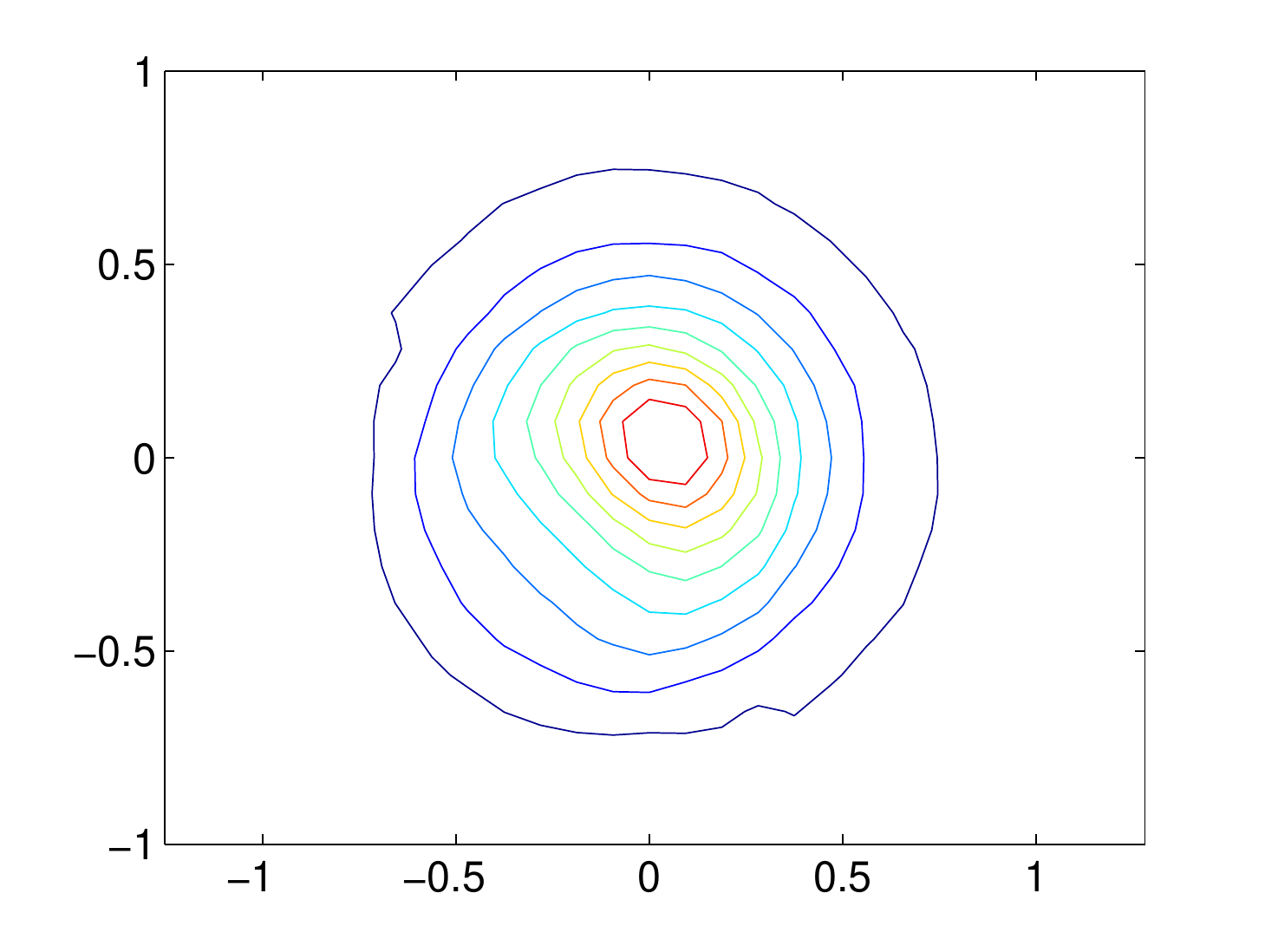}

$|u(t,x,y,0)|_{t=0.00819}|$\hspace{1.4cm} contour plot\hspace{1.4cm} $|u(t,x,y,0)|_{t=0.020475}|$\hspace{1.45cm} contour plot\hspace{0.5cm}\vspace{5mm}

\includegraphics[width=40mm]{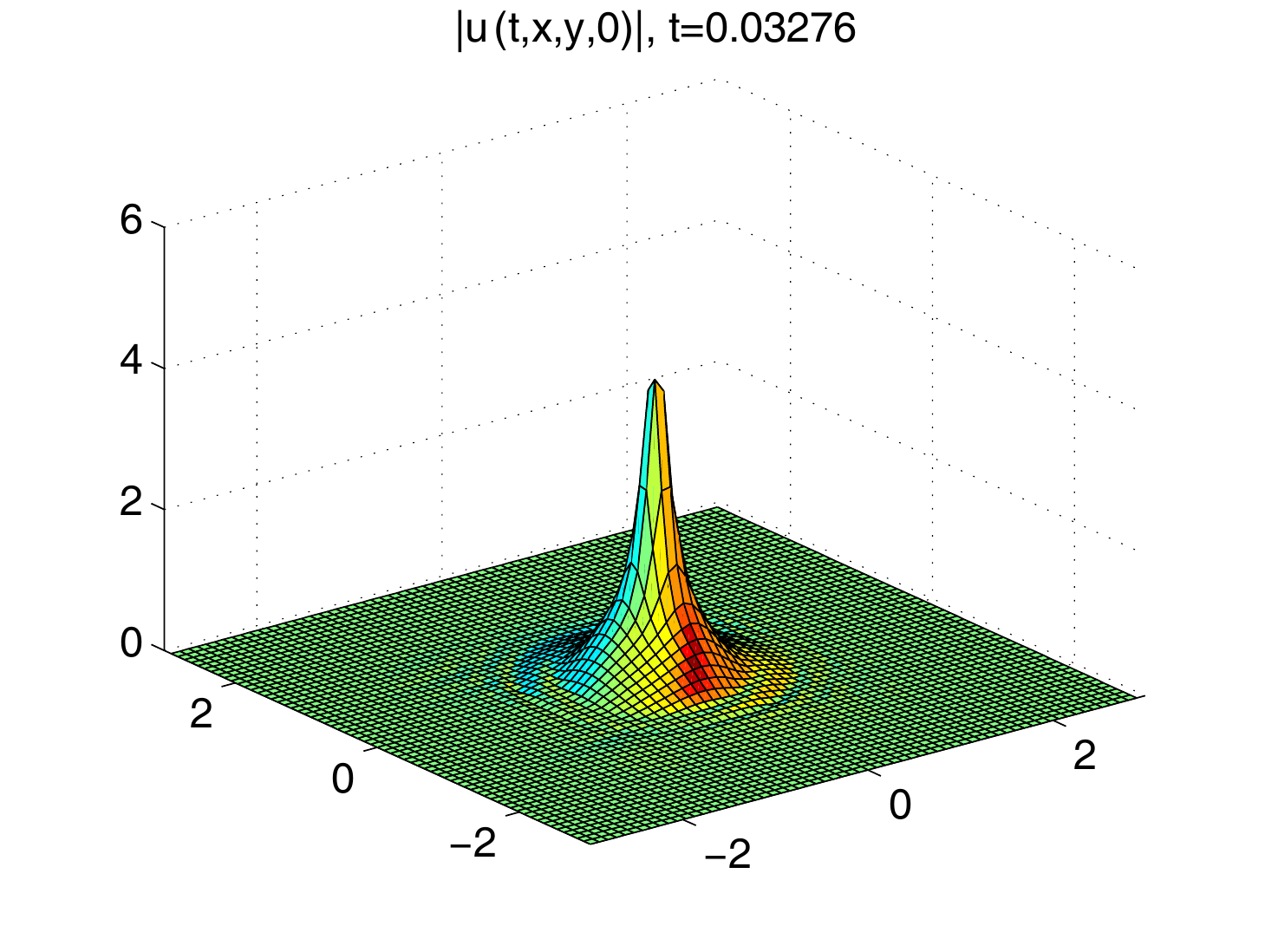}
\includegraphics[width=40mm]{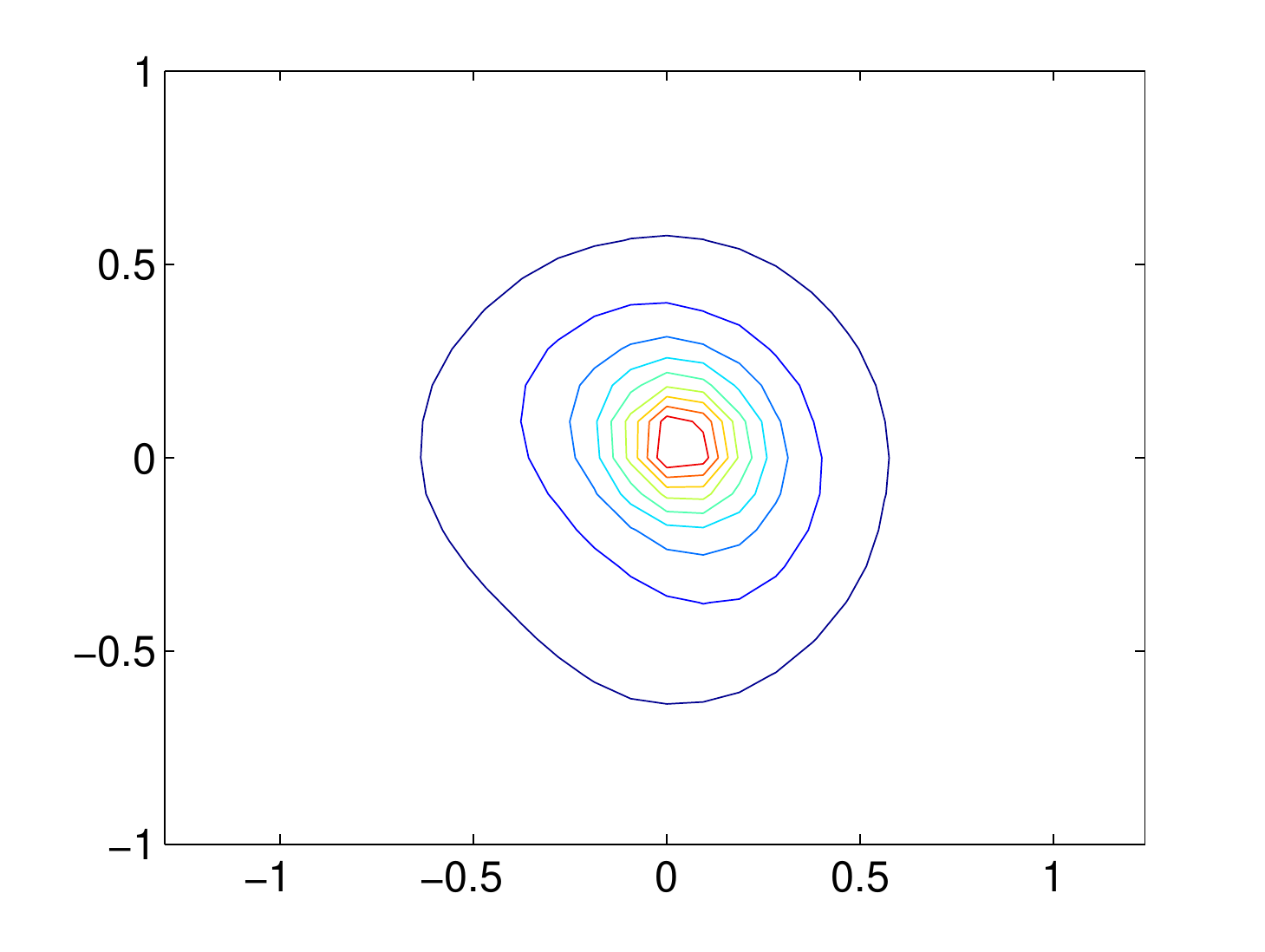}
\includegraphics[width=40mm]{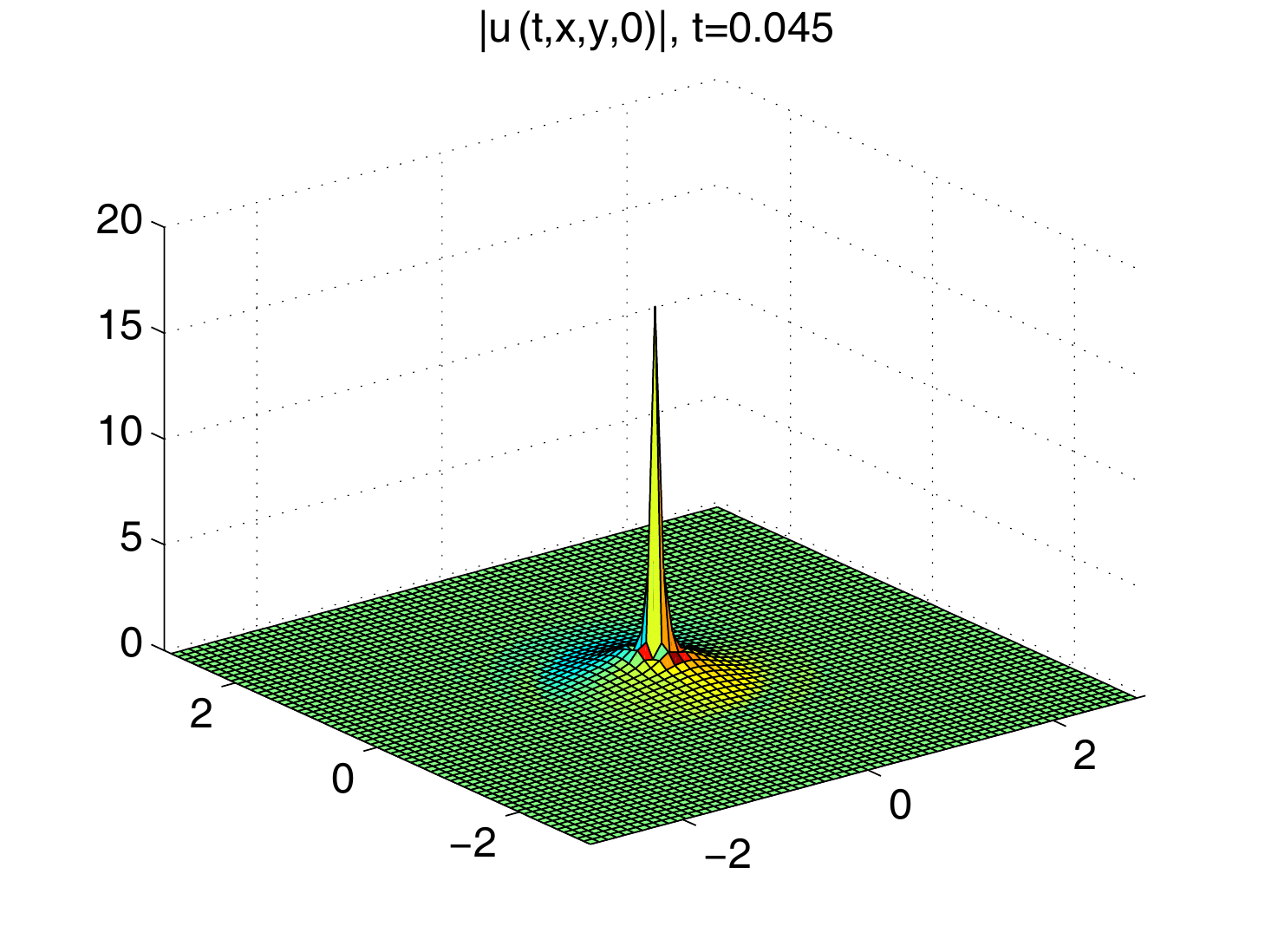}
\includegraphics[width=40mm]{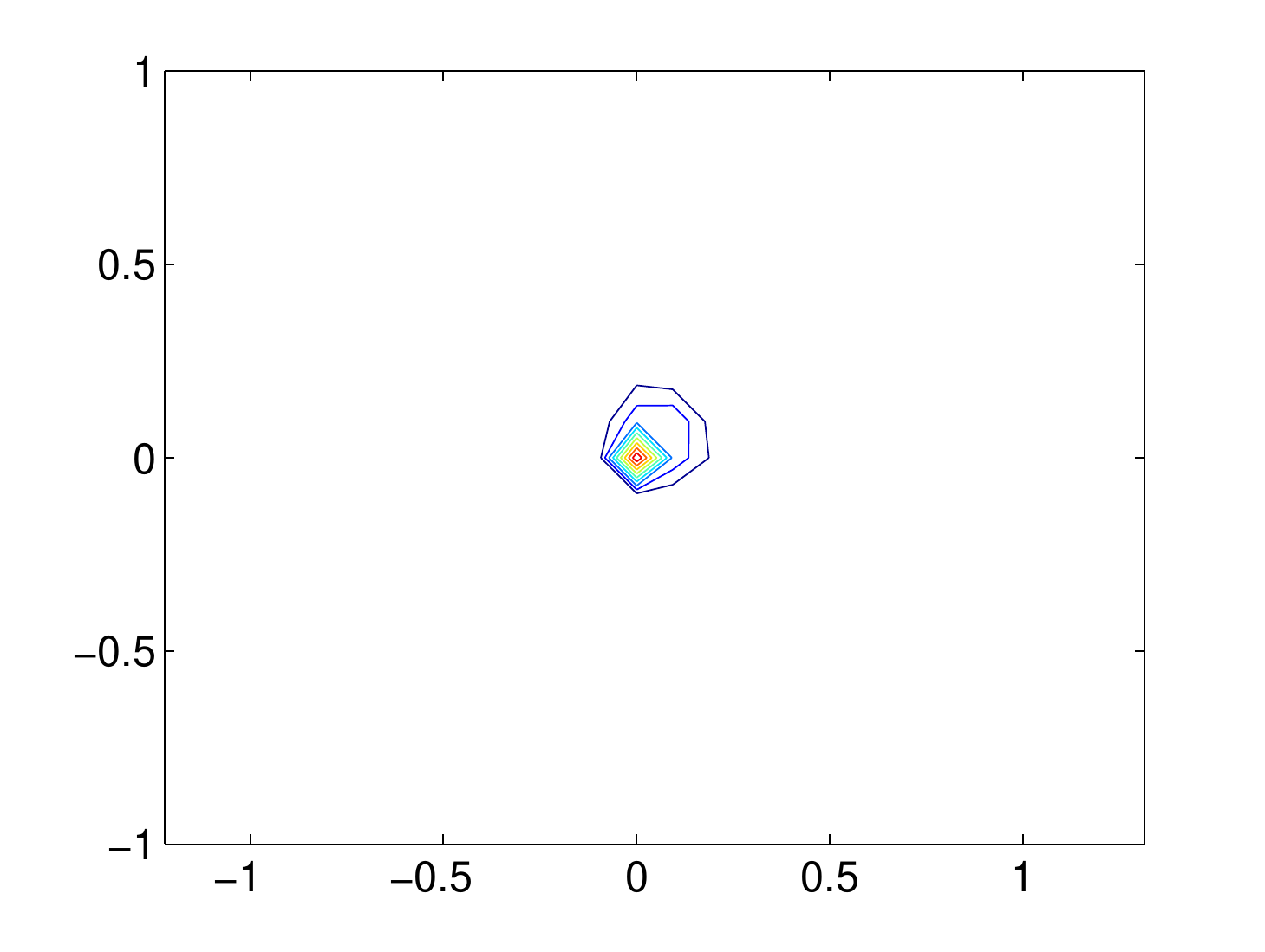}

$|u(t,x,y,0)|_{t=0.03276}|$\hspace{1.45cm} contour plot\hspace{1.5cm} $|u(t,x,y,0)|_{t=0.045}|$\hspace{1.65cm} contour plot\hspace{0.6cm}

\caption{Example \ref{ex2}: Graphs of the fast-potential model, the parameters are:
$\eps=1$, $e(t)=(-1,-1,-1)^T$, $\sigma=1.35$, $a=40$, $C_1=200$, $c=10$, $\omega=5$.}
\label{Fig009-1}
\end{figure}

Furthermore, from Figures \ref{Fig009}--\ref{Fig009-1}, we can see that the Coulomb potential significantly impacts the wave packet at the early stage. This becomes particularly clear with slower blow-up and smaller frequency $\og$.

\subsubsection{Interaction of the magnetic field with the particle}
\begin{example}\label{ex:shot} Here we are interested in the patterns of interaction of the magnetic field with the particle, especially when a particle is shot towards the field, and scattered by it. We set $e=(0,0,1)^T$ in eq. (\ref{eq:fast-p}), and choose an initial datum of the form
\[
u_0(\xb)=e^{-4|\xb-(1,1,0)^T|^2+ i (x_1-x_2)}, \quad \mbox{ for } \xb=(x_1,x_2,x_3)^T.
\]
In this example, we add a harmonic trap potential
\[V(\xb)=50|\xb|^2,\]
\ie the equation in \eqref{eq:fast} becomes
\begin{equation}
i\d_tu^\omega=-\Delta u^\omega+V\(\xb-b(t)\)u^\omega+V^\omega u^\omega
+C_1(|\cdot|\ast|u^\omega|^2)u^\omega-a|u^\omega|^\sigma u^\omega
\end{equation}
\end{example}

\begin{figure}[h!t]
\centering
\includegraphics[width=50mm]{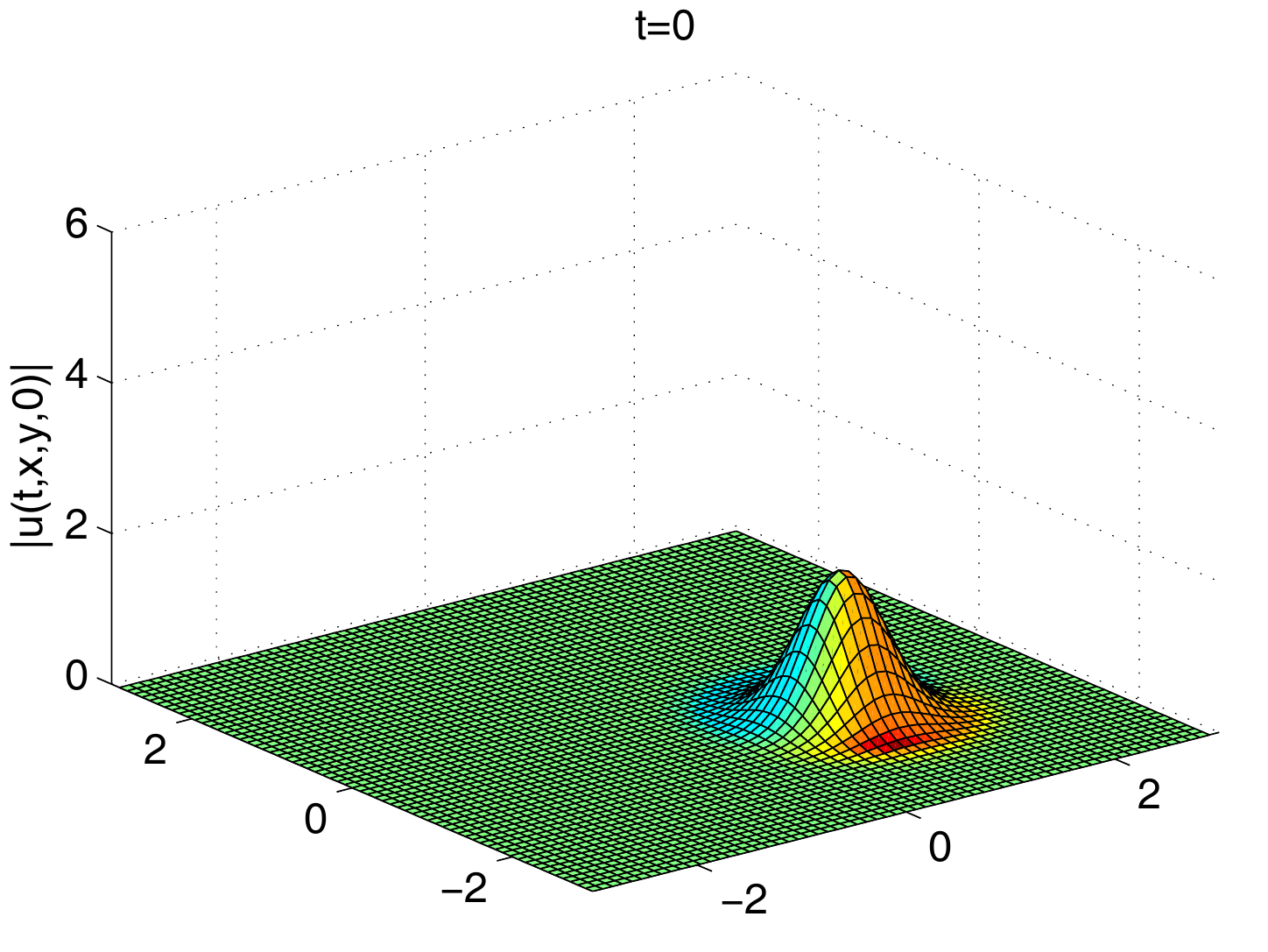}
\includegraphics[width=50mm]{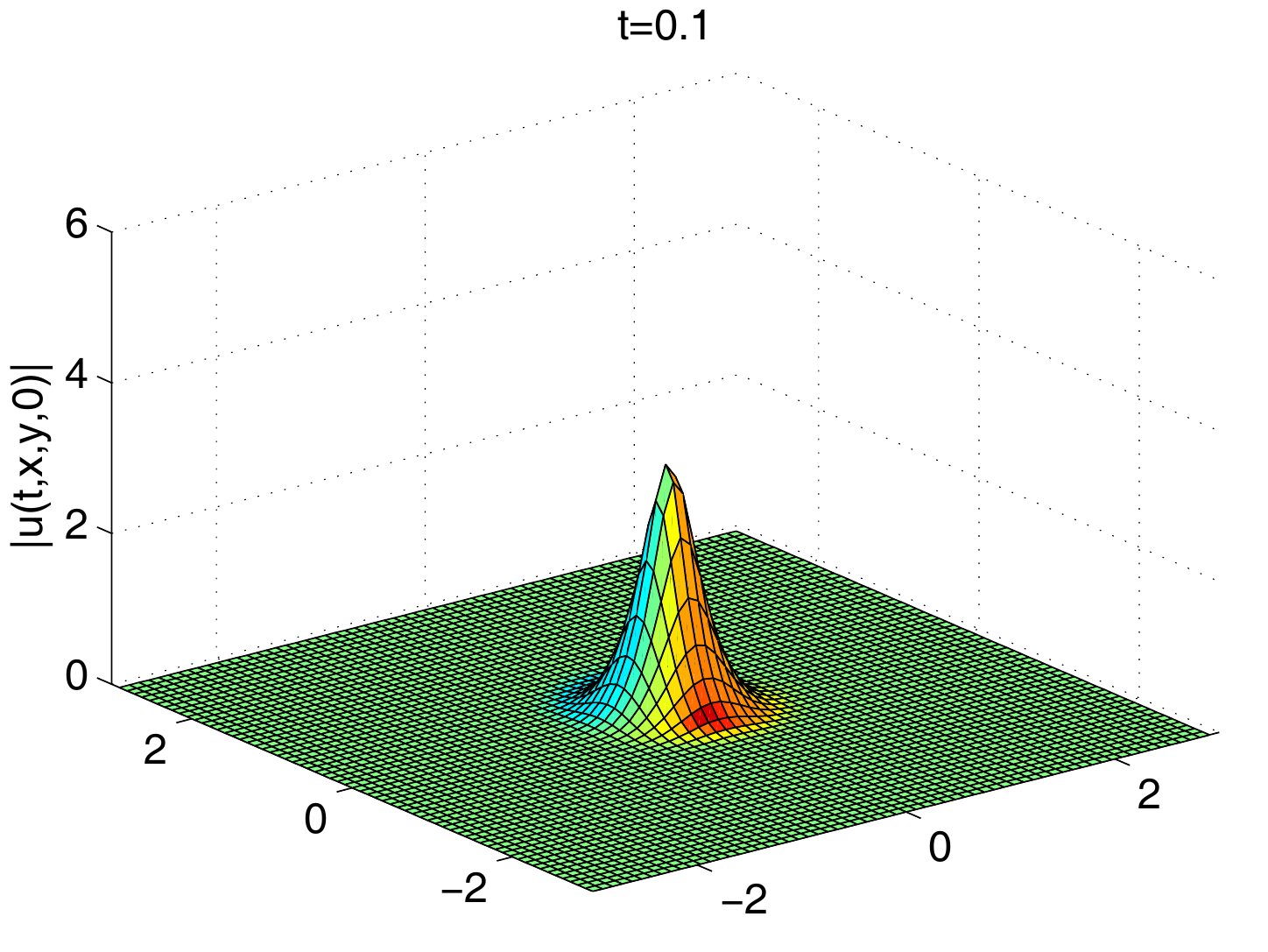}
\includegraphics[width=50mm]{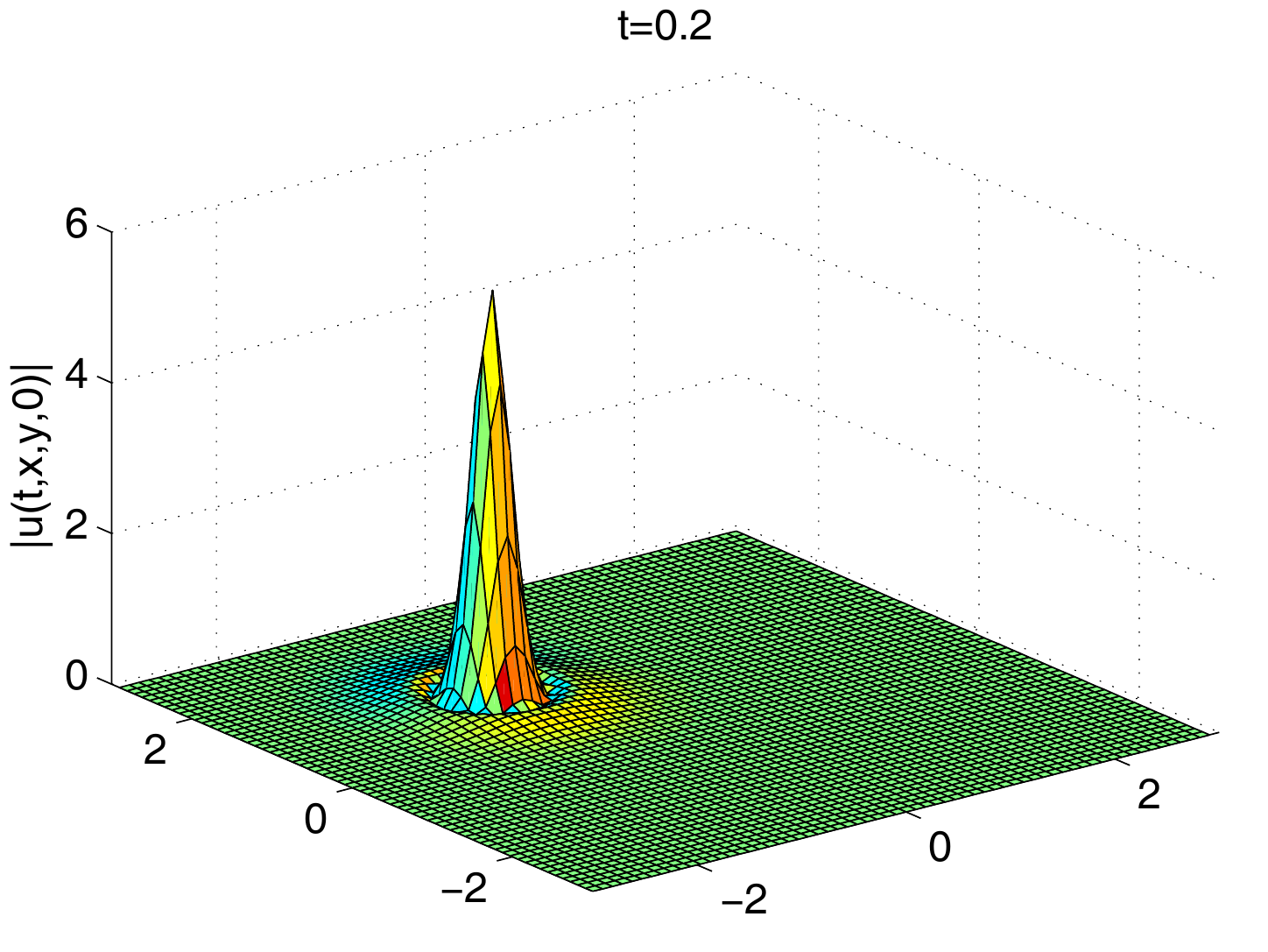}\vspace{1mm}

\includegraphics[width=50mm]{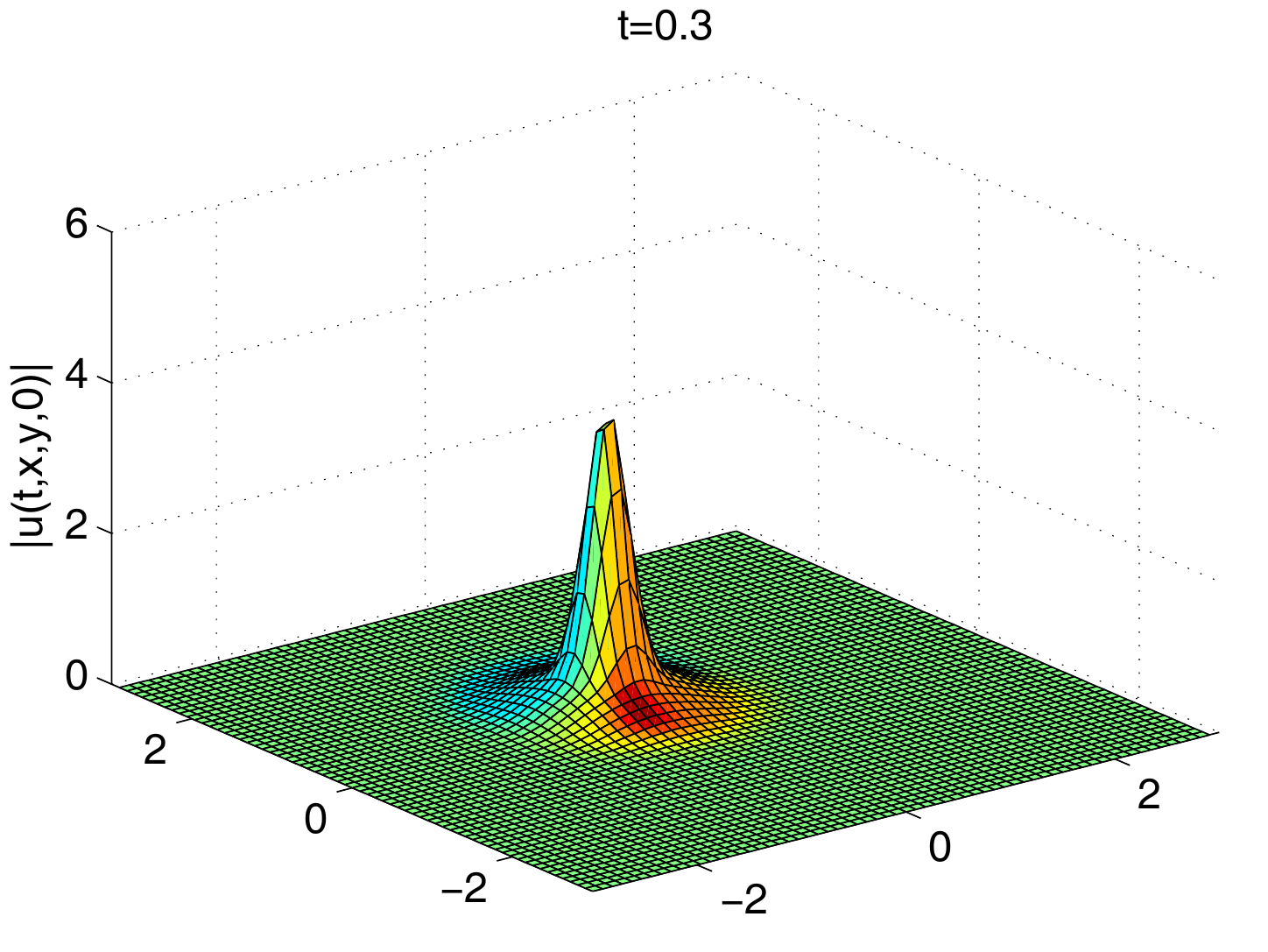}
\includegraphics[width=50mm]{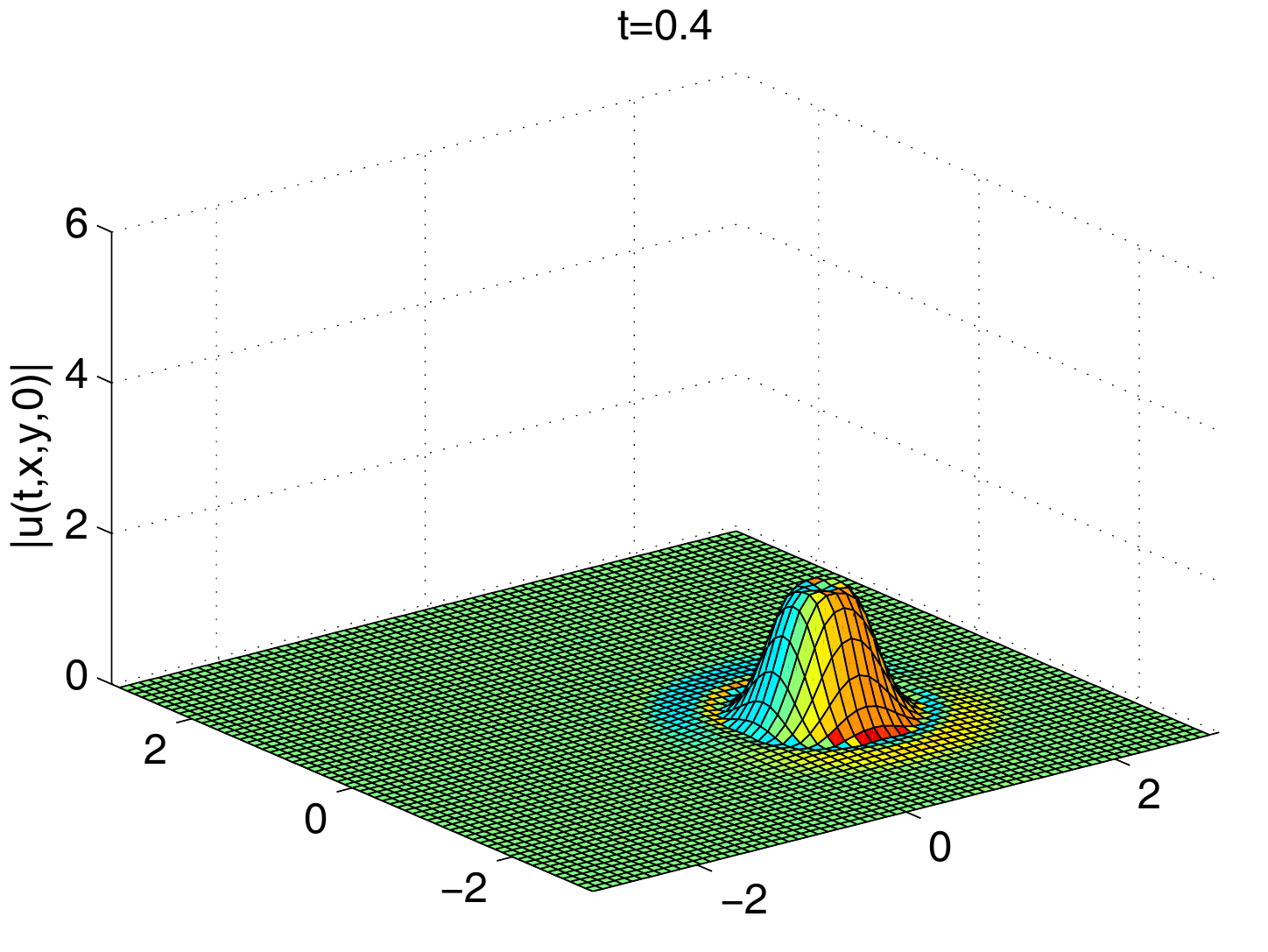}
\includegraphics[width=50mm]{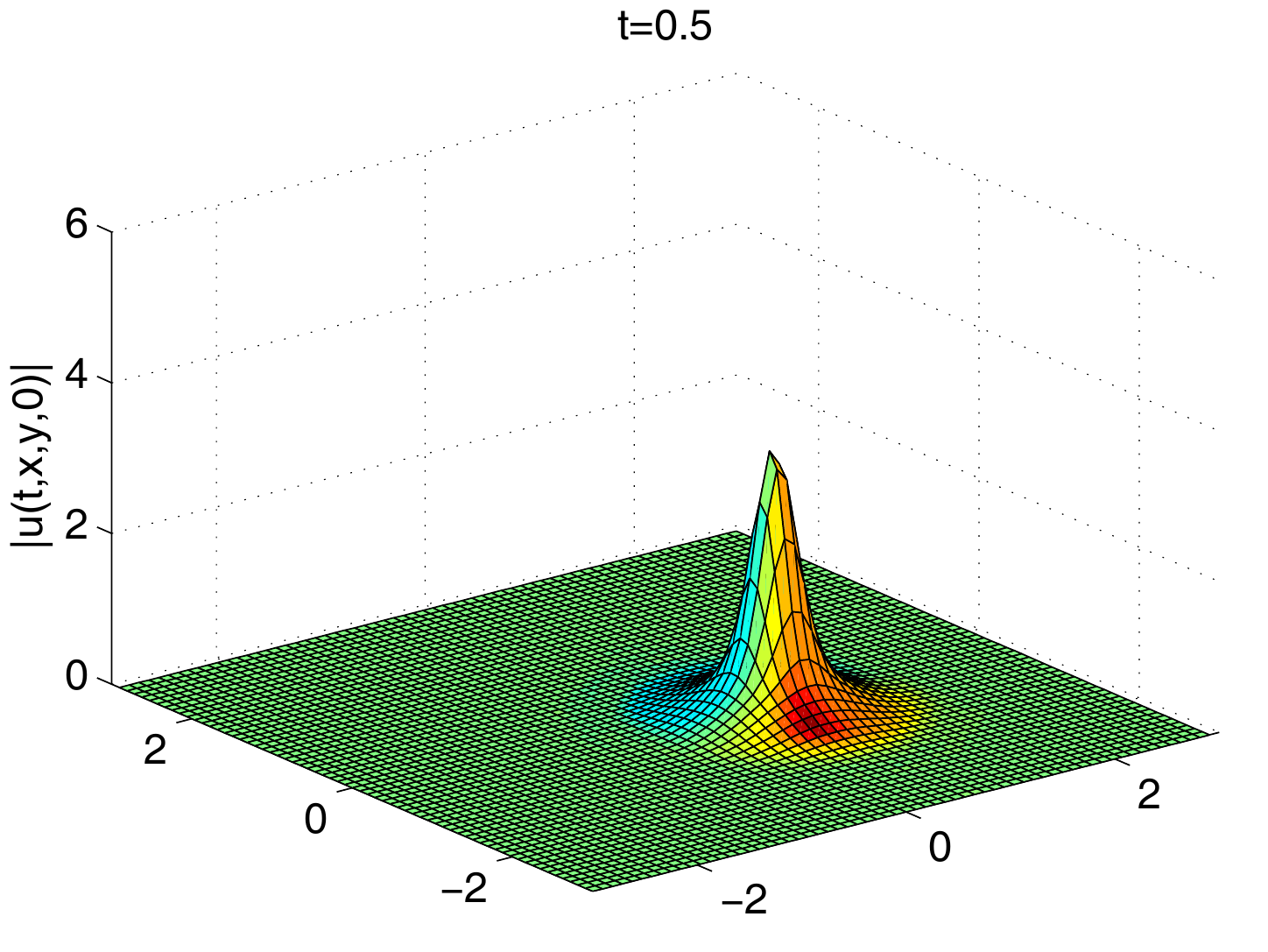}\vspace{1mm}

\includegraphics[width=50mm]{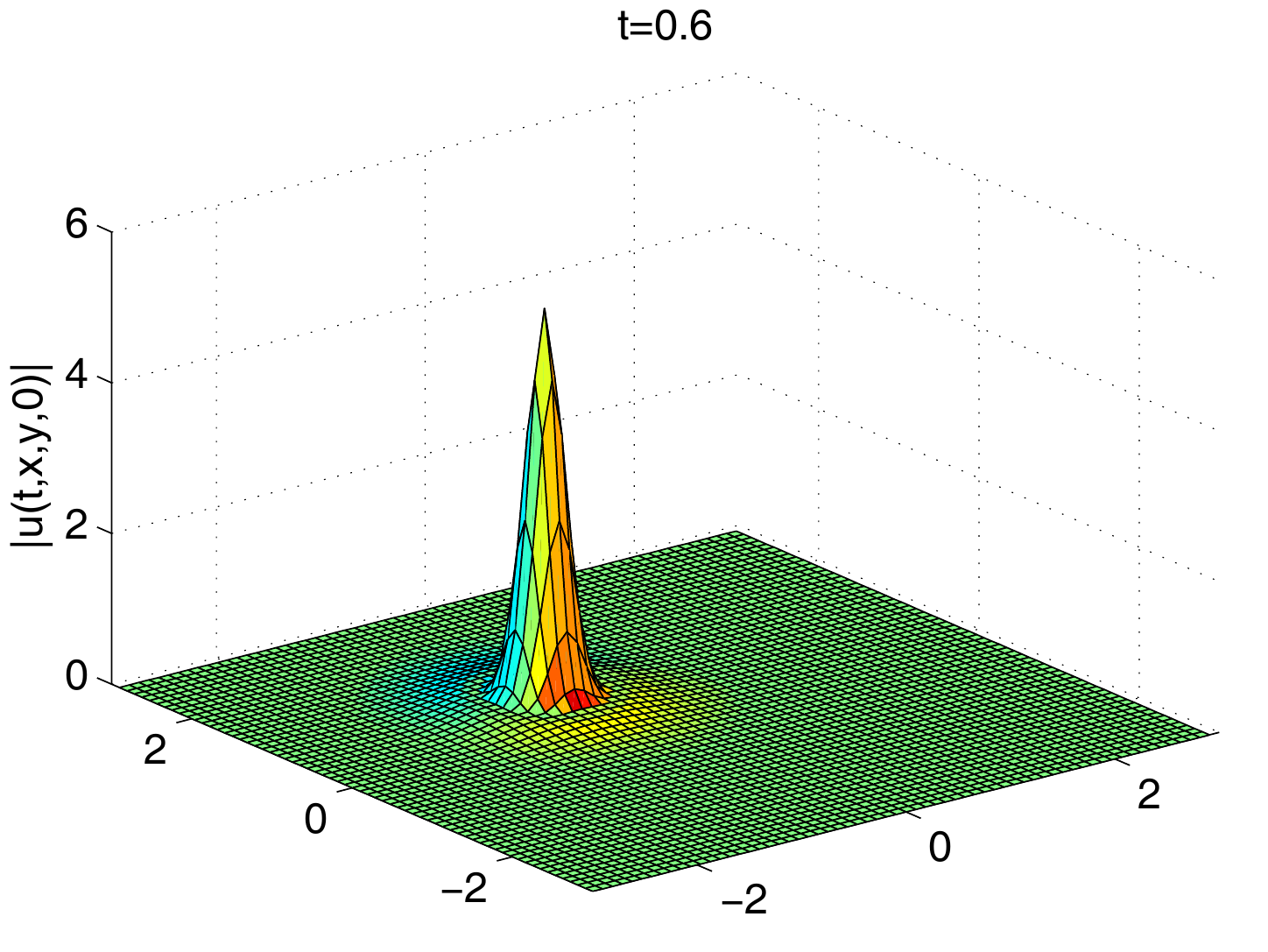}
\includegraphics[width=50mm]{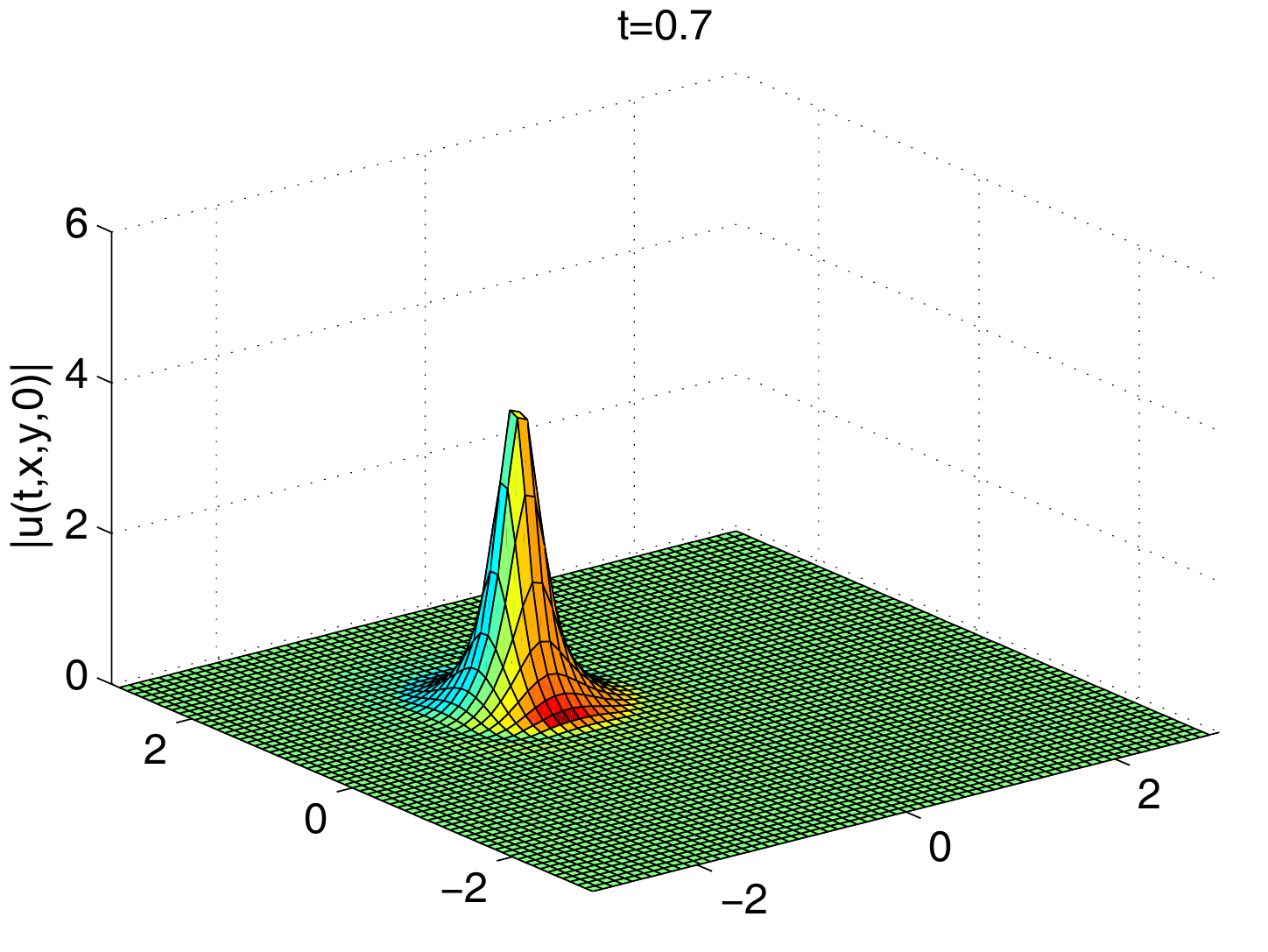}
\includegraphics[width=50mm]{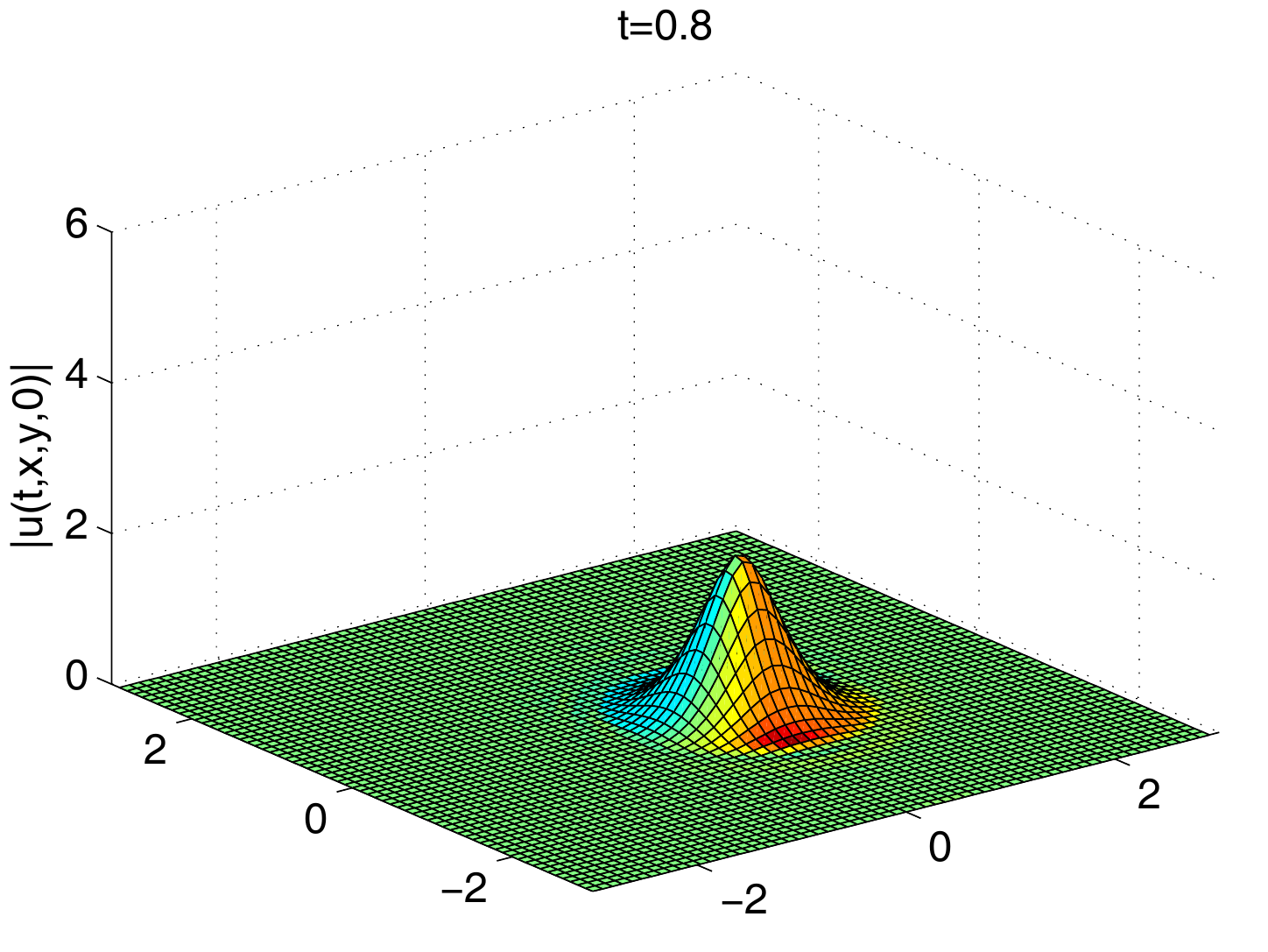}\vspace{1mm}

\includegraphics[width=50mm]{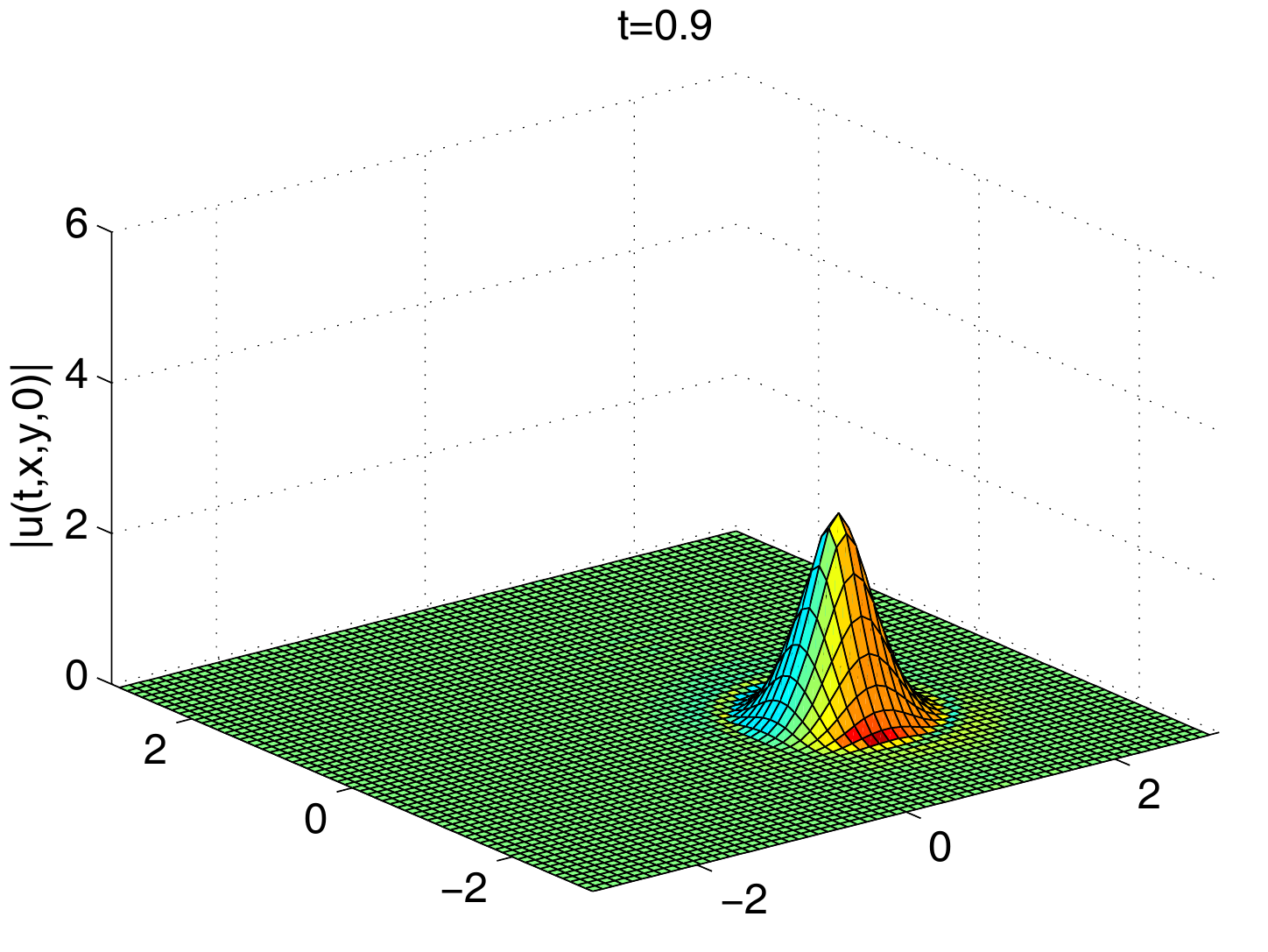}
\includegraphics[width=50mm]{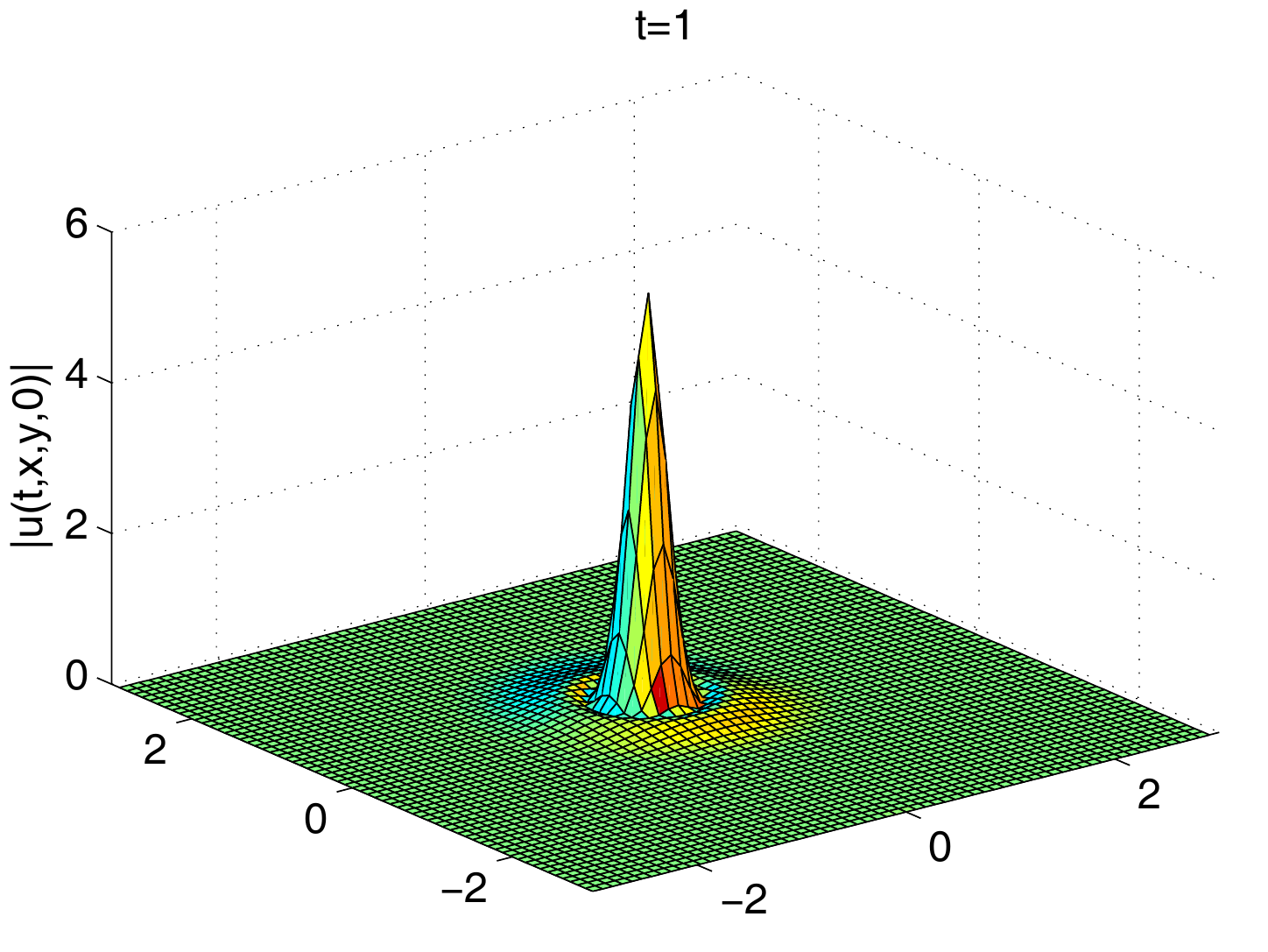}
\includegraphics[width=50mm]{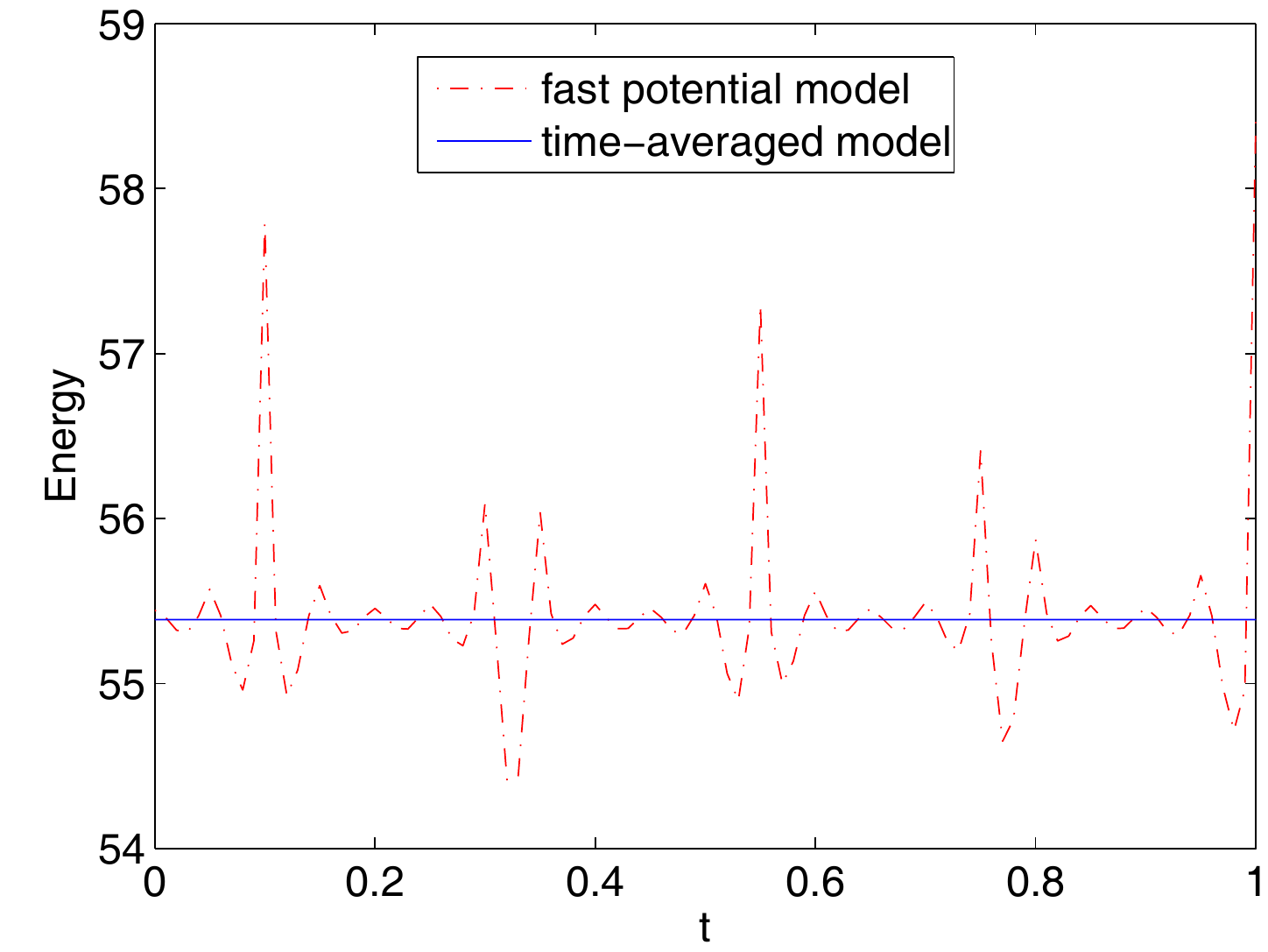}
\caption{Example \ref{ex:shot}: (fast-potential model), the parameters are:
$\eps=1$, $e(t)=(0,0,1)^T$, $\sigma=\frac{2}{3}$, $a=50$, $C_1=20$, $c=1$, $\omega=10^2$.}
\label{Fig010}
\end{figure}

The results are shown in Figure \ref{Fig010}. We can see that the wave packet moves under the interaction with the trap $V(\xb)$.

\subsubsection{Time-dependent vectors $e(t)$}
\begin{example}\label{ex:td} We consider a time-dependent vector $e(t)=(-1,-1,-1)\sin 2\pi t$ in eq. (\ref{eq:fast-p}), and an initial datum of the form
\[
u_0(x)=e^{-4|\xb|^2}.
\]
\end{example}

\begin{figure}[h!t]
\centering
\includegraphics[width=40mm]{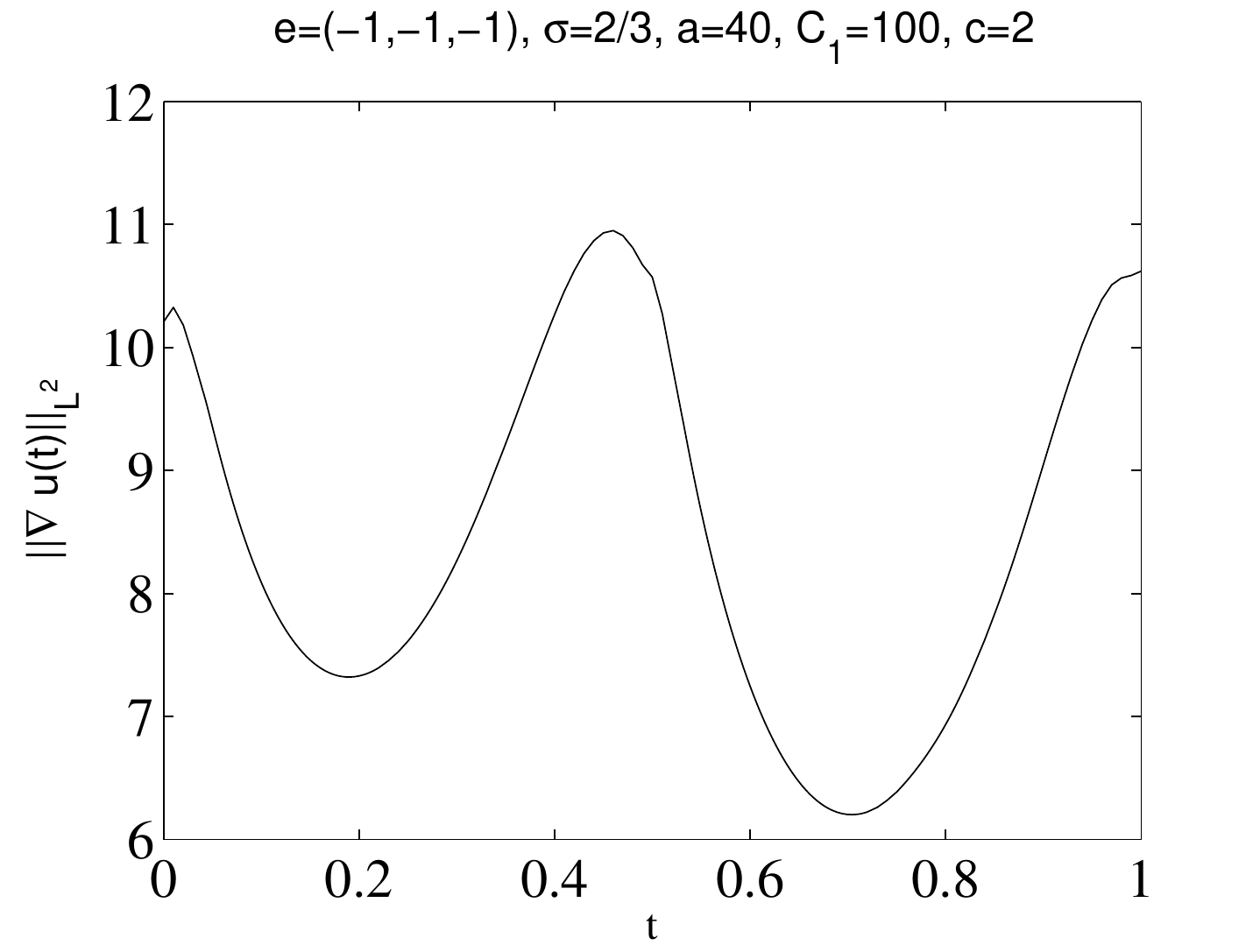}
\includegraphics[width=40mm]{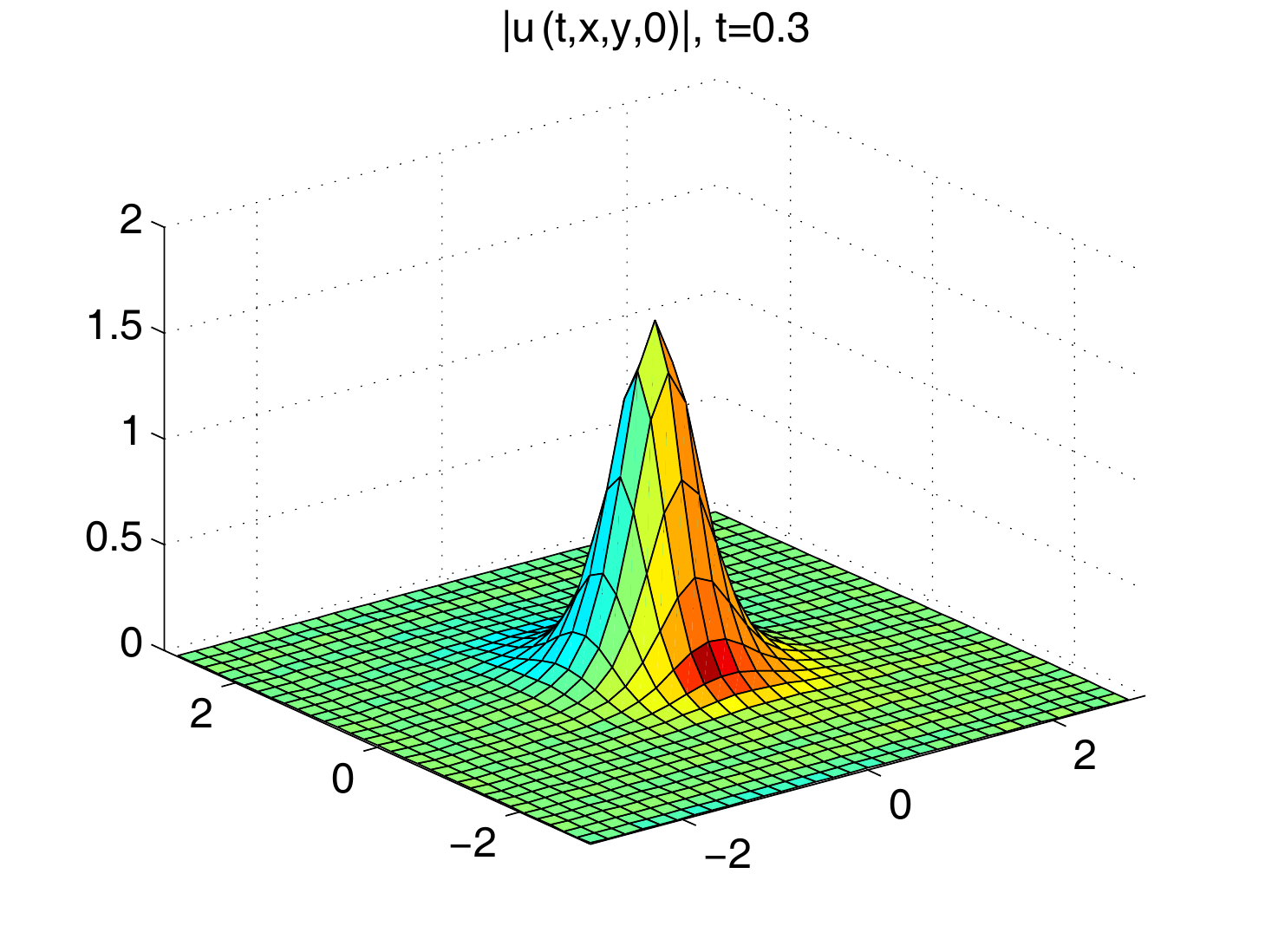}
\includegraphics[width=40mm]{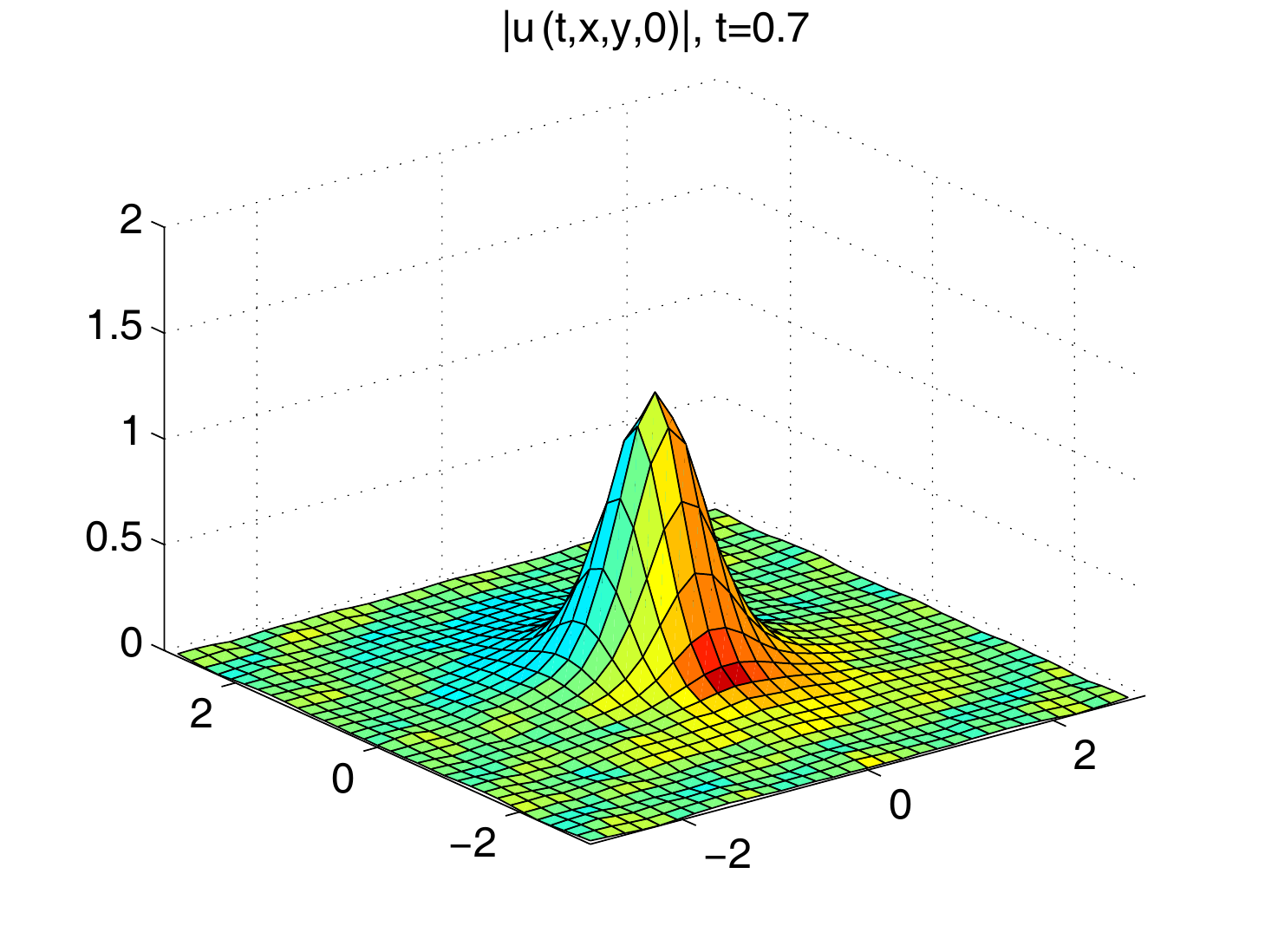}
\includegraphics[width=40mm]{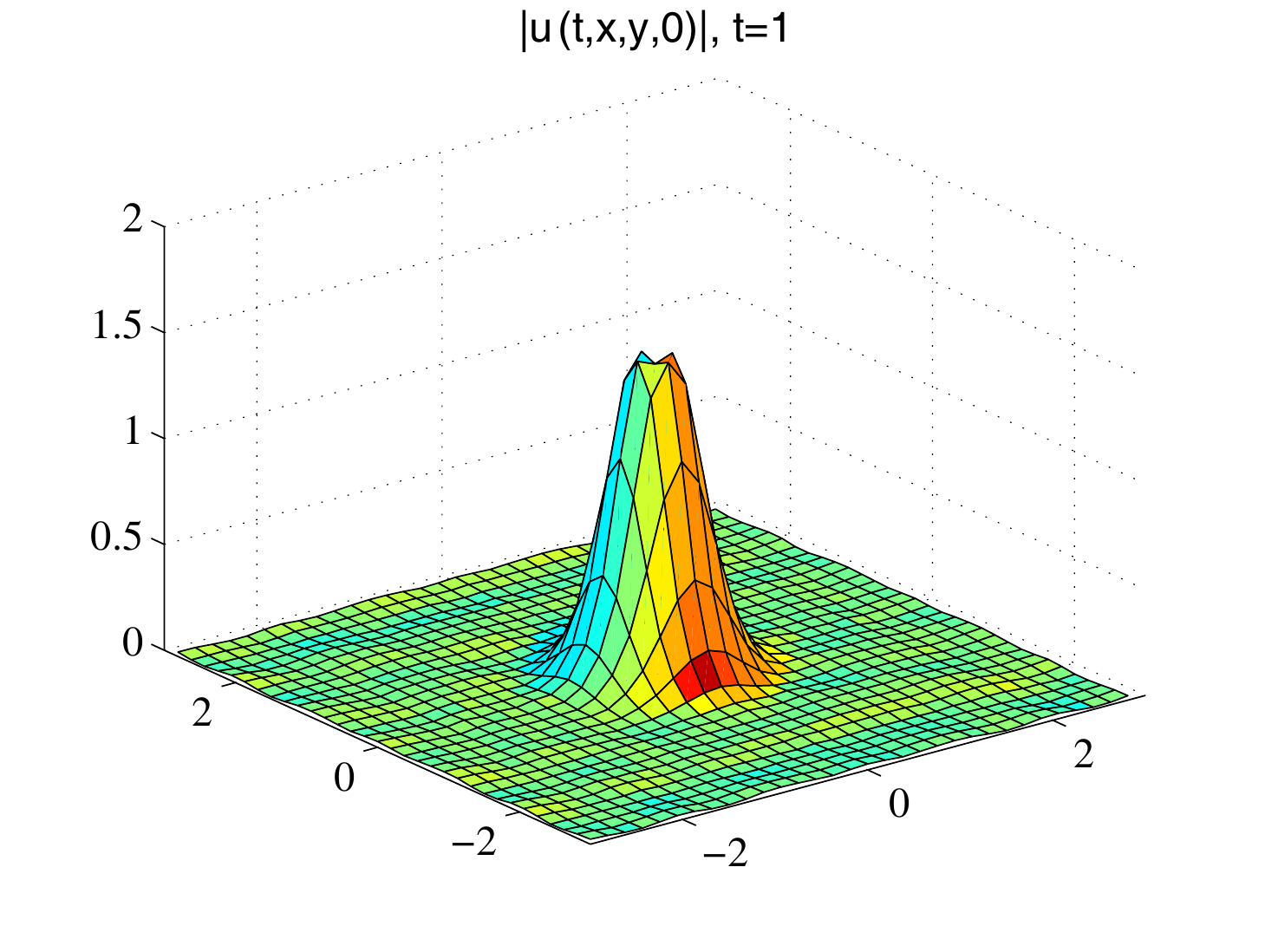} \vspace{5mm}

\includegraphics[width=40mm]{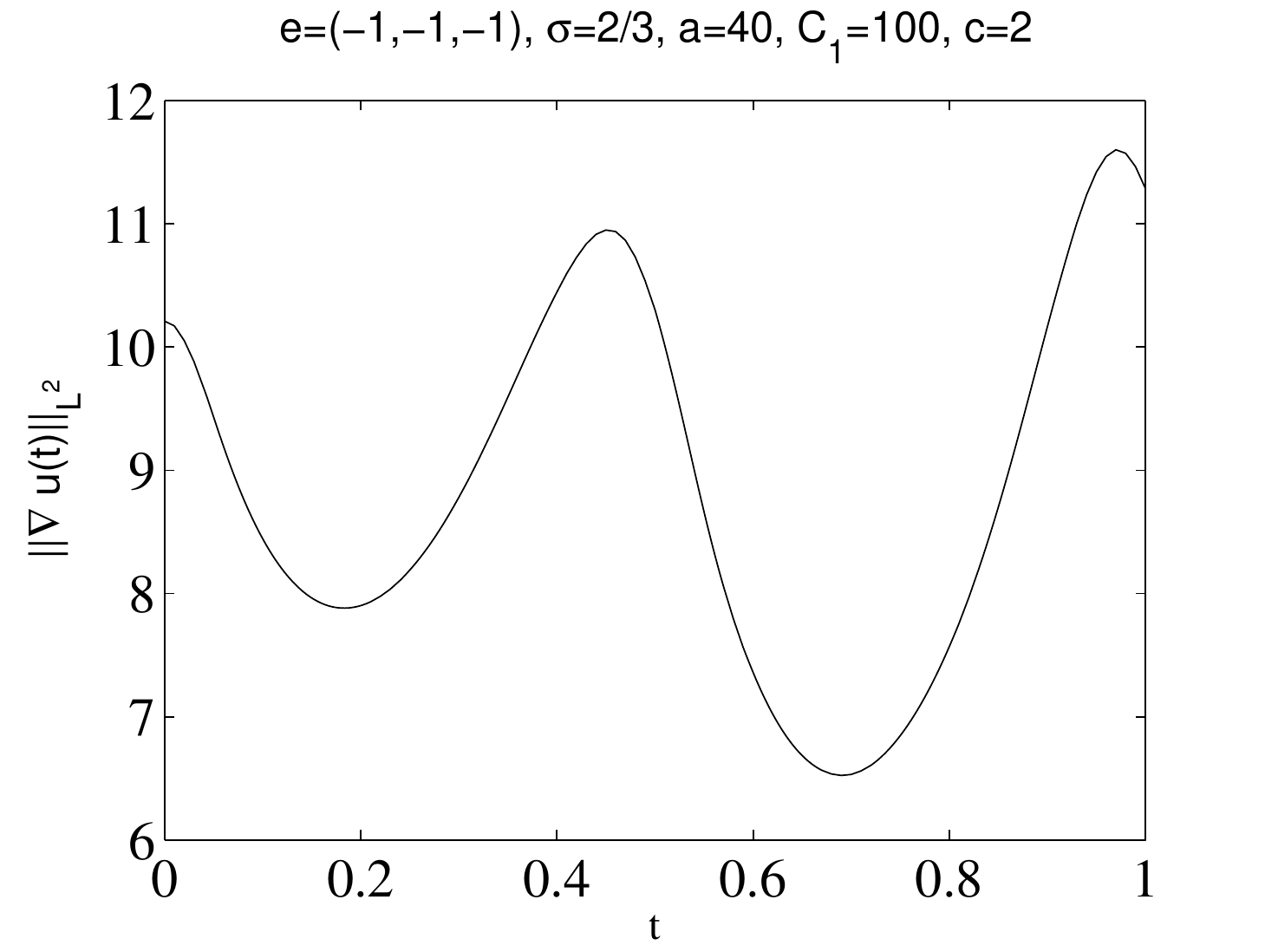}
\includegraphics[width=40mm]{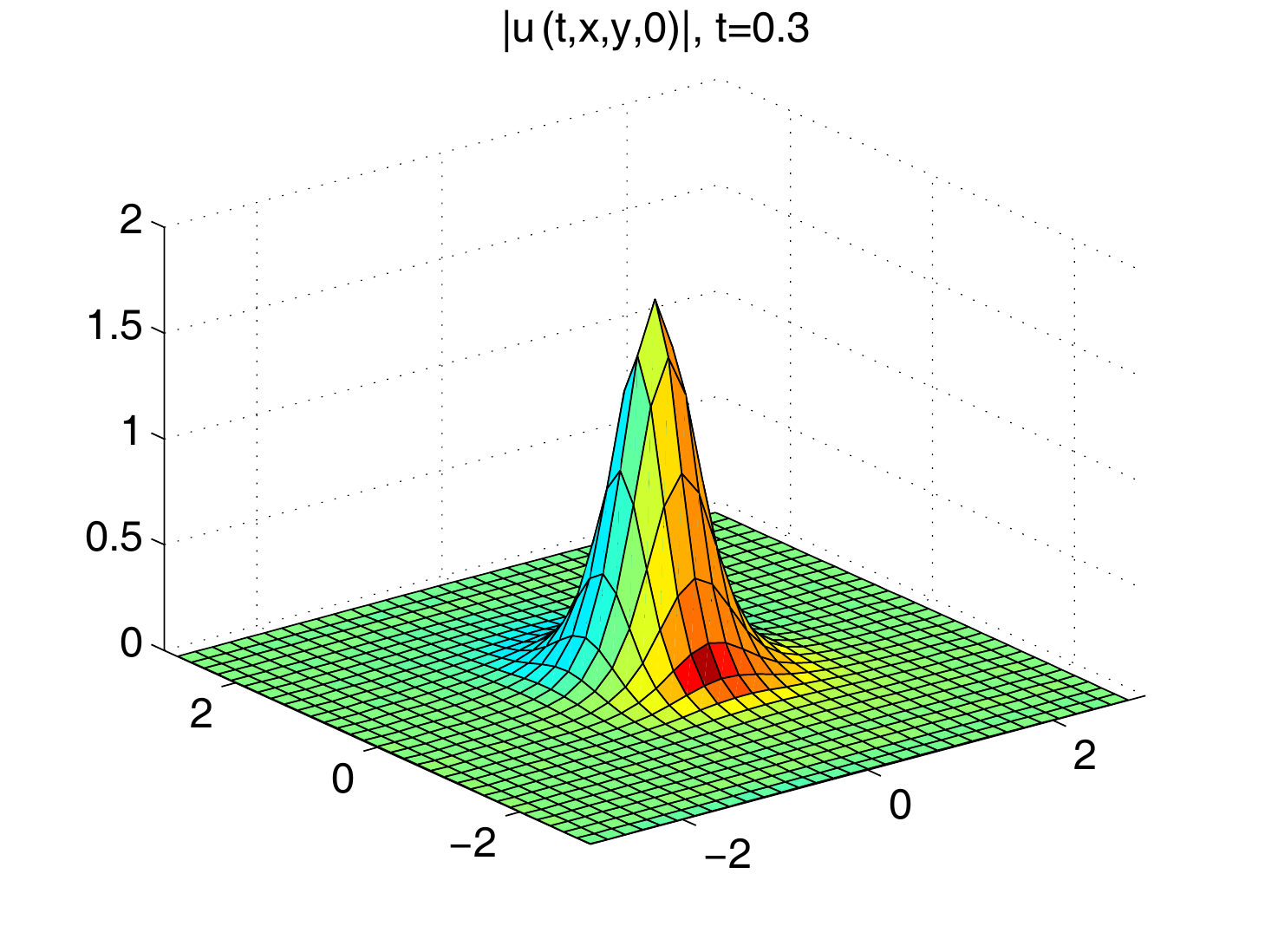}
\includegraphics[width=40mm]{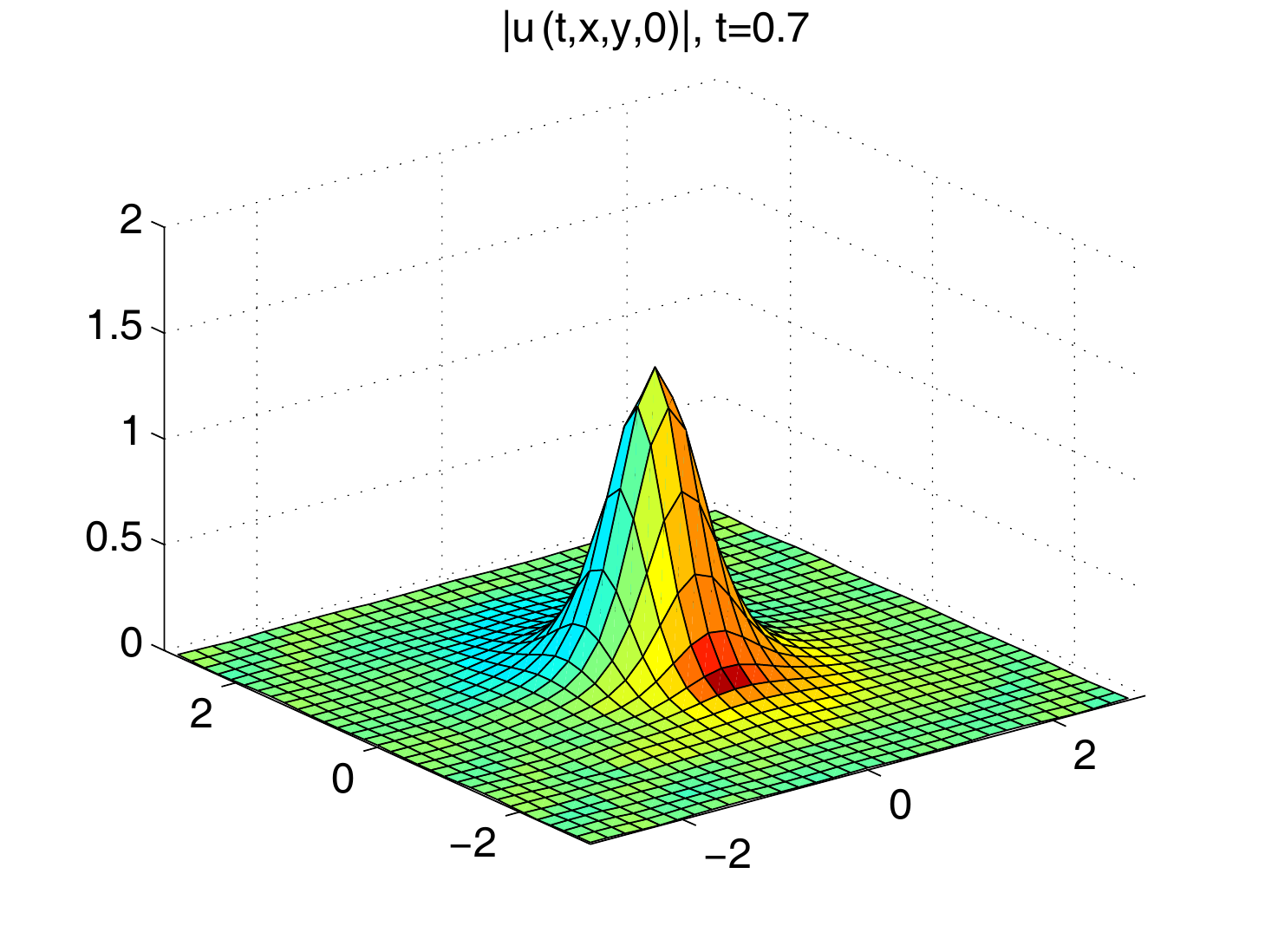}
\includegraphics[width=40mm]{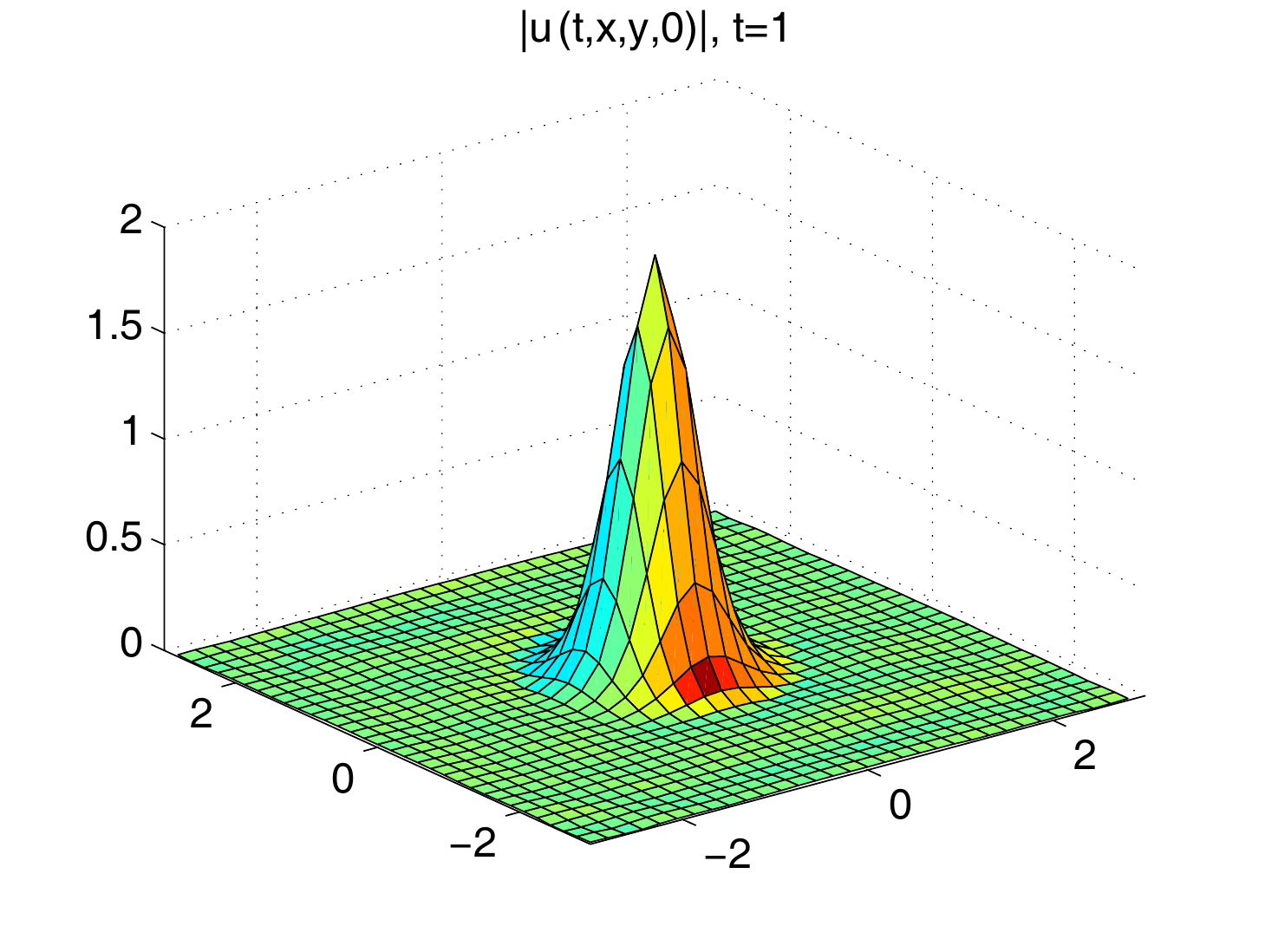}
\caption{Example \ref{ex:td}: Graphs of the example \ref{ex:td} for fast-potential model (first line) and time-averaged model (second line), the parameters are:
$\eps=1$, $e(t)=(-1,-1,-1)^T\sin 2\pi t$, $\sigma=\frac{2}{3}$, $a=40$, $C_1=100$, $c=2$, $\omega=10^2$.}
\label{Fig011}
\end{figure}

The results are shown in Figure \ref{Fig011}. As $e(t)$ is slowly varying in time, its effects are not very pronounced except when $e(t)$ is close to 0.

\subsubsection{Interaction with periodic lattice}
\begin{example}\label{ex:pl} Finally, we consider a case with periodic lattice,
\ie the equation in \eqref{eq:fast} becomes
\begin{equation}
i\d_tu^\omega=-\Delta u^\omega+V^\omega u^\omega+V_\Gm\(\xb-b(t)\)u^\omega
+C_1(|\cdot|^{-1}\ast|u^\omega|^2)u^\omega-a|u^\omega|^\sigma u^\omega,
\end{equation}
and its time-averaged problem is
\begin{equation}
i\d_tu=-\Delta u+\langle V\rangle u+\langle V_\Gamma\rangle u+C_1(|\cdot|^{-1}\ast|u|^2)u-a|u|^\sigma u.
\end{equation}
We also take the initial datum of the form
\[
u_0(x)=e^{-4|\xb|^2}.
\]
Here $V_\Gm(\yb)$ is periodic w.r.t to a regular lattice $\Gm$,
\[V_\Gm(\yb+\gm)=V_\Gm(\yb), \quad \fa \gm\in\Gm,\ \yb\in \R^3.\]
For sake of simplicity, we choose \[V_\Gm(\yb)=\sum_{l=1}^3 \sin^2(\og_l y_l), \quad \mbox{ for }\yb=(y_1,y_2,y_3)^T\in\R^3, \]
with $\og_1=\og_2=\og_3=2\pi$.
\end{example}

\begin{figure}[h!t]
\centering
\includegraphics[width=50mm]{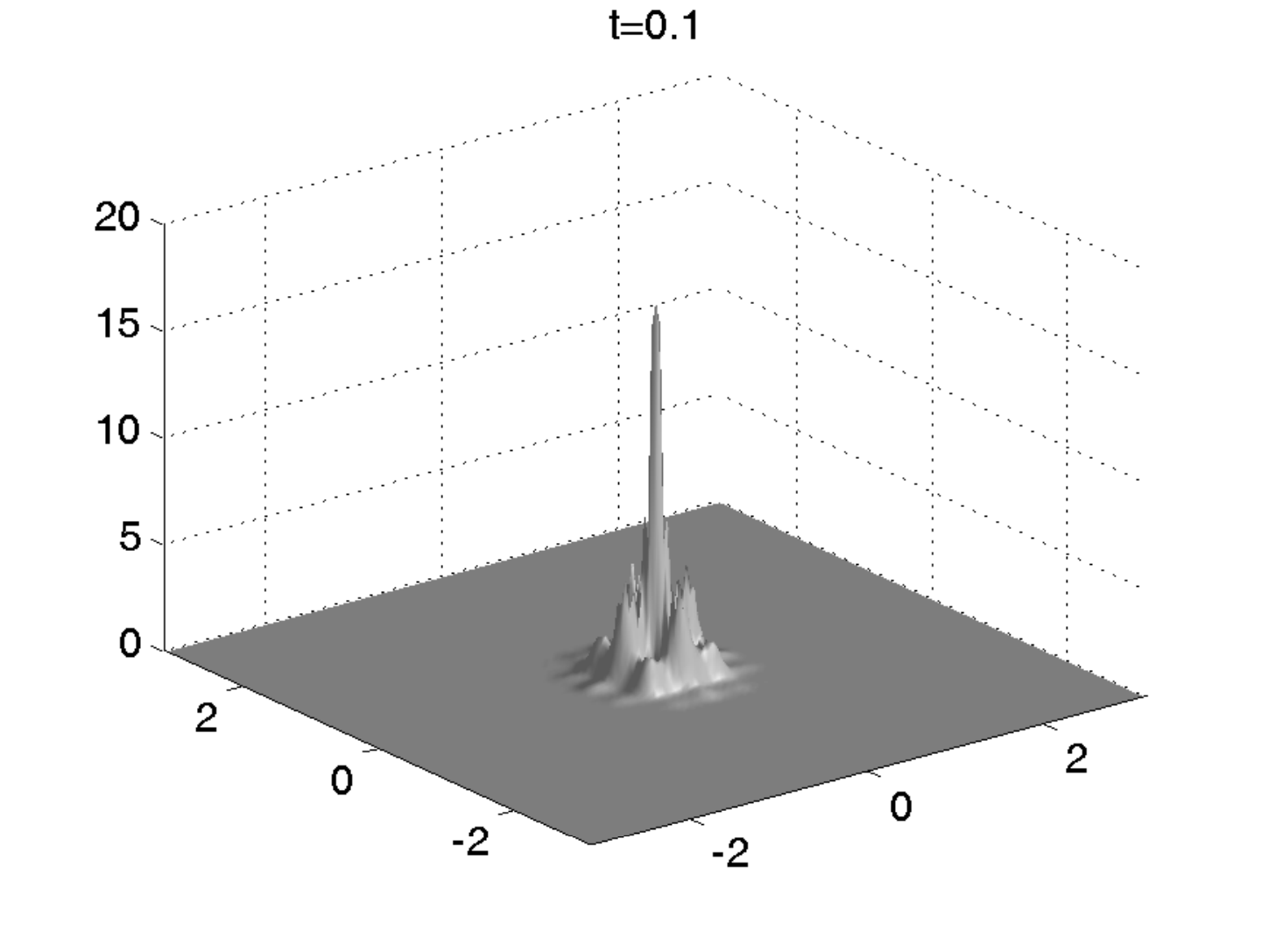}
\includegraphics[width=50mm]{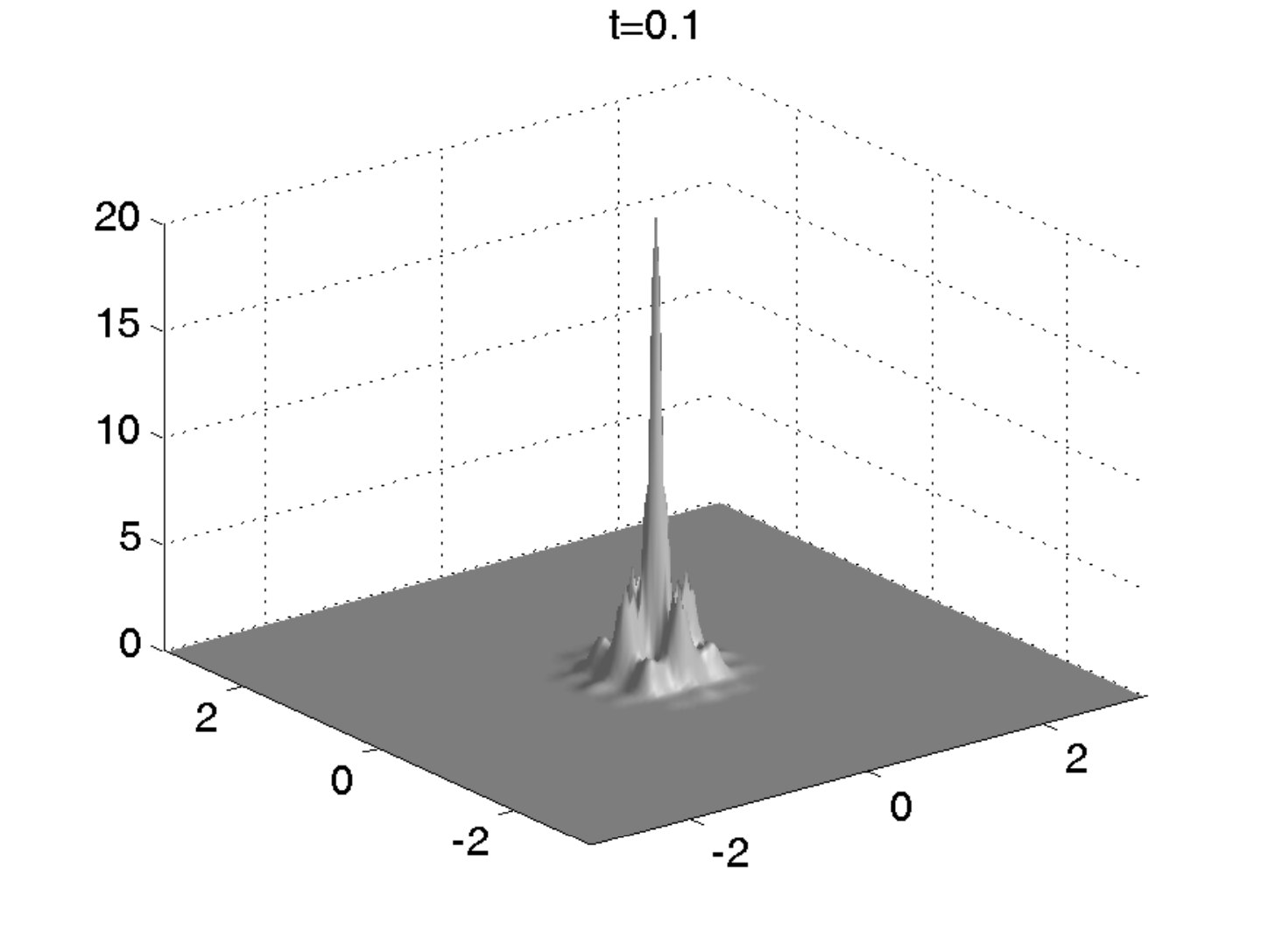}
\includegraphics[width=50mm]{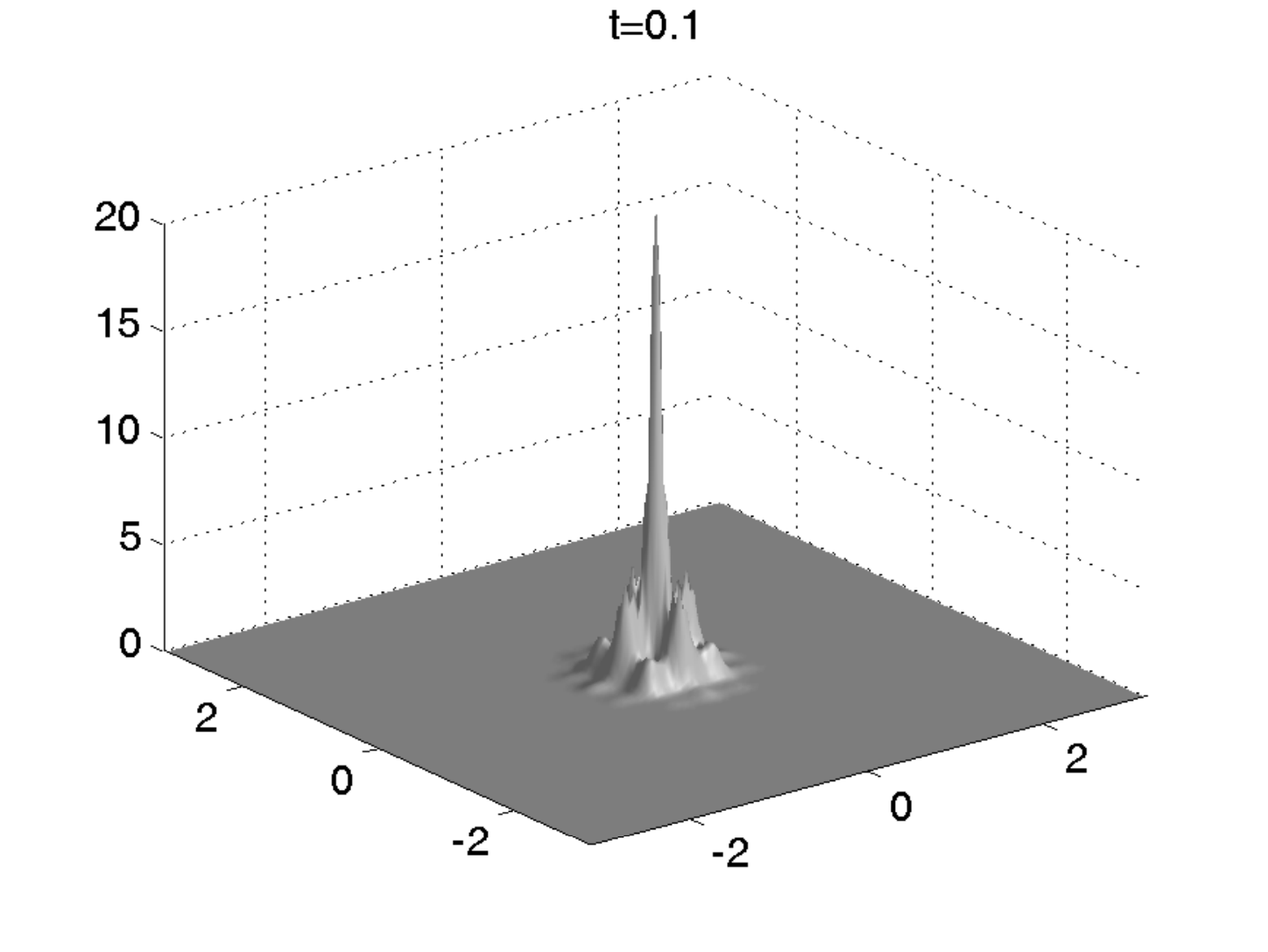}
\caption{Graphs of $|u^\og(t,x,y,0)|^2$ and $|u(t,x,y,0)|^2$ of the example \ref{ex:pl} for fast-potential model (left $\og=80$, middle $\og=160$)
and time-averaged model (right) by our Bloch-decomposition based algorithm given in section \ref{sec:2.2}, the parameters are:
$\eps=\frac{1}{8}$, $e(t)=(0,0,1)^T$, $\sigma=\frac{2}{3}$, $a=1$, $C_1=1$, $c=2$, mesh size $h=\frac{1}{128}$.}
\label{Fig012}
\end{figure}

\begin{figure}[h!t]
\centering
\includegraphics[width=50mm]{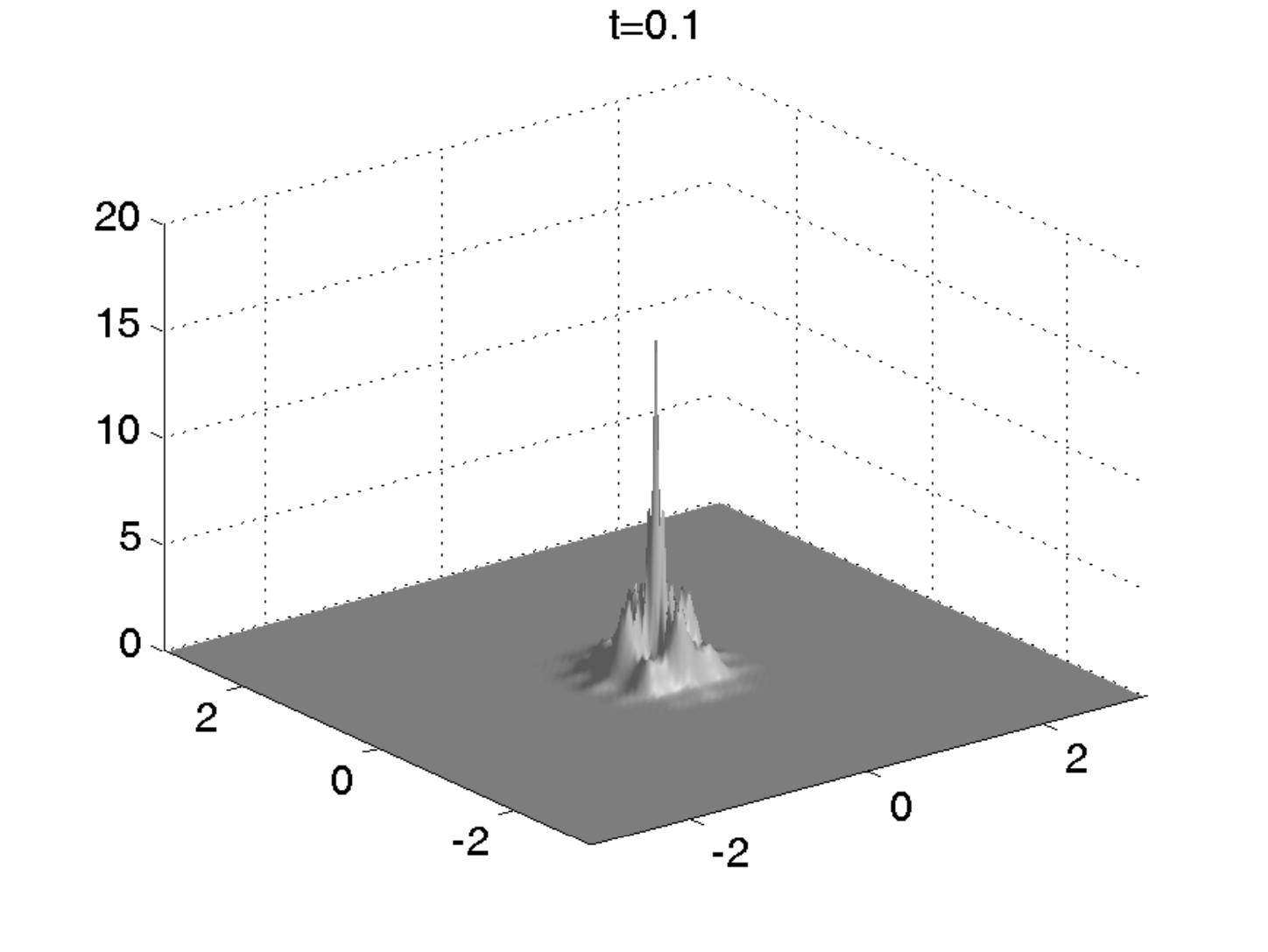}
\includegraphics[width=50mm]{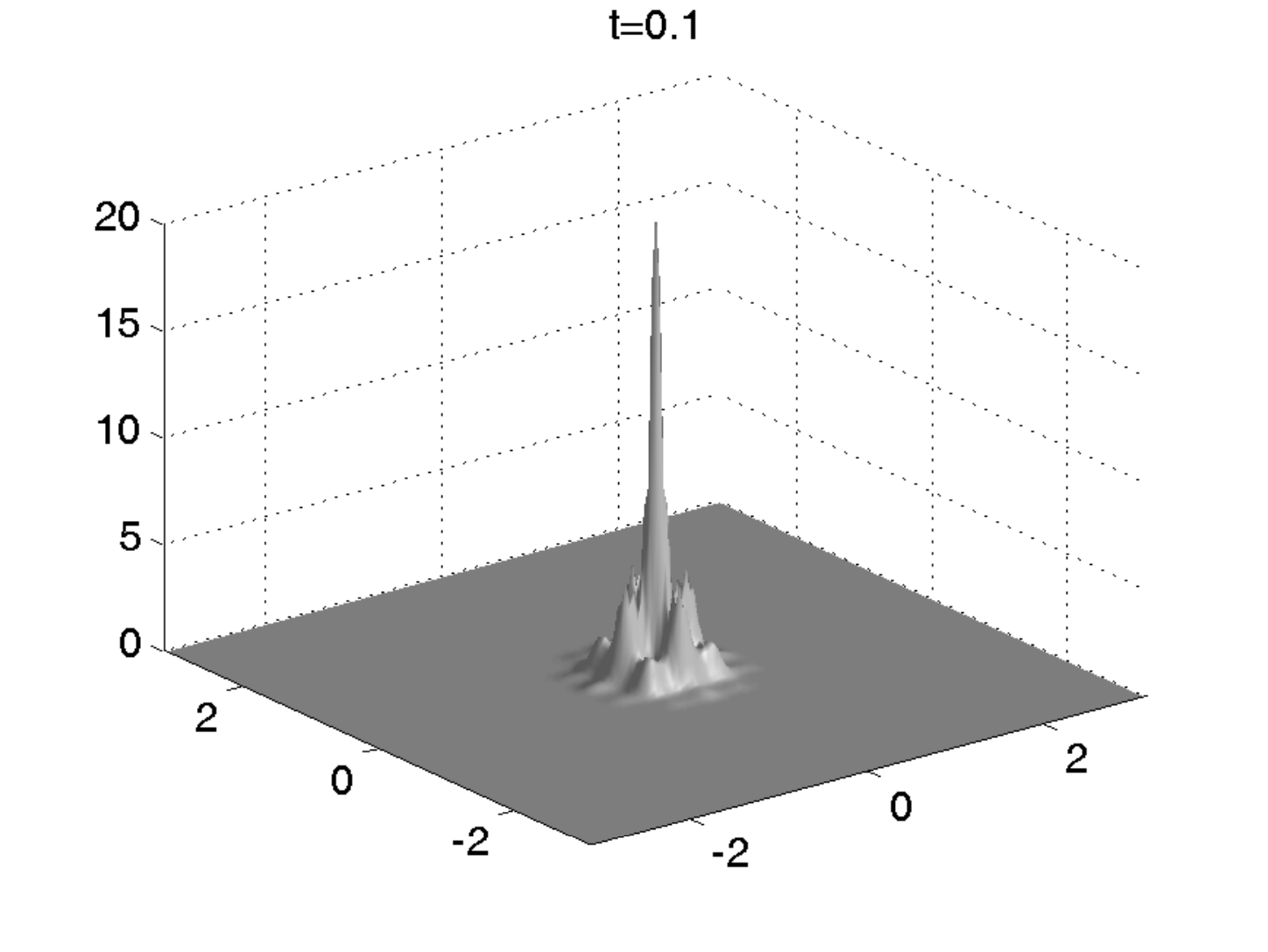}
\caption{Graphs of $|u(t,x,y,0)|^2$ of the example \ref{ex:pl} for fast-potential model by the typical pseudo-spectral method given in section \ref{sec:2.1} with $\og=160$,
left: mesh size $h=\frac{1}{128}$, right: $h=\frac{1}{1024}$.}
\label{Fig013}
\end{figure}

The results are shown in Figures \ref{Fig012}-\ref{Fig013}. We can see the interaction between the lattice potential and the wave packet.
Certainly, the wave function of the fast model is more peaked than the solution of the time-averaged model. From these figures, we can also
find that the Bloch-decomposition based algorithm is more efficient than traditional pseudo-spectral method in this case.

\section{Stability estimate}\label{anal:sect}
In this Section we present the analytical results needed in order to rigorously justify the approximation used in Section 2. Because of the singularities of the Coulomb potential, we use a mollified version in the numerical simulation. Here we discuss the additional errors introduced by the smoothing.

In order to avoid technicalities we only state a Proposition which shows that equation \eqref{eq:fast} is stable under perturbations of the potential. The proof is indeed quite standard, as similar results are well known in the literature, see below for the references.
The proof of stability of equation \eqref{eq:fast} with respect to perturbations of the potential is done by adapting in a straightforward way the local well-posedness result. For furhter details we refer the reader to \cite{Kato87}, \cite{Caz} Section 4.
\begin{theorem}
Let $V, \tilde V$ be two potentials such that $V, \tilde V\in L^\infty(\R; L^{\frac{3}{1+3\delta}}(\R^3)+L^{\frac{3}{1-3\delta}}(\R^3))$,
$\nabla V, \nabla\tilde V\in L^\infty(\R; L^{\frac{3}{2+3\delta}}(\R^3)+L^{\frac{3}{2-3\delta}}(\R^3))$, for some $\delta>0$ small.
Let $u, \tilde u$ be the solutions of the Cauchy problems
\begin{equation*}
\left\{\begin{aligned}
&i\partial_tu=-\Delta u+V u+C_1(|\cdot|^{-1}\ast|u|^2)u-a|u|^\sigma u\\
&u(t,\xb)\big|_{t=0}=u_0(\xb),
\end{aligned}\right.
\end{equation*}
  with the potential $V$ and respectively with $V$ replaced by $\tilde V$, and initial data $u_0$ and, respectively $\tilde u_0\in H^1(\R^3)$, respectively. Let $0<T<\min\{T_{max}, \tilde T_{max}\}$, where $T_{max}, \tilde T_{max}$ are the maximal existence times for $u, \tilde u$, respectively. Then, for any admissible pair $(q, r)$, we have
\begin{equation}\label{eq:stab1}
\|u-\tilde u\|_{L^q((0, T); W^{1, r}(\R^3))}\leq C\pare{
\|u_0-\tilde u_0\|_{H^1}+\|V-\tilde V\|_{L^\infty_t(L^{\frac{3}{1+3\delta}}_x+L^{\frac{3}{1-3\delta}}_x)}
+\|\nabla(V-\tilde V)\|_{L^\infty_t(L^{\frac{3}{2+3\delta}}_x+L^{\frac{3}{2-3\delta}}_x)}}.
\end{equation}
\end{theorem}

\section{Conclusion}
In this paper, we propose the time-splitting (Bloch-decomposition based) pseudo-spectral
methods to simulate the XFEL Schr\"odinger equation (with and without periodic potential). Our simulation results go far beyond known analytical results \cite{xfel1}. In particular we demonstrate numerically that the time-averaging procedure works in mass supercritical/energy subcritical cases, even up to blow-up time. Moreover, we show numerically that the energy of the oscillatory problem converges weakly in time to the energy of the time averaged problem (which is constant when the vector $e(t)$ is not dependent on time). Also, we demonstrate the impact of external trapping potentials and periodic lattice potentials on the electron beam wave function and numerically verify the time average procedure in those cases.

\end{document}